\documentclass[a4paper,11pt]{article}  
\usepackage{multirow}         
\usepackage{mathtools}
\usepackage{subcaption}
\usepackage{amsfonts,amssymb,amsmath,amsthm,latexsym,bbm}   
\usepackage{enumerate} 
\usepackage{epsf,epsfig} 
\usepackage{xcolor,colortbl,color}
\usepackage{graphicx,graphics}  
\usepackage{caption}  
\usepackage[utf8]{inputenc}
\usepackage{tikz}
\usepackage[utf8]{inputenc}
\usepackage[english]{babel}
\usepackage{soul}
\usepackage{comment}
\usepackage{hyperref} 
\makeatletter \@addtoreset{equation}{section} \makeatother
\makeatletter \@addtoreset{enunciato}{section} \makeatother
\newcounter{enunciato}[section]

\newtheorem{ittheorem}{Theorem}
\newtheorem{itlemma}{Lemma}
\newtheorem{itproposition}{Proposition}
\newtheorem{itdefinition}{Definition}
\newtheorem{itremark}{Remark}
\newtheorem{itclaim}{Claim}
\newtheorem{itfact}{Fact}

\newtheorem{itconjecture}{Conjecture}
\newtheorem{itcorollary}{Corollary}
\newtheorem{example}{Example}

\newenvironment{theorem}{\addtocounter{enunciato}{1}
\begin{ittheorem}}{\end{ittheorem}}
\newenvironment{lemma}{\addtocounter{enunciato}{1}
\begin{itlemma}}{\end{itlemma}}
\newenvironment{proposition}{\addtocounter{enunciato}{1}
\begin{itproposition}}{\end{itproposition}}
\newenvironment{definition}{\addtocounter{enunciato}{1}
\begin{itdefinition}}{\end{itdefinition}}
\newenvironment{remark}{\addtocounter{enunciato}{1}
\begin{itremark}}{\end{itremark}}

\newenvironment{corollary}{\addtocounter{enunciato}{1}
\begin{itcorollary}}{\end{itcorollary}}

\newcommand{\rN}{\mc{N}}
\newcommand{\rW}{\mc{W}}
\newcommand{\rS}{\mc{S}}
\newcommand{\vs}{v_\mu}
\newcommand{\varz}{{a_\alpha}} 
\newcommand{\err}{{\mathfrak{E}}}
\newcommand{\rWl}{{\mc{W}_\lambda}}
\newcommand{\limdist}{\Longrightarrow}
\newcommand{\sZ}{\sigma_0} 
\newcommand{\Lconv}{\Lambda_{\text{conv},s}} 
\newcommand{\Law}{\overset{\text{Law}}{=}}	
\newcommand{\hlg}{\hat{\lambda}_{g}}

\newcommand{\gep}{\varepsilon}

\newcommand{\mc}[1]{{\mathcal #1}}

\newcommand{\bb}[1]{{\mathbb #1}}
\newcommand{\bbm}[1]{{\mathbbm #1}}
\newcommand{\R}{\ensuremath{\mathbb{R}}}

\newcommand{\Z}{\ensuremath{\mathbb{Z}}}
\newcommand{\N}{\ensuremath{\mathbb{N}}}
\newcommand{\PP}{\ensuremath{\mathbb{P}}}

\newcommand{\EE}{\ensuremath{\mathbb{E}}}
\newcommand{\dd}{\ensuremath{\mathrm{d}}}

\newcommand{\prt}[1]{\left( #1 \right)} 
\newcommand{\crt}[1]{\left[ #1 \right]} 

\newcommand{\bPP}{\ensuremath{\overline{\mathbb{P}}}}
\newcommand{\bEE}{\ensuremath{\overline{\mathbb{E}}}}

\newcommand{\bigo}{\mathcal{O}}

\DeclareMathOperator{\Var}{Var}
\newcommand{\fl}[1]{\lfloor #1 \rfloor}  
\newcommand{\ind}[1]{ \mathbf{1}_{ \{ #1 \} } } 
\newcommand{\tc}{K_0} 

\begin{document}  

\title{Gaussian, stable, tempered stable and mixed
limit laws for 
  random walks in cooling random environments.
}

\author{\renewcommand{\thefootnote}{\arabic{footnote}}
Luca Avena\footnotemark[1]\,,
Conrado da Costa\footnotemark[2]\,,
Jonathon Peterson\footnotemark[3]}

\date{\today}

\footnotetext[1]{Mathematical Institute,
  Leiden University,
P.O.\ Box 9512, 2300 RA Leiden, The Netherlands}
\footnotetext[2]{Department of Mathematical Sciences,
  Durham University,
South Road, Durham DH1 3LE, UK}
\footnotetext[3]{Purdue University,
  Department of Mathematics, 
150 N. University Street, West Lafayette, IN 47907}

\maketitle


\begin{abstract}Random Walks in Cooling Random Environments (RWCRE)
is a model of random walks in dynamic random environments
where the entire environment is resampled along a fixed sequence of times,
called the “cooling sequence”,
and is kept fixed in between those times.
This model interpolates between that of a homogenous random walk,
where the environment is reset at every step,
and Random Walks in (static) Random Environments (RWRE),
where the environment is never resampled.
In this work we focus on the limiting distributions of
one-dimensional RWCRE in the regime where the fluctuations of the corresponding (static) RWRE
is given by a $s$-stable random variable with $s \in (1,2)$.
In this regime, due to the two extreme cases
(resampling every step and never resampling, respectively),
a crossover
from Gaussian to stable limits for sufficiently regular 
cooling sequence was previously conjectured. 
Our first result  
answers
affirmatively this conjecture by making clear  
critical exponent, norming sequences and limiting laws associated with
the crossover which 
demonstrates 
a change from Gaussian to $s$-stable limits,
passing at criticality through a certain generalized tempered stable distribution.
We then explore the resulting RWCRE scaling limits for general cooling sequences.
On the one hand,
we offer sets of operative sufficient conditions
that guarantee asymptotic emergence of either Gaussian,
 $s$-stable or generalized tempered distributions from a certain class.
On the other hand,
 we give explicit examples and describe how to construct irregular cooling
sequences for which the corresponding limit law is characterized
by mixtures of the three above mentioned laws.
To obtain these results, we need and derive a number
of refined asymptotic results for
the static RWRE with $s \in (1,2)$
which may be of independent interest.

\\
{\it Acknowledgment:} 
L.A. was supported through NWO Gravitation Grant NETWORKS-024.002.003.
J.P. was supported through a Simons Foundation Collaboration Grant \#635064.  
We are grateful to Gennady Samorodnitsky for his help regarding the statements and proofs of the results in Appendix \ref{A}. 

\end{abstract}

\newpage

\tableofcontents

\newpage

\section{Context and overview.} \label{Intro}

\paragraph{Perturbation of a frozen media through resetting.}
RWRE (Random Walks in Random Environments) is a well-known model,
of central relevance within the theory of disordered systems,
for particles moving in media with impurities.
It consists of a Markov chain with random transition kernels
determined by an underlying field of variables,
referred to as random environment,
which is sampled at time zero from a given law and stays
``frozen'' during the evolution of the Random Walk (RW).

Rigorous studies on RWRE can be traced back to the 1970s~\cite{S75}
and along the years the model has been widely investigated
on $d$-dimensional integer lattices.
This setup poses many challenges and still several questions
remain open when $d\geq 2$, see~\cite{Z04}.
Unlike the higher dimensional setup, for $d=1$,
RWRE is reversible and,
by the analysis of the associated hitting times and the so-called potential,
a fairly complete picture of its limiting properties has been obtained along the years.

Depending on the choice of the law of the environment,
strong spatial local effects
lead to substantial qualitative differences with respect
to a standard homogeneous RW.
Indeed, due to the spatial inhomogeneities,
trapping and slow-down phenomena can give rise to
a variety of rich behaviors such as
suballistic transience~\cite{S75}, non-Gaussian limiting distributions~\cite{KKS75,Sin82},
sub-exponential large deviation probabilities~\cite{CGZ00,GdH94},
aging~\cite{ESZ09,Z04}, etc. 

In the recent~\cite{AdH17}, the authors introduce a model,
referred to as RWCRE (Random Walks in Cooling Random Environments),
which can be thought of as a perturbation of RWRE,
obtained by resampling the environment in an independent fashion
over a prescribed sequence of times.
This sequence is described by a function,
which is referred to as {\em cooling map}.
RWCRE is thus an example of \emph{RW in dynamic random environments}
in which depending on the choice of the cooling map,
one can flexibly ``tune''
the intensity of space-time correlations. In particular, RWCRE interpolates between RWRE, corresponding to 
no resetting, and a homogeneous RW after averaging over the random field, which corresponds to 
resetting the environment at every time unit.

The overall goal is to see what sort of limiting behavior
can emerge for different cooling maps.
In this respect, it is worth mentioning that
the study of RW models in dynamic random environments has witnessed various
interesting progresses over the past decade. Yet, unlike in the RWCRE, in most of
this literature where the medium changes over time,
limiting results are obtained for models with good mixing properties leading to behaviors such as the ones observed for a standard homogeneous RW, see e.g.~\cite{ABF17,HKT20} and references therein for a recent account.
For RWCRE
one can not only recover some of the non-Gaussian limiting distributions of the RWRE model, but also obtain some new limiting distributions that were not obtained either in the homogeneous RW or RWRE models \cite{ACdCdH20}.

\paragraph{State of the art of RWCRE.} 
The study of RWCRE in one-dimension has been pursued
in a sequence of recent works~\cite{AdH17,ACdCdH19,ACdCdH20,AdC20,Yon19}
in various regimes which we next briefly describe. 

A general recurrence criterion is still open,
although for diverging cooling sequences, as shown in~\cite{ACdCdH20},
it can be related to the classical (non-local) recurrence criterion in~\cite{S75} for RWRE.

For the law of large numbers for the RW displacement,
Thm 1.10 in~\cite{ACdCdH19} and the various general statements in~\cite{AdC20},
show that the limit speed is deterministic and can be characterized
in full generality~\cite[Section 3]{AdC20}.
In particular, its value coincides with the RWRE speed for cooling maps that diverge in a Ces\`aro sense.

For large deviations of the empirical speed,
if increments between consecutive resettings diverge,
it is shown in~\cite{ACdCdH19} that the (quenched)
asymptotic costs for deviations are exactly as in RWRE,
regardless of the speed of divergence of the resettings.
Which is to say, somewhat surprisingly in light of the fluctuation results,
that large deviation rate functions for the empirical
speed of RWRE are left unchanged under a wide class of
perturbation induced by the cooling map.

When we consider fluctuations and scaling limits,
the picture is much more delicate and heavily depends
on the law of the corresponding RWRE that one is perturbing.
Let us briefly recall that for transient RWRE, there is a certain parameter
$s>0$ associated to the law of the environment
(see~\eqref{parameters} below), which captures
essentially four different classes of possible scaling limits:
\begin{enumerate}
 \item \textbf{Recurrent:} Non-Gaussian limiting distribution with strongly sub-diffusive scaling $(\log n)^2$.  Limit distribution is a non-trivial functional of Brownian motion \cite{Sin82,Kes86}. 
 \item \textbf{Transient, $s\in (0,1)$:} Limiting distribution with no centering and sub-linear scaling $n^s$. Limit distribution is a transformation of an $s$-stable law \cite{KKS75}. 
 \item \textbf{Transient, $s\in [1,2)$:} Limiting distributions are $s$-stable with superdiffusive scaling $n^{1/s}$; linear centering when $s>1$ \cite{KKS75}.
 \item \textbf{Transient, $s \geq 2$:} Gaussian limiting distribution; diffusive scaling when $s>2$ \cite{KKS75}. 
\end{enumerate}

So far, fluctuations results for RWCRE have been obtained  only when the underlying environment is in classes 1.\ and 4.\ (for $s>2$). 

\paragraph{Cooling in the Sinai regime.} 
A recurrent RWRE is sometimes referred to as a Sinai walk due
to Sinai's derivation of the limiting distribution for this case \cite{Sin82},
and the corresponding limiting distribution is called the Sinai-Kesten
law due to Kesten's derivation of the density \cite{Kes86}.
As shown first in~\cite{AdH17} for some regular cooling maps,
and then in great generality in~\cite{ACdCdH20},
convergence in distribution depends on the regularity and speed
of the chosen cooling map.
Sub-sequential limits can be characterized in general and may lead to mixtures
of Gaussian and 
Sinai-Kesten laws~\cite[Thm.\ 2]{ACdCdH20}. 
Fluctuations are controlled by the total variance
of the RWCRE and are essentially always sub-diffusive,
and the limiting distribution is Gaussian
only for cooling maps in which increments between resetting do not grow more than exponentially,
see~\cite[Cor.\ 1]{ACdCdH20}.

In this regime, the recent~\cite{Yon19} investigates convergence
of the full RWCRE process for polynomially and exponentially growing
cooling increments leading, respectively, to a time-scaled
Brownian motion~\cite[Thm. 1]{Yon19} and to a (degenerate)
random constant distributed as a standard Gaussian~\cite[Thm. 2]{Yon19}.

\paragraph{Cooling in the CLT regime, $s> 2$.} 
The other well-understood and actually the easiest regime,
corresponds to  $s>2$ in the class 4.\ mentioned above.
In this case, 
the limits are Gaussian for any cooling map and if the increments between resettings diverge then the scaling is of the form $c\sqrt{n}$~\cite[Thm.\ 3, Cor.\ 2]{ACdCdH20}.

\paragraph{New results at glance: fluctuations in the stable regime $s\in(1,2)$} \label{fluctintro}
What happens when $s\in(0,2)$
(i.e. classes 2. and 3. above) is still open, and in this paper we investigate
fluctuations for $s\in(1,2)$.
In this case, refered to as the stable regime,
RWRE is transient and ballistic,
with stable limit laws after scaling by $n^{1/s}$.
Here, unlike 
the Sinai regime or the CLT regime,
the variance no longer determines the scaling
for the RWCRE and hence the analysis depends on
the regularity assumptions on the cooling map.
Our main results are summarized as follows.

The first result, Theorem~\ref{PolyTransition}
describes the scaling limit of RWCRE when the cooling map
has {\em polynomially} growing increments.
In this case, the system presents three possible limiting scenarios.
The critical scenario occurs when the exponent of polynomial growth equals $1/(s-1)$,
and in this case the fluctuations are of the order $n^{1/s}$ and the limiting
distribution is neither Gaussian nor stable but instead a type of distribution
which we call generalized tempered $s$-stable.
In the supercritical regime, fluctuations stay of order $n^{1/s}$
but the limiting law coincides with the stable one for static RWRE.
On the other hand, in the subcritical case
the limiting distribution is Gaussian and the fluctuations 
are scaled by $n^\beta$, where $\beta \in (1/2,1/s)$
depends explicitly on the exponent of the polynomial growth of the cooling increments.
Such a crossover from Gaussian to stable limits for polynomial cooling
was conjectured in~\cite{AdH17} on the basis of the
fluctuations of the RWRE hitting times.
Thus, Theorem~\ref{PolyTransition} not only settles affirmatively
this conjecture but also identifies the precise order of fluctuations,
the critical exponent where this crossover occurs, and the limiting distribution at criticality.

We then explore limit distributions
for general but sufficiently regular cooling sequences.
In particular, we give  operative conditions on the cooling sequence
to obtain Gaussian (Theorem~\ref{GeneralConditions-N}), stable 
or generalized tempered stable (Theorem~\ref{GeneralConditions-S})
limit distributions. 
Statements in Theorem~\ref{PolyTransition} for polynomial cooling maps
are in fact special case of these general theorems.
Then by constructing an ``interweaving'' of multiple polynomial cooling maps
we show in Theorem~\ref{arbmix} that one can construct cooling maps
for which the limiting distribution can be an arbitrary linear combination
of Gaussian, $s$-stable, and a member of a broad class of
generalized tempered stable random variables.


The proofs of these results make use of a variety of
standard techniques such as  Lindeberg conditions and the
characterization of stable laws via Poisson point processes.
Yet their implementation is non-trivial and
require a number of precise estimates for RWRE which we derive specifically for the proofs.
Among these technical RWRE estimates
are some moment asymptotics and bounds for RWRE in the stable
regime that may be of independent interest, see Theorem~\ref{lem:VarZn}.
 
\paragraph{Structure of the paper.}
The next section is devoted to model definitions,
notation and basic results.
In particular, the classical RWRE is introduced in
Section~\ref{setupRWRE} along with the main assumptions and asymptotic
results in the stable regime which represents
our point of departure.
RWCRE is then defined in Section~\ref{model}.
Our main results mentioned above are collected and discussed in Section~\ref{S3}.

We then start all the proofs.
Those about RWRE are presented in Section~\ref{S4}
together with a number of other large deviation estimates
which will be instrumental for the analysis of the cooling model.
Proofs of the RWCRE scaling limits for regular cooling maps
are given in Sections~\ref{S5} and~\ref{S7},
in which we prove, respectively, emergence of Gaussian, and of
stable or generalized tempered stable distributions. 
In Section~\ref{S8} we treat the non-regular
maps which lead to mixture of different limiting laws,
and in Section~\ref{rEx} we give examples of some highly irregular
cooling maps which demonstrate both how the techniques of this paper
can be extended to obtain subsequential limiting distributions
not contained in our main results and also how the techniques
of this paper can be applied to obtain limiting distributions
even when the cooling maps do not satisfy the regularity conditions
of our main results.

We conclude with three appendices: Appendix~\ref{A} devoted to facts about stable laws, Appendix~\ref{Regeneration}, that recalls the construction of regeneration times
for RWRE and some related results, and Appendix~\ref{Acp} which contains simple technical lemmas used in the proofs.

\section{Setting and Background}\label{S2}
\subsection{RWRE: stable regime \texorpdfstring{$s \in (1,2)$}{}}\label{setupRWRE}
Throughout the paper we use the notation $\N_0 = \N \cup \{0\}$ with $\N = \{1,2,\dots\}$.
The classical one-dimensional (static) RWRE model is defined as follows.
Let $\omega=(\omega_x)_{x\in\Z}$ be an i.i.d.\ sequence with law
\begin{equation}\label{iid}
\mu = \alpha^{\Z},
\end{equation} 
for some probability distribution $\alpha$ on $(0,1)$.
We write $\langle\cdot\rangle$ to denote the expectation w.r.t.\ $\alpha$.
\begin{definition}[{\bf RWRE}]\label{RWREdef}
  \text{}\\ 
{\rm Let  $\omega$ be an (i.i.d.) environment sampled from $\mu$.
We call \emph{Random Walk in Random Environment} the Markov chain
$Z = (Z_n)_{n\in\N_0}$ with state space $\Z$ and transition probabilities}
\begin{equation} 
\label{Ker}
P^{\omega}(Z_{n+1} = x + e \mid Z_n = x) 
= \left\{
\begin{array}{ll}
\omega_x &\mbox{ if } e = 1,\\  
1 - \omega_x &\mbox{ if } e = - 1,
\end{array}
\right. 
\qquad n \in \N_0.
\end{equation}
{\rm We denote by} $P_x^{\omega}(\cdot)$
{\rm the \emph{quenched} law of the Markov chain identified by
the transitions in~\eqref{Ker} starting from $x\in\Z$, and by
\begin{equation}\label{annealed}
P_{x}^{\mu}(\cdot) = \int_{(0,1)^\Z} P_x^{\omega}(\cdot)\,\mu(\dd\omega),
\end{equation}
the corresponding \emph{annealed} law.}
\end{definition}
One-dimensional RWRE is by now well understood,
both under the quenched and the annealed law.
It exhibits very different limiting behaviors
(asymptotic speed, scaling limits and large deviations)
depending on the choice of $\mu$
(or $\alpha$ in the present i.i.d. setting captured by~\eqref{iid}).
For a general overview, we refer the reader to the
lecture notes by Zeitouni~\cite{Z04}.
Here we collect some basic facts and definitions that
will be needed throughout the paper.
We will focus on the annealed stable regime
as introduced below and first studied by
Kesten, Kozlov and Spitzer~\cite{KKS75}.

Let us start with some assumptions on $\mu$ (or $\alpha$).
A crucial quantity to characterize the asymptotic properties
of RWRE is the ratio of the transition probabilities
to the left and to the right at the origin
(or any other vertex due to the i.i.d. assumption~\eqref{iid})
$\rho_0 = \frac{1 - \omega_0}{\omega_0}$.
For the remainder of the paper, we assume that
\begin{equation}\label{trans}
\langle \log \rho_0\rangle <0,
\end{equation}
which, as shown in~\cite{S75},
guarantees {\bf right transience}. 
In what follows we restrict ourselves to the regime where  
\begin{equation}\label{parameters}
  \exists\; s \in (1,2) \text{ such that }
  \quad \langle \rho_0^s\rangle = 1.
\end{equation}
This condition characterizes what we call \textbf{stable regime},
and as captured in the next proposition,
guarantees \emph{ballisticity} of the walk 
(see \eqref{speedann}). 
We further assume that 
\begin{equation}\label{lat}
\log \rho_0 \text{ is {\bf non lattice}},
\end{equation}
this is a technical assumption required in~\cite{KKS75}
to characterize emergence of limiting stable laws, see Eq.~\eqref{StableLeftSkew} below. Finally we require the following  {\bf ellipticity condition}
\begin{equation}\label{ellipticity}\langle \rho_0^{s+\epsilon} \rangle < \infty, \quad\text{for some }\epsilon>0,
 \end{equation}
 which is needed for the tail estimate in \eqref{precisetail}. We will consider $\mu$'s that satisfy all the above conditions,
which we summarize in the following definition.

\begin{definition}[{\bf $s$-canonical $\mu$ for the stable regime}] \rm{}
  \label{goodenvironment}
\text{}\\ 
We say that $\mu$  is \emph{$s$-canonical} if it satisfies conditions~\eqref{iid}, ~\eqref{trans},~\eqref{parameters}~\eqref{lat} and~\eqref{ellipticity}. 
\end{definition} 

The next proposition represents our point of departure.
In this statement, and in the sequel, we denote convergence in distribution
of an arbitrary sequence of random variables $(Y_n)_{ n \in \bb{N}}$
to a random variable $Y_*$ as $n\to\infty$ by $Y_n \limdist Y_*$.

\begin{proposition}[{\bf RWRE: speed, limit law and deviations for $s\in(1,2)$}]\label{RWREscaling} 
\text{}\\ \rm{}
Let $\mu$ be any $s$-canonical law with $s \in (1,2)$ and consider the RWRE
process $Z$ with environment sampled from $\mu$.
Then:
\begin{itemize}
\item {\bf (LLN)} under the annealed law (and under the quenched too),
  $Z$ is almost surely right-transient and admits deterministic limiting speed:
\begin{equation}\label{speedann}
P_{0}^{\mu} \left(\lim_{n\to\infty} \frac{Z_n}{n}=v_\mu \right)=1,  \quad \text{with } v_\mu = \frac{ 1 - \langle \rho_0 \rangle }{1 + \langle \rho_0 \rangle}>0.
\end{equation}
\item {\bf (Fluctuations)} under the annealed law $P_{0}^{\mu}$, there exists $b>0$,
such that 
\begin{equation}\label{StableLeftSkew}
 \frac{Z_n-v_\mu n}{ n^{1/s}} \limdist \rS_s, \quad \text{as } n\to\infty,
\end{equation}
where $\rS_s$ is the stable (mean zero totally skewed to the left)
random variable with characteristic function
\begin{equation}\label{stablecf}
E\left[ e^{i u \rS_s} \right]
=\exp\left[-b |u|^s \left(1+i\frac{u}{|u|}\tan\left(\frac{s\pi}{2}\right)\right)\right], \qquad u \in \R.
\end{equation}


\item {\bf (Moderate slow-down deviations)} there exists a constant $\tc>0$ such that 
\begin{equation}\label{precisetail}
\lim_{n\to\infty} \sup_{ t \in \tilde{I}_n} \left| \frac{P_{0}^{\mu}(Z_n - n\vs < -t )}{(n\vs-t)t^{-s}} - \tc \right| = 0, 
\end{equation}
where  $\tilde{I}_n : = [n^{1/s}(\log n)^3, n\vs - \log n ]$.
\end{itemize} 
\end{proposition}

The LLN  in~\eqref{speedann} was proved in~\cite{S75} and in particular it does not need assumptions ~\eqref{parameters},~\eqref{lat} and~\eqref{ellipticity}  in Def \ref{goodenvironment}.
The stable law convergence in~\eqref{StableLeftSkew} was proved in~\cite{KKS75} under a slight weaker assumption than the one in~\eqref{ellipticity}.
The latter is in fact only needed to show the limit in~\eqref{precisetail} which was proved in~\cite{BD18}.
We remark that the constants $b$ in~\eqref{stablecf} and $\tc$ in~\eqref{precisetail} are related by\footnote{It follows from the proof of \eqref{precisetail} in \cite{BD18} that $\tc$ can be expressed in terms of another constant $C_3$ which appears in a tail asymptotic result in \cite[Lemma 3.2]{BD18}. This same tail asymptotic result is also given in \cite[Lemma 6]{KKS75}, and it follows from this that one can derive a formula for $b$ in \eqref{stablecf} in terms of $C_3$ also. 
}
\begin{equation}\label{relate}
b = \tc \vs \Gamma(1-s) \cos(\tfrac{\pi s}{2}),
\end{equation}
where 
$\Gamma(1-s) = \frac{\Gamma(2-s)}{1-s} =  \frac{1}{1-s}\int_0^\infty e^{-t}t^{1-s}\, dt$. 

\subsection{RWCRE: Cooling} \label{model}

The cooling random environment is the \emph{space-time} random environment
built by partitioning $\N_0$ into a sequence of intervals,
and assigning independently to each interval an environment sampled from $\mu$. 
Formally, let $(T_k)_{k\in\N}$ be an increment sequence such that $T_k\in\N$,
we will refer to this sequence as {\bf cooling increment sequence}. 
We denote further by $\tau(k):=\sum_{i=1}^{k}T_i $ the $k$-th cooling time,
i.e. the time at which a new environment is freshly sampled  from $\mu$.
We will refer to $\tau$ as the {\bf cooling map}.

\begin{definition}[{\bf Random Walks in Cooling Random Environments (RWCRE)}]\label{RWCREdef}
\text{}\\ 
Consider a cooling increment sequence $(T_k)_{k\in\N}$  and a distribution $\mu$ on environments.
For a fixed $n\in\N$
set 
\begin{equation}\label{defellT}
\ell_n+1: = \inf\{\ell \colon \tau(\ell) > n\},
\qquad \text{ and }\qquad \bar{T}_n : = n - \sum_{k = 1}^{\ell_n} T_k.
\end{equation}
Let $\bar{\omega} = \{\omega^{(k)}\}_{k\geq 1} = \{ (\omega_x^{(k)})_{x \in\Z} \}_{k\geq 1}$
be an i.i.d.\ sequence of environments with $\omega^{(k)} \sim \mu$.
We define the RWCRE sequence $X=(X_n)_{n\in\N_0}$
in the sequence of environments $\bar{\omega}$ and with cooling map $\tau$ by
\begin{equation}\label{sumpiece}
X_n : = \sum_{k = 1}^{\ell_
n} Z^{(k)}_{T_k} + Z^{(\ell_n + 1)}_{\bar{T}_n} \qquad n\in\N_0, 
\end{equation}
where for $k\geq 1$, $Z^{(k)}_\cdot:=\left(Z^{(k)}_n\right)_{n \in \N_0}$
is distributed as a RWRE process with underlying environment 
$\omega^{(k)} = (\omega_x^{(k)})_{x \in\Z}$,
and the sequence of random walks $\{Z^{(k)}_\cdot\}_{k\geq 1}$ are independent. 
\end{definition}
This process corresponds to a discrete-time RW evolving
in a random environment with law $\mu$ which is resampled
in an independent fashion along the sequence of times $\tau(k)$
determined by the cooling map.
We notice that for $T_1= \infty$ this model reduces to RWRE,
while for $T_k\equiv 1$ it reduces to 
a homogeneous RW (under the annealed measure) with local drift $\EE[Z_1]$. 
The name cooling comes from the fact that 
when dealing with maps for which the increments
$T_k$ eventually diverges, the environment will be resampled
less and less, and hence, depending on the growth of $T_k$,
the corresponding motion will resemble the random walk in
the static or ``frozen" random environment.
Notice that as defined in \eqref{defellT},  $\ell_n+1$ denotes the index of the
increment in which $n$ belongs to, and that $\ell_n$ counts
the total number of resettings up until time $n$.

We will analyze the model under the \emph{annealed law} that starts from the origin. 
Formally this refers to the path measure obtained by the average
with respect to $\mu^\N$ of the quenched path measure, say $P_0^{\bar{\omega},\tau}$,
associated to the kernel $P^{\bar{\omega},\tau}(X_{n+1} = x + 1 \mid X_n = x)$. 
In what follows, to lighten the notation, we will simply denote by
\begin{equation}\label{annlaw}
\bb{P}(\cdot) :=  \mu^{\bb{N}}\otimes  P_0^{\bar{\omega},\tau}(\cdot) = \int P_0^{\bar{\omega},\tau}(\cdot) \, d\mu^{\bb{N}}(\bar{\omega}),
\end{equation}
such an {\bf annealed measure} and, in each statement,
we will specify that we consider RWCRE associated to a given $s$-\emph{canonical law}
$\mu$ and a given \emph{cooling map} $\tau$. 
As a slight abuse of notation when discussing just a single RWRE process $(Z_n)_{n\geq 0}$
we will also use $\PP$ for the annealed measure instead of $P^\mu$
since as noted above a RWRE can be seen as a RWCRE with $T_1 = \infty$.

\subsection{Relevant Distributions}
In the next section we state our results about
the limit behavior of RWCRE for perturbations of RWRE in the stable regime.
As we will see, depending on the choice of the cooling map,
we will encounter the following type of limit laws,
or possibly mixtures of them.
These correspond to:
\begin{itemize}\item the \emph{Standard  Gaussian}, denoted by $\cal N$;
\item the (\emph{mean-zero totally skewed to the left}) {\em Stable} defined
by its characteristic function in Eq.~\eqref{stablecf}, and denoted by $\rS_s$;
\item a third special type of random variable $\rW_\lambda$
defined below and referred to as (\emph{mean-zero totally skewed to the left})
\emph{generalized tempered $s$-stable laws}.
\end{itemize}
\begin{definition}[{\bf Generalized mean-zero left-skewed tempered $s$-stable laws}] \label{gtslaw}
\text{}\\ 
For $s \in (1,2)$, a random variable $\rW_\lambda$ is said to have
a generalized mean-zero left-skewed tempered $s$-stable law,
if it has characteristic function 
\begin{equation}\label{Wlambdacritical}
E[ e^{iu\rWl} ] = \exp\left\{ \int_{-\infty}^0 (e^{iux}-1-iux) \lambda(x) \, dx \right\}, 
\end{equation}
where the function $\lambda(x)$ is of the form $\lambda(x) = c |x|^{-s-1} a(x)$
for some $c>0$ and some non-decreasing, continuous function $a$ on $(-\infty,0]$
with $a(0) = 1$ and $\lim_{x\to-\infty} a(x) = 0$. 
\end{definition}

\begin{remark}[\textbf{Relation of the three laws}]\label{3laws}\rm{}
We note that the family of random variables $\rW_\lambda$ interpolates between the standard Gaussian $\cal N$ and the $s$-stable $\rS_s$, in the sense that they belong 
to the closure, with respect to weak convergence, of
the vector space of generalized tempered $s$-stable laws.
Indeed, for any $c,r>0$ let 
\begin{equation}\label{lcr}
 \lambda_{c,r}(x) =  c|x|^{-s-1}(1+x/r)_+, \qquad x<0.  
\end{equation}
Then $\rW_{\lambda_{c,r}}$ converges in distribution as $r\to \infty$ to the random variable $\rS_s$ with characteristic function as in \eqref{stablecf} with $b=-c\Gamma(-s)\cos(\frac{\pi s}{2})$.
On the other hand, as $r\to 0^+$ the random variables $r^{\frac{s}{2}-1} \rW_{\lambda_{c,r}}$ converge in distribution to a centered Gaussian.  One can check these claims by taking limits of the corresponding characteristic functions. 
\end{remark}

We use the term \emph{generalized} in Definition~\ref{gtslaw}
because left-skewed tempered $s$-stable laws are the special case when
$\lambda(x) = c |x|^{-s-1} e^{\theta x}$ for some $\theta>0$.
Tempered stable distributions and the corresponding
L\'evy processes
(also called L\'evy flights~\cite{Ism95} and the CGMY model~\cite{CGMY03})
have been the subject of interest recently in financial modeling~\cite{CT04,FL21}, 
but we are not aware of any other results where tempered stable laws arise naturally
as the limit of a discrete process as is the case with our results below.

\section{Results} 
\label{S3}
Our first theorem shows the mentioned crossover,
from normal to stable limit, passing through an intermediate
critical generalized tempered stable
law, as one changes the polynomial cooling rate in~\eqref{polygrow} below.

\begin{theorem}[{\bf Trichotomy: phase transition for polynomial cooling}]
  \label{PolyTransition}
\text{}\\ 
Let $X$ be  a RWCRE  associated to a given
$s$-canonical law $\mu$, as in
Def.~\ref{goodenvironment}, and consider a cooling
map $\tau$ with (eventual) polynomial growth, that is,
such that
\begin{equation}\label{polygrow}
  \lim_{k\to\infty} \frac{T_k}{A k^a} = 1,
  \qquad \text{ for some }
  A,a \in (0,\infty).
\end{equation}
Then, the following three limiting scenarios are possible:
\begin{itemize}
\item {\bf (Normal)}
  For $a <\frac{1}{s-1}$,
\begin{equation}\label{gauss_slow}
  \frac{X_n - \EE[X_n]}{Bn^{\beta} }
  \limdist \rN,
\end{equation}
where $\beta := \frac{a(3-s)+1}{2(a+1)}$ and 
$B^2  := \frac{2\tc \vs^{3-s}A^{\frac{2-s}{a+1}}(a+1)^{\frac{a(3-s)+1}{a+1}}}{(2-s)(3-s)(a(3-s)+1)}$.
\item {\bf (Critical)} 
  For $a = \frac{1}{s-1}$,
\begin{equation}\label{Xnstableliml}
  \frac{X_n - \EE[X_n]}{n^{1/s}}
  \limdist \rW_{\lambda_{c,r}}, 
\end{equation}
where $\lambda_{c,r}$ is defined as in \eqref{lcr} with $c=\tc \vs s$ and 
$r = \vs\left(\frac{s}{s-1} \right)^{1/s} A^{\frac{s-1}{s}} $.
\item {\bf(Stable)}
  For $a> \frac{1}{s-1}$,
  \begin{equation}\label{Xnstablelim}
    \frac{X_n - \EE[X_n]}{n^{1/s}}
    \limdist \rS_s.
\end{equation}
\end{itemize}  
\end{theorem}
The three statements in
Theorem~\ref{PolyTransition} are special cases of the
following two more general theorems which give
sufficient conditions for respectively, normal, and a class of generalized tempered stable laws as in Def.
\ref{gtslaw} which includes $\rS_s$ and $\rW_{\lambda_{c,r}} $ from Theorem ~\ref{PolyTransition}.
See Remark~\ref{3laws} for their relations. 
\begin{theorem}[{\bf Sufficient conditions for
      pure Gaussian limits}]
\label{GeneralConditions-N}
\text{}\\ 
Let $X$ be  a RWCRE associated to an $s$-canonical
law $\mu$, as in Def.~\ref{goodenvironment}. If the
cooling map $\tau$ is such that 
\begin{equation}\label{Gauss_iff}
  \lim_n\sup_{k \leq  n}
  \frac{T_k}{\prt{\sum_{k=0}^n (T_k)^{(3-s)} }^{1/2}}=0,
\end{equation}
then
\begin{equation}\label{gauss_slow2}
  \frac{X_n - \EE[X_n]}{\sqrt{\Var(X_n)}}
  \limdist \rN.
\end{equation}
\end{theorem}

The statement above is proven in Section \ref{S5} by checking the classical Lindeberg's conditions.
We notice in particular that the norming sequence in \eqref{gauss_slow2}
is determined by the standard deviation and in particular its asymptotic
behavior varies as the growth of the cooling increment sequence varies.
This variation of the scaling as a function of the cooling growth
can be appreciated in \eqref{gauss_slow}. On the other hand,
as stated in Theorem \ref{PolyTransition},
if the polynomial increments start to grow too much,
it is a signature of exiting the Gaussian world and in particular we see that
for the emergence of both non-Gaussian laws in \eqref{Xnstableliml} and \eqref{Xnstablelim},
the corresponding norming sequences are not a function of the power in the polynomial
cooling growth and are given by $n^{1/s}$ rather than the standard deviation.
In the next theorem we offer sufficient regularity conditions on the cooling
map which guarantee emergence of a sub-class of the generalized stable distributions
in Def. \eqref{gtslaw}, which in particular include the two limiting random variables in
\eqref{Xnstableliml} and \eqref{Xnstablelim}.
This regularity of the cooling map is expressed in terms of the existence of a limit,
see \eqref{cag}, which captures the asymptotic stability for the empirical
distribution of the increments that are \emph{large},
meaning that they have non-negligible contribution on the scale ($n^{1/s}$) of
the global running time to the power $1/s$.
 
\begin{theorem}[{\bf Sufficient conditions for
     generalized s-Stable limits}]
\label{GeneralConditions-S}
\text{}\\ 
Let $X$ be  a RWCRE associated to an $s$-canonical
law $\mu$, as in Def.~\ref{goodenvironment}. 
Assume that the following limit exists
\begin{equation}\label{cag}
  \lim_{n} \frac{ \sum_{k=1}^{n} T_k \ind{T_k < x \tau(n)^{1/s}} }{\tau(n)} = g(x), \quad \text{for all } x \in (0,\infty), 
\end{equation} 
with $g$ being a continuous function on $[0,\infty)$ with $g(0) = 0$ and $g(\infty) := \lim_{x\to\infty} g(x) \in [0,1]$.
\begin{itemize}

\item[{\bf(S1)}] If 
$\sup_n \sum_{k=1}^n \frac{T^{1/s}_k}{\tau(n)^{1/s}}< \infty$,
\quad and \quad 
$\lim_{n\to\infty} \sum_{k=1}^n \frac{T^{1/s}_k}{\tau(n)^{1/s}} \bbm{1}_{T_k<m}=0$ for all $m<\infty$, 
then
\begin{equation}\label{Xnstablelim1}
  \frac{X_n - \EE[X_n]}{n^{1/s}}  \limdist \rS_s. 
\end{equation}
\item[{\bf(S2)}] If
$\lim_{n\to\infty} \frac{\max_{k\leq n} T_k(\log T_k)^{4s}}{\tau(n)} = 0$,  then

\begin{equation}\label{Xnstablelim2}
\frac{X_{n} - \EE[X_{n}]}{n^{1/s} } \limdist 
\begin{cases}
\rS_s & \text{if } g(\infty) = 0, \\
\rW_{\lambda_g} + (1-g(\infty))^{1/s} \rS_s  & \text{if } g(\infty) \in(0,1],
\end{cases}
\end{equation}
where $\rW_{\lambda_g}$ is the random variable with characteristic function as in Def.  \eqref{gtslaw}
with $\lambda_g:(-\infty,0) \to [0,\infty)$  given by
\begin{equation}\label{lgdef}
 \lambda_g(-t) = \tc t^{-s} \int_{t/\vs}^\infty \left( \frac{\vs s}{t} - \frac{s-1}{x} \right) \, g(dx), \qquad t>0, 
\end{equation}
and $\rS_s$ has characteristic function \eqref{stablecf} and is independent of $\rW_{\lambda_g}$.
\end{itemize}
\end{theorem}

\begin{remark}[\textbf{Regularity of the cooling \& $g$ function}]
The $g$ function characterizes the density of increments at scale $\tau(n)^{1/s}$,
and if the cooling map is regular enough to satisfy \eqref{cag} for some continuous $g$ 
with $g(0)=0$,
the above theorem suggests that a generalized tempered stable or a pure stable component
should be expected in the limit. The extra conditions {\bf(S1)} or {\bf(S2)}
are in particular sufficient to guarantee convergence to these types of laws. 
Theorem~\ref{GeneralConditions-S} says nothing about possible Gaussian components
for which increments are on scales smaller than $\tau(n)^{1/s}$.
Moreover the Poisson point process approach used in the proof of
Theorem~\ref{GeneralConditions-S} is not well suited for proving Gaussian limits. 
We further remark that while it is tempting to conjecture that if \eqref{cag}
holds with $g(x) \equiv 1$ then the limit is Gaussian this is not true as
can be seen by Example \ref{mixedEx} in Section \ref{rEx} when $2^{s/2} \leq r < 2$. 
\end{remark}
The previous results give sufficient conditions for convergence to Gaussian, stable, or generalized tempered stable distributions. 
Our next result shows that one can also obtain arbitrary linear combinations of these three types of distributions at least within a certain subclass of generalized laws defined as follows.
Let $\Lconv$ be the class of functions of the form $\lambda(x) = c |x|^{-s-1} a(x)$, where $c>0$ and $a:(-\infty,0]\to [0,1]$ is a convex, non-decreasing function with $\lim_{x\to -\infty} a(x) = 0$  and $a(0) = 1$.

\begin{theorem}[{\bf Mixed laws}]
\label{arbmix}
Let $\mu$ be a fixed $s$-canonical law.
Given $a_1,a_2,a_3 \geq 0$ and a function $\lambda \in \Lconv$,
there exist a cooling map $\tau$
and 
constants $b>0$ and $\beta \in [1/2,1/s]$ such that
the RWCRE $X$ associated to the law
$\mu$ with cooling map $\tau$ satisfies
\begin{equation}\label{triplemix}
  \frac{X_n -\bb{E}[X_n]}{ b n^\beta }
  \limdist a_1 \rN + a_2 \rW_\lambda + a_3 \rS_s,
\end{equation}
with $\rW_\lambda$ as in~\eqref{Xnstableliml}.
\end{theorem}

The proof of the above statement, presented in Section~\ref{S8}, is split
into several steps which in particular offer a constructive procedure to build
the map $\tau$.
This construction is such that the scaling exponent $\beta < 1/s$ only
when there is a Gaussian component in the limit (i.e., when $a_1 > 0$).
However, as can be seen by Example \ref{mixedEx} in Section \ref{rEx}
with $r>2^{s/2}$ this relation between the scaling exponent and the limiting distribution
isn't necessarily true for general cooling maps.

We conjecture that Theorem \ref{arbmix} identifies
all possible limiting distributions that can be obtained
for this model of RWCRE, but if one also allows
for subsequential limits then there are limiting distributions
not covered by Theorem \ref{arbmix}
(see Example \ref{ex:exotic} in Section \ref{rEx}).

Since RWCRE is built upon finite pieces of RWRE,
precise estimates on $Z_n$ are needed in the proofs of the previous results.
We collect in the next theorem the most relevant
such precise estimates which, to the best of our
knowledge, are new and interesting for the analysis of RWRE for $s\in(1,2)$.
The proof of the theorem is given in Section~\ref{S4}
where other concentration estimates for RWRE are derived.

\begin{theorem}[{\bf Stable RWRE: asymptotic $s-$moment,  mean and variance}]
\label{lem:VarZn}
\text{}\\ 
Let $Z$ be  a RWRE with a given $s$-canonical law
$\mu$ as in Def.~\ref{goodenvironment}, with $s \in (1,2)$. Then 
 \begin{equation}\label{LpMom}
   \sup_n \EE\left[ \left| \frac{Z_n-\vs}{n^{1/s}} \right|^p \right] < \infty,
   \qquad \forall p \in (0,s), 
\end{equation}
\begin{equation}\label{meanzn}
 \EE[Z_n] = n\vs + o(n^{1/s}),
\end{equation}
and
\begin{equation}\label{varzn}
  \Var(Z_n) = \sZ^2 n^{3-s}+o(n^{3-s}),
\end{equation}
where $\sZ^2 := 2\tc \vs^{3-s}/(2-s)(3-s)$.
\end{theorem}

The limiting distributions for RWCRE stated above are all given
with centering $\EE[X_n]$ rather than with a linear centering $n\vs$
as in the case of RWRE in \eqref{StableLeftSkew}.  
However, in certain cases once a limiting distribution is obtained
when centered by the mean one can then use \eqref{meanzn}
to show that the same limiting distribution holds when centered by $n\vs$. 
In particular, if condition {\bf{(S1)}} holds then one can check that
$\EE[X_n] - n\vs = o(n^{1/s})$ so that $\frac{X_n - n\vs}{n^{1/s}} \limdist \rS_s$. 
Another consequence of \eqref{meanzn} is that
the stable limit law in \eqref{StableLeftSkew} also holds with centering
$\EE[Z_n]$, that is $\frac{Z_n - \EE[Z_n]}{n^{1/s}} \limdist \rS_s$,
and we will use this fact in the proof of \eqref{Xnstablelim1}.

\section{Proofs: RWRE asymptotics}\label{S4}

The aim of this section is to prove Theorem~\ref{lem:VarZn} and some
other preparatory statements for RWRE related to large and moderate deviations
in the stable regime. In particular, we start in the next two sections with the proofs
of~\eqref{LpMom} and~\eqref{meanzn}, respectively.
Right and left tail estimates are then stated and proven in Section~\ref{LDP} and with
the help of the latter, we derive in Section~\ref{vproof} the asympotics of
the variance in~\eqref{varzn}.
The statements in this section assume without explicit mention that $Z$
is an RWRE with environment law given by an $s$-canonical law $\mu$,
as in Def.~\ref{goodenvironment}.
In many of these proofs, we will make use of the classical RWRE regeneration times
sequence defined via~\eqref{Rk} in Appendix~\ref{Regeneration}.
\subsection{RWRE \texorpdfstring{$L^p$}{Lp} moments estimate\texorpdfstring{: proof of~\eqref{LpMom}}{.}}\label{Lpproof}
The claim in~\eqref{LpMom} is equivalent to
$\EE\left[ \left| Z_n-n\vs \right|^p \right] = \bigo(n^{p/s})$,
and this is what we show below. 
Also, without loss of generality we can assume below that $p\in[1,s)$.
Let $R_{k}$, $k \in \bb{N}_0$, be regeneration times defined in  Appendix~\ref{Regeneration}
and let $k(n)$ be the number of regeneration times by time $n$;
that is $R_{k(n)} \leq n < R_{k(n)+1}$.
Let $\bEE$ denote expectation of RWRE with respect to $\bPP$,
where $\bPP$ is the probability $\bb{P}$ conditioned on the event $\{\inf_{n\geq 0} Z_n = 0\}$.
Also, recall~\eqref{Rk} and let
\begin{equation}\label{R1M}
  c_*:= \frac{1}{\bEE[R_1]}=\frac{1}{\EE[R_2-R_1]}.
\end{equation} 
Then, using the inequality $|a+b+c|^p \leq 3^{p-1}(|a|^p + |b|^p + |c|^p)$ we obtain that
\begin{align}
 \EE\left[ |Z_n-n\vs|^p \right] 
\leq  3^{p-1} &\biggl\{ \EE\left[ \left| Z_n-Z_{R_{k(n)}} - (n-R_{k(n)}) \vs \right|^p \right] \label{uLp1} \\
&\qquad +  \EE\left[ \left| Z_{R_{k(n)}} - Z_{R_{\fl{c_* n}}} - (R_{k(n)}-R_{\fl{c_* n}} ) \vs \right|^p \right] \label{uLp2} \\
&\qquad +  \EE\left[ \left| Z_{R_{\fl{c_* n}}} - R_{\fl{c_* n}} \vs \right|^p \right] \biggr\} \label{uLp3}. 
\end{align}
To complete the proof, in the following paragraphs we prove that
the term in~\eqref{uLp3} is of order $\bigo(n^{p/s})$ and that each of 
the remaining terms is of order $o(n^{p/s})$.

\paragraph{Bound on the term in~\eqref{uLp3}.}
We remark that 
\[
Z_{R_{\fl{c_* n}}} - R_{\fl{c_* n}} \vs = \sum_{k=1}^{\fl{c_* n}} \left( Z_{R_k}-Z_{R_{k-1}} - (R_k-R_{k-1})\vs \right)
\]
is the sum of independent random variables, all of which are i.i.d. except the first.
In view of Lemma~\ref{lem:taumoments} the first term  $Z_{R_1} -R_1\vs$
is negligible for this sum.
It follows from \eqref{vR1}, \eqref{Xtau-tail}, and \eqref{tau-tail} that for $k\geq 2$ the random variables $Z_{R_k}-Z_{R_{k-1}} - (R_k-R_{k-1})\vs$ are zero mean random variables with exponential tails to the right and left tails that are regularly varying of index $-s$.
Thus we can apply  Corollary~\ref{LpSn} to conclude that
the expectation in~\eqref{uLp3} is $\bigo(n^{p/s})$.

\paragraph{Bound on the term in~\eqref{uLp1}.}
For the expectation in~\eqref{uLp1}, note that by the definition of $k(n)$
and the fact that the walk is a nearest neighbor walk, we have that 
\[
\EE\left[ \left| Z_n-Z_{R_{k(n)}} - (n-R_{k(n)}) \vs \right|^p \right] \leq (\vs+1)^p \EE[(R_{k(n)+1}-R_{k(n)})^p].
\]
To control the expectation above, we partition the total probability on
the possible values that $k(n)$ and $R_{k(n)}$ may attain
and then use the i.i.d.\ structure of regeneration times. Explicitly
\begin{align*}
\EE[(R_{k(n)+1}-R_{k(n)})^p]
&=  \EE[R_1^p \ind{R_1 > n}] + \sum_{k=1}^{n} \sum_{m=0}^{n-k} \EE[(R_{k+1}-R_k)^p \ind{R_k=n-m, \, R_{k+1}-R_k > m} ] \\
&=  \EE[R_1^p \ind{R_1 > n}] + \sum_{k=1}^{n} \sum_{m=0}^{n-k} \PP(R_k=n-m) \bEE[R_1^p \ind{R_1 > m}] \\
&= \EE[R_1^p \ind{R_1 > n}] + \sum_{m=0}^{n-1} \left( \sum_{k=1}^{n-m} \PP(R_k=n-m) \right) \bEE[R_1^p \ind{R_1 > m}] \\
&= \EE[R_1^p \ind{R_1 > n}] + \sum_{m=0}^{n-1} \PP(\exists k: \, R_k = n-m) \bEE[R_1^p \ind{R_1 > m}] \\
&\leq \EE[R_1^p \ind{R_1 > n}] + \sum_{m=0}^{n-1} \bEE[R_1^p \ind{R_1 > m}]. 
\end{align*}
The first term in the right hand side is asymptotically vanishing
thanks to Lemma~\ref{lem:taumoments}. Because~\eqref{tau-tail} implies that
$\bEE[R_1^p \ind{R_1 > m}] \sim C' m^{-s+p}$, the sum on the right is $\bigo( n^{1-s+p})$. 
Since $1-s+p < \frac{p}{s}$ when $p<s$,
it follows that the expectation in~\eqref{uLp1} is $o(n^{p/s})$. 
\paragraph{Bound on the term in~\eqref{uLp2}.}
To ease notation we let $W_k := Z_{R_k}-R_k \vs$ for $k\geq 1$. 
Now, fix $\beta \in (1/s,1)$ and $p'$ such that $1\leq p<p'<s$. Then,
\begin{equation}\label{maxWp}
\begin{aligned}
  & \EE\left[ | W_{k(n)}-W_{\fl{c_* n}} |^p  \right]
  \\&\quad
  \leq \EE\left[ \max_{k:|k-c_* n| \leq  n^\beta} |W_k-W_{\fl{c_* n}}|^p \right]
  + 2 \EE\left[ \max_{k\leq n}|W_k|^p \ind{|k(n)-c_* n| > n^\beta}  \right]
  \\
  &\quad \leq 2 \bEE\left[ \max_{k\leq 2 n^\beta} |W_k|^{p'} \right]^{\frac{p}{p'}}
  +  2   \left( \EE\left[ \max_{k\leq n} |W_k|^{p'} \right] \right)^{\frac{p}{p'}} \PP(|k(n)-c_* n| > n^{\beta})^{1-\frac{p}{p'}}
  \\ 
  &\quad \leq  C \left( \bEE\left[ |W_{\fl{2n^\beta}} |^{p'} \right]\right)^{\frac{p}{p'}}
  + C  \left( \EE\left[ |W_n|^{p'} \right] \right)^{\frac{p}{p'}} \PP(|k(n)-c_* n| > n^{\beta})^{1-\frac{p}{p'}},
\end{aligned}
\end{equation}
where in the second inequality we used the i.i.d.\ structure of the regeneration times
for the first term and H\"older's inequality for the second term,
and in the last inequality we used the $L^p$-maximal inequality
for martingales with $C = 2\left( \frac{p'}{p'-1} \right)^p$. 
As noted in the analysis of~\eqref{uLp3},
the two expectations in the last line above can be bounded using
Corollary~\ref{LpSn} (and also Lemma \ref{lem:taumoments} for the second expectation).
Thus, we get 
\begin{equation}\label{e:3}
  \EE\left[ | W_{k(n)}-W_{\fl{c_* n}} |^p  \right]
  = \bigo( n^{\beta p/s}) + \bigo( n^{p/s} ) \PP(|k(n)-c_* n| > n^{\beta})^{1-\frac{p}{p'}}. 
\end{equation}
Therefore, since $\beta<1$ and $1 - s\beta<0$, to finish it is enough to prove that 
\begin{equation} \label{knmd}
 \PP(|k(n)-c_* n| > n^{\beta}) = \bigo(n^{1-s\beta}).
\end{equation}
Since $\PP(|k(n)-c_* n| > n^{\beta})\leq \PP(k(n) < c_* n -n^{\beta}) + \PP(k(n) > c_* n + n^{\beta})$, we  have that
\begin{align}
  & \PP(|k(n)-c_* n| > n^{\beta}) \leq  \PP(R_{\fl{c_* n-n^\beta}+1} > n ) + \PP(R_{\lceil{c_* n+n^\beta}\rceil} \leq n) \\
&\leq \PP\left( R_1 > \frac{n^{\beta}}{2c_*} \right)  +  \bPP\left(R_{\fl{c_* n-n^\beta}} > n - \frac{n^{\beta}}{2c_*} \right) + \bPP(R_{\lceil{c_* n+n^\beta}\rceil-1} \leq n).\label{go}
\end{align}
To complete the proof we estimate the three terms in the right hand side above as follows.
The first term is 
$\bigo(n^{-\beta q})$ for any $q<s$
as can been seen by the Markov inequality
and Lemma~\ref{lem:taumoments}.
For the second term if we let $m = \fl{c_* n - n ^\beta}$ we obtain
\begin{equation}\label{e:n1}
\begin{aligned}
  & \bPP\left(R_{m} > n - \frac{n^{\beta}}{2c_*} \right) 
  = \bPP\left( \sum_{k=1}^{m} \left( R_k-R_{k-1} - \bEE[R_1] \right) > n-\frac{n^{\beta}}{2c_*}-
\frac{1}{c_*} m \right) \\
& 
\leq \bPP\left( \sum_{k=1}^{m} \left( R_k-R_{k-1} - \bEE[R_1] \right) > Cn^{\beta-1/s}    
m^{1/s}  \right) = \bigo(n^{1-s\beta}),
\end{aligned}
\end{equation}
where in the inequality we used that  for some $C>0$
\begin{equation}\label{cbeta}
n-\frac{n^{\beta}}{2c_*}-\frac{1}{c_*} m > \frac{3}{2 c_*}n^\beta =  \frac{3}{2 c_*}n^{\beta-1/s} n^{1/s} >C n^{\beta-1/s} m^{1/s},
\end{equation}
and for the last step in \eqref{e:n1} we have used the Lemma~\ref{Sntail}
for the i.i.d.\ sequence of zero-mean variables
$(R_k-R_{k-1} - \bEE[R_1])_{k\geq 1}$.

For the third probability in~\eqref{go},
using that $\bEE[R_1] = 1/c_*$ and $m = \lceil c_* n-n^\beta\rceil -1$ by the same argument in~\eqref{cbeta}
we have  that
\begin{align*}
 \bPP(R_{m} \leq n)
&= \bPP\left( \sum_{k=1}^{m} \left( R_k-R_{k-1} - \bEE[R_1] \right) < n - \frac {m}{c_*} \right) \nonumber \\
&\leq  \bPP\left( \sum_{k=1}^{m} \left( R_k-R_{k-1} - \bEE[R_1] \right) < -C n^{\beta-1/s} m^{1/s}  \right),
\end{align*}
we can again invoke the tail decay of regeneration times in~\eqref{tau-tail},
which together with the bound in Lemma~\ref{Sntail-left} gives that
this last probability is bounded above by 
$C' e^{-c n^{\frac{\beta s-1}{s-1}}}$ for some constants $c,C'>0$,
which concludes the proof.
\qed

\subsection{RWRE mean estimate}\label{mproof}
In this section we prove~\eqref{meanzn}.
Consider the regeneration times sequence
and as in the previous section let $k(n)$ denote the number of regeneration times by time $n$. 
For any $c_*>0$ and $n\in\N$ we can write
\begin{equation}\label{meanDec}
\begin{split}
 \EE[Z_n]
= & n\vs  + \EE[Z_{R_{\fl{c_* n}}}- R_{\fl{c_* n}} \vs]
 + 
 \EE\left[Z_n - Z_{R_{k(n)}} - (n-R_{k(n)}) \vs\right]\\
 &+ \EE\left[Z_{R_{k(n)}}- Z_{R_{\fl{c_* n}}} - (R_{k(n)} - R_{\fl{c_* n}}) \vs\right]. 
\end{split}
\end{equation}

Letting  $c_* = \frac{1}{\bEE[R_1]}$ as in~\eqref{R1M}
and using~\eqref{vR1} together with~\eqref{tme} and~\eqref{ime},
we then have that the second term in the right hand side above stays bounded,
that is
\[
\EE[Z_{R_{\fl{c_* n}}}- R_{\fl{c_* n}} \vs]=\bigo(1).
\]
On the other hand, by arguing as in the proof in Section~\ref{Lpproof}
for the terms~\eqref{uLp1} and~\eqref{uLp2}, respectively, with $p=1$,
we also have that the last 
two terms in the right hand side of~\eqref{meanDec} are $o(n^{1/s})$.\qed

\subsection{RWRE tail estimates}\label{LDP}

The main results in this section are right and left tail estimates for the RWRE which range from the limiting distribution scale all the way to the large deviation scale. 
We begin with estimates on the right tail. 


\begin{lemma}[{\bf Right tail estimate}]\label{Znrighttail}\text{}\\ 
There exist constants $a,c,C>0$ such that for all $n$ large enough and $0<t<a n^{1-\frac{1}{s}}$
\begin{equation}\label{e:Zn}
\PP(Z_n - n\vs > tn^{1/s}) \leq C e^{-c t^{s/(s-1)}}.
\end{equation}
\end{lemma}
\begin{proof}
Note that it is enough to prove~\eqref{e:Zn}
for $\delta \leq t < a n^{1-\frac{1}{s}}$ since we may extend the bound to $0<t<\delta$ if we take
the constant $C$ in front of the exponential large enough.
Thus, for the remainder of the proof we will assume that $ \delta \leq t < a n^{1-\frac{1}{s}}$. 
  
First of all, note that for any choice of $m \in \bb{N}$ 
\begin{equation}\label{mdspeedup1}
\begin{aligned}
&\PP(Z_n > n\vs + t n^{1/s} ) 
 \leq \PP( Z_{R_m} > n\vs + t n^{1/s} ) + \PP(R_m < n ) \\
& \leq \PP\left(Z_{R_1} > \frac{(1-\vs)t n^{1/s}}{2} \right) + \bPP\left( Z_{R_{m-1}} > n\vs + \frac{(1+\vs)t n^{1/s}}{2} \right) + \bPP(R_{m-1} < n ). 
\end{aligned}
\end{equation}
Since $Z_{R_1}$ has an exponential tail due to~\eqref{ZR1tail},
the first probability on the right is
bounded by $C_1 e^{-c_1 t n^{1/s}}$ for some constants $C_1,c_1>0$. 

For the analysis of the last two terms in~\eqref{mdspeedup1}
we let $m = m(n,t) =1 + \fl{c_*}(n+t n^{1/s})$,
where $c_*= \frac{1}{\bEE[R_1]}$ as in~\eqref{R1M}. 
Using that  $\bEE[Z_{R_1}] = \vs \bEE[R_1] = \vs/{c_*}$ we have 
\begin{align*}
  & \bPP\left( Z_{R_{m-1}} > n\vs + \frac{(1+\vs)t n^{1/s}}{2}  \right)
  \leq \bPP\left( \sum_{k=1}^{m -1} \left( Z_{R_k}-Z_{R_{k-1}} - \bEE[Z_{R_1}]  \right) > \frac{1-\vs}{2} tn^{1/s}  \right). 
\end{align*}
Since the random variables in the sum inside the last probability are
i.i.d.\ with exponential tails (Corollary~\ref{cor:regtails}),
it follows from the large deviation estimates in~\cite[Thm. III.15]{Pe75},
that there exist constants $a, c_2>0$ so that this probability
is bounded above by $e^{-c_2 t^2 n^{\frac{2}{s}-1}}$ for all $t \leq a n^{1-\frac{1}{s}}$. 

For the third probability in~\eqref{mdspeedup1},
since $\bEE[R_1] = 1/c_*$,  for $n$ large enough we have that
\begin{equation}\label{e:n2}
\begin{aligned}
\bPP(R_{m(n,t)-1} < n )
&\leq  \bPP\left( \sum_{k=1}^{\fl{c_*(n+t n^{1/s})}} \left( R_k-R_{k-1} - \bEE[R_1] \right) < \frac{-t n^{1/s}}{2c_*}  \right).
\end{aligned}
\end{equation}
It follows from the tail decay of regeneration times in~\eqref{tau-tail}
and the large deviation bound in Lemma~\ref{Sntail-left} that
there is a constant $c_3>0$ such that this last probability is bounded above
by $e^{-c_3 t^{\frac{s}{s-1}}}$ for all $t \leq a n^{1-\frac{1}{s}}$. 

Combining the above upper bounds for the three terms in~\eqref{mdspeedup1}, we have that for $n$ large enough
\begin{equation}\label{e:4}
 \PP(Z_n - n\vs > t n^{1/s}) \leq C_1 e^{-c_1 t n^{1/s}} + e^{-c_2 t^2 n^{\frac{2}{s}-1}} + e^{-c_3 t^{\frac{s}{s-1}}}
\leq C e^{-c t^{\frac{s}{s-1}} }, 
\end{equation}
where again in the last equality we used that $t \leq a n^{1-\frac{1}{s}}$.
\end{proof}


\begin{corollary}[{\bf Asymptotics on positive part of the variance}]\label{cor:L2}
\text{}\\ 
For any $s\in(1,2)$, the following asymptotics in $n$ is valid: 
$$\EE[((Z_n-\EE[Z_n])_+)^2] = \bigo(n^{2/s}).$$ 
\end{corollary}

\begin{proof}
By~\eqref{meanzn} it is enough to prove that $\EE[((Z_n-n\vs)_+)^2] = \bigo(n^{2/s})$. 
By Lemma~\ref{Znrighttail} and the fact that $Z_n \leq n$ we have that 
\begin{equation}\label{PosB}
 \begin{aligned}
  \EE[((Z_n-n\vs)_+)^2]
&\leq \int_0^{an} x \PP( Z_n -n\vs > x ) \, dx + n^2 \PP( Z_n - n\vs \geq a n ) \\
&= n^{2/s} \int_0^{an^{1-\frac{1}{s}}}  t \PP( Z_n -n\vs > t n^{1/s} ) \, dt + n^2 \PP( Z_n - n\vs \geq a n ) \\
&\leq n^{2/s} \int_0^\infty  t C e^{-c t^{\frac{s}{s-1}}} \, dt + C n^2 e^{-c a^{\frac{s}{s-1}} n} 
= \bigo(n^{2/s}).
 \end{aligned}
\end{equation}
\end{proof}

We next turn our attention to left tail estimates for the RWRE. 
Note that \eqref{precisetail} gives very precise left tail asymptotics,
but over a region that doesn't quite cover all of the moderate devations we are interested in.
The following Lemma gives a weaker bound but over a scale that covers the entire moderate deviation regime. 

\begin{lemma}[{\bf General left tail estimates}]\label{lem:Znmoddev}\text{}\\
There exist constant $C,C'<\infty$ such that for large enough $n$
\begin{align}
 \PP(Z_n - n\vs \leq -t n^{1/s}) &\leq C t^{-s}, \quad \forall t \leq \frac{\vs}{2} n^{1-\frac{1}{s}}, \label{Znmd} \\
\text{and}\qquad
 \PP(Z_n - \EE[Z_n] \leq -t n^{1/s}) &\leq C' t^{-s}, \quad \forall t \leq \frac{\vs}{2} n^{1-\frac{1}{s}}. \label{Znmd-c}
\end{align}
\end{lemma}
\begin{proof}[Proof of Lemma~\ref{lem:Znmoddev}]
First of all, since $\EE[Z_n] =  n\vs + o(n^{1/s})$, we only need to prove~\eqref{Znmd}. 
Moreover, we take $C\geq 1$, so it suffices to prove~\eqref{Znmd} for $1 \leq t \leq  \frac{\vs}{2} n^{1-\frac{1}{s}}$.

As with the proof of Lemma \ref{Znrighttail} we will once again use regeneration times.
For any $m \geq 1$, since $R_m \leq n$ implies $Z_{R_m} \leq Z_n$ we have that
\begin{align}
& \PP ( Z_n - n\vs \leq -t n^{1/s} ) 
\leq \PP( R_m > n ) + \PP(Z_{R_m} \leq n\vs - t n^{1/s} ) \nonumber \\
&\qquad\qquad\leq \PP\left( R_1 > \frac{t}{2} n^{1/s} \right) + \bPP(R_{m-1} > n-\frac{t}{2}n^{1/s} ) + \bPP( Z_{R_{m-1}} \leq n\vs - t n^{1/s} ) \label{mdslowdown}
\end{align}
For the first term in \eqref{mdslowdown}, note that Lemma \ref{lem:taumoments}
  implies that $\EE[R_1^{s-1}]<\infty$ and therefore
\begin{equation}\label{mdslowdown1}
 \PP\left( R_1 > \frac{t}{2} n^{1/s} \right) \leq C t^{-s+1}n^{-1+\frac{1}{s}} \leq C' t^{-s}, \quad \forall t \leq \frac{\vs}{2} n^{1-\frac{1}{s}}
\end{equation}
To bound the last two terms in \eqref{mdslowdown}, we will let $m = 1 + \fl{c_*(n-tn^{1/s})}$ where again $c_* = \frac{1}{\bEE[R_1]}$ so that 
for $t\geq 1$ and $n$ sufficiently large we have 
\begin{align}
& \bPP(R_{m-1} > n-\frac{t}{2}n^{1/s} ) + \bPP( Z_{R_{m-1}} \leq n\vs - t n^{1/s} ) \nonumber \\
&\quad \leq \bPP\left( \sum_{k=1}^{\fl{c_*(n-tn^{1/s})}} \left( R_k - R_{k-1} - \bEE[R_1] \right) > \frac{t}{2} n^{1/s} \right) \label{mdslowdown2} \\
&\quad\qquad + \bPP\left( \sum_{k=1}^{\fl{c_*(n-tn^{1/s})}} \left( Z_{R_k} - Z_{R_{k-1}} - \bEE[Z_{R_1}] \right) < -\left(\frac{1-\vs}{2} \right) t n^{1/s}   \right).  \label{mdslowdown3}
\end{align}
Thanks to \eqref{tau-tail}, we may apply Lemma \ref{Sntail}
  to obtain that the probability in \eqref{mdslowdown2} is bounded by $C t^{-s}$ for $n$ large,
while since \eqref{Xtau-tail} implies the random variables
inside the sum in \eqref{mdslowdown3} have exponential tails
we can again use \cite[Thm. III.15]{Pe75} to bound this last probability by 
$e^{-c t^2 n^{\frac{2}{s}-1}} \leq e^{-c' t^{\frac{s}{s-1}} }$, where the last inequality
holds since $t \leq \frac{\vs}{2} n^{1-\frac{1}{s}}$. 
Finally, since $e^{-c' t^{\frac{s}{s-1}} } \leq C t^{-s}$ for some $C>1$ and all
$t\geq 1$ this completes the proof of the lemma. 
\end{proof}

The following corollary gives a simple extension of the precise left tail asymptotics from ~\eqref{precisetail}
when we center $Z_n$ with the mean rather than by $n\vs$. 
Combined with the more general left tail bound in Lemma \ref{lem:Znmoddev} this then gives
a truncated second moment bound \eqref{truncVar} that is instrumental for the proofs to come.  

\begin{corollary}[{\bf Moderate slow-down deviations centering with mean}]\label{cor:prectail-c}\text{}\\
Set $I_n : = [n^{1/s}(\log n)^4, n\vs -n^{1/s} \log n]$ as in~\eqref{precisetail}, then
\begin{equation}\label{moderatecenter}
 \lim_{n\to\infty}  \sup_{ t \in I_n}  \left| \frac{\PP(Z_n - \EE[Z_n] < -t )}{(n\vs-t)t^{-s}} - \tc \right| = 0. 	
\end{equation}
Furthermore, there exists a constant $C<\infty$ such that for $n$ sufficiently large
\begin{equation}\label{truncVar}
  \EE[(Z_n - \EE[Z_n])^2 \ind{Z_n - \EE[Z_n] \in (-t,0)} ] \leq C n t^{2-s},
  \qquad \forall t \leq n\vs - n^{1/s}(\log n).  
\end{equation}
\end{corollary}
\begin{proof}
 First of all, note that 
\begin{equation}\label{centervn}
 \PP(Z_n - \EE[Z_n] < -t ) = \PP\left(Z_n - n \vs < -(t-\EE[Z_n]+n\vs) \right). 	
\end{equation}
If $n^{1/s}(\log n)^4 \leq t \leq n\vs-(\log n) n^{1/s}$,
then for $n$ sufficiently large from~\eqref{meanzn} it follows that 
\begin{equation}\label{adjust}
(\log n)^3 n^{1/s} \leq t-\EE[Z_n] + n\vs \leq n \vs -\frac{1}{2} (\log n) n^{1/s} < n\vs - \log n. 
\end{equation}
Therefore, $I_n \subset \tilde{I}_n$ and  we can apply
the tail asymptotics~\eqref{precisetail} with $t-\EE[Z_n]+n\vs$ in place of $t$. 
That is, we may write $ \PP(Z_n - \EE[Z_n] < -t )/[(n\vs-t)t^{-s}]$ as 
\begin{equation}\label{3frac}
\begin{aligned}
  & \frac{\PP\left(Z_n - n \vs < -(t-\EE[Z_n]+n\vs) \right) }{ (\EE[Z_n]-t)(t-\EE[Z_n]+n\vs)^{-s} }
  \left(\frac{\EE[Z_n] - t}{n\vs-t} \right)\left( \frac{t}{t-\EE[Z_n]+n\vs} \right)^s.
\end{aligned}
\end{equation}
To complete the proof  of~\eqref{moderatecenter}, note that  as $n\to\infty$
the first term on~\eqref{3frac} converges to $\tc$ and the last two terms converge to 1 uniformly in  $t \in I_n$.
Note that for the convergence of the last two terms to 1 we again use~\eqref{meanzn}.

We next show~\eqref{truncVar}.
By the tail estimate in equation~\eqref{moderatecenter}
and Lemma~\ref{lem:Znmoddev}, we see that 
$\PP(Z_n - \EE[Z_n] < -x ) \leq Cn x^{-s}$ for $n$ large enough and
$0 < x \leq n\vs - n^{1/s}(\log n)$.
Therefore, if $t\leq n\vs - n^{1/s}(\log n)$
and $n$ is large enough we have
\begin{equation}\label{e:n4}
\begin{aligned}
&  \EE[(Z_n - \EE[Z_n])^2 \ind{Z_n - \EE[Z_n] \in (-t,0)} ] 
= \int_0^t 2x \PP\left( -t < Z_n - \EE[Z_n] \leq -x \right) \, dx \\
&\qquad\leq \int_0^t 2x \PP\left( Z_n - \EE[Z_n] \leq -x \right) \, dx 
\leq 2Cn \int_0^t x^{1-s} \, dx
= 2Cn t^{2-s}. 
\end{aligned}
\end{equation}

\end{proof}

\subsection{RWRE variance asymptotics}\label{vproof}
In this section we prove~\eqref{varzn}.
By \eqref{meanzn} and then \eqref{PosB} we have that 
\begin{align*}
 \Var(Z_n)
&= \EE[(Z_n-n\vs)^2] + o(n^{2/s}) \\
&= \EE[\left((Z_n-n\vs)_-\right)^2] + \bigo(n^{2/s}) 
= 2\int_0^\infty t \PP(Z_n-n\vs < -t) \, dt + o(n^{3-s}), 
\end{align*}
where the last equality follows from the fact that $3-s > \frac{2}{s}$ when $s \in (1,2)$. 
It remains to show that the integral term, when multiplied by $n^{s-3}$ converges to $\sZ^2 = \frac{2\tc \vs^{3-s}}{(2-s)(3-s)}$ as $n\to\infty$. 
To this end, fixing a $\delta\in(0,\vs/2)$, we have that $n^{s-3}$ times this integral term can be decomposed as
\begin{equation}\label{breakparts}
\begin{aligned}
& 2n^{s-3}\int_0^{n\delta} t \PP(Z_n-n\vs < -t)\, dt 
+
2n^{s-3}\int_{n\delta}^{n(\vs-\delta)} t \PP(Z_n-n\vs < -t)\, dt \\
&\qquad +
2n^{s-3}\int_{n(\vs-\delta)}^{n(\vs+1)} t \PP(Z_n-n\vs < -t)\, dt
\\&=:
I+II+III.
\end{aligned}
\end{equation}
The truncation of the integrals up to $t\leq n(\vs+1)$ is due to the fact that $|Z_n|\leq n $.
We will show below that the main contribution to the sum in \eqref{breakparts}
will come from $II$ while $I$ and $III$ will be vanishingly small as $\delta \to 0$.
For $II$, 
we see that for large enough $n$ the interval $[{n\delta},{n(\vs-\delta)}]$
is contained in $\tilde{I}_n$ from \eqref{precisetail}.
Therefore by~\eqref{precisetail} we obtain that
\begin{equation}\label{IId}
\begin{aligned}
II&= 
2n^{s-3}\int_{n\delta}^{n(\vs-\delta)}t (K_0+o(1))(n\vs-t)t^{-s} \, dt \\
&= 2K_0 \left[ \vs\frac{\left(\vs-\delta\right)^{2-s}- \delta^{2-s}}{2-s}
  - \frac{\left(\vs-\delta\right)^{3-s}- \delta^{3-s}}{3-s}  \right] + o(1).
\end{aligned}
\end{equation}
Therefore, $II$  goes to $\frac{2K_0\vs^{3-s}}{(2-s)(3-s)}$ as first $n\to \infty$ and then $\delta \to 0$. 

We now show the negligibility of the other terms  in~\eqref{breakparts}
as $n\to \infty$ and then $\delta\to 0$.
For the first term $I$,
using a substitution $z=tn^{-1/s}$ and then applying ~\eqref{Znmd} 
we obtain for large enough $n$ the bound
\begin{align*}
 I &= 2 n^{s-3+\frac{2}{s}} \int_0^{\delta n^{1-\frac{1}{s}}} z \PP(Z_n - nv \leq -z n^{1/s}) \, dz 
\leq C n^{s-3+\frac{2}{s}} \int_0^{\delta n^{1-\frac{1}{s}}} z^{1-s} \, dz
= \frac{C \delta^{2-s}}{2-s}. 
\end{align*}
For the term $III$ in~\eqref{breakparts}, since the probabilities in the integral are decreasing in $t$ we have for $n$ large enough that \eqref{precisetail} implies
\begin{align*}
III
&\leq 2n^{s-3}\PP(Z_n-n\vs < -n(\vs-\delta)) \int_{n(\vs-\delta)}^{n(\vs+1)} t \, dt\\
&\leq 4n^{s-3}\tc (n\delta)(n\vs-n\delta)^{-s}\frac{n^2\left(\vs+1\right)^{2}}{2}
= \frac{2\tc \delta (\vs+1)^2}{(\vs-\delta)^s}. 
\end{align*}
Since the bounds of $I$ and $III$ above hold for $n$ sufficiently large and vanish as $\delta \to 0$, this completes the proof of the asymptotics of $\Var(Z_n)$. 
\qed

\section{Proofs: pure Gaussian limits}\label{S5}

We first prove Theorem \ref{GeneralConditions-N} and then treat the Gaussian limits in Theorem \ref{PolyTransition} as a subcase.

\subsection{Gaussianity for well-behaving cooling maps}
In this section we 
prove Theorem \ref{GeneralConditions-N}.
We will first prove the limiting distribution along the subsequence of times $\tau(n)$
\begin{equation}\label{Xtaun-Gauss}
 \frac{X_{\tau(n)} - \EE[X_{\tau(n)}]}{\sqrt{\Var(X_{\tau(n)})}} \limdist \rN, \quad \text{as } n\to \infty,
\end{equation}
and then extend the result to all times. 

\paragraph{Gaussian limits for $X_{\tau(n)}$.}
For ease of notation, let 
\begin{equation}\label{varscale}
\mc{Z}^k_n:=\frac{Z_{T_k}^{(k)} - \bb{E}[Z_{T_k}^{(k)}]}{\sqrt{\Var(X_{\tau(n)})}}, 
\end{equation}
so that \eqref{Xtaun-Gauss} becomes $\sum_{k=1}^n \mc{Z}^k_n \limdist \rN$.
Now we note that $\{\mc{Z}^k_n\}_{ k \leq n}$ is a triangular array
composed of independent mean $0$ random variables such that $\sum_{k=1}^{n} \Var(\mc{Z}^k_n) = 1$ for all $n\in \N$. 
To prove \eqref{Xtaun-Gauss} we will check the Lindeberg condition~\cite[Thm 7.3.1, p. 307]{AshDol00} for this triangular array.
In particular, we need to check that the triangular array is uniformly asymptotically negligible, 
 \begin{equation}\label{uan}
   \lim_{n\to\infty} \sup_{k\leq n } \bb{P}\crt{|\mc{Z}^k_n|\geq \varepsilon} = 0, \qquad\forall \varepsilon>0, 
 \end{equation}
and also that Lindeberg's condition holds 
\begin{equation}\label{Lindeberg}
\lim_{n\to\infty} \sum_{k=1}^{n} \EE\left[\prt{\mc{Z}^k_n}^2 \ind{|\mc{Z}^k_n| > \epsilon}\right] = 0, \quad \forall\;\varepsilon>0.
\end{equation}
Since 
the walk is nearest neighbor and since
\eqref{varzn} implies $\Var(X_{\tau(n)}) \geq c \sum_{k=1}^n T_k^{3-s}$ for some $c>0$, it follows that 
\[
\left| \mc{Z}^k_n \right| = \left| \frac{Z^{(k)}_{T_k} - \EE[Z^{(k)}_{T_k}]}{\sqrt{\Var(X_{\tau(n)})}} \right|
\leq \frac{2 T_k}{ \sqrt{c \sum_{k=1}^n T_k^{3-s}}}, \quad 1\leq k\leq n.  
\]
The assumption \eqref{Gauss_iff} implies that the right side vanishes uniformly in $k\leq n$ as $n\to\infty$. 
In particular, this implies that for any fixed $\epsilon>0$ and $n$
sufficiently large we have with probability 1 that $|\mc{Z}^k_n|\leq \epsilon$ for all $k\leq n$. 
Thus, \eqref{uan} and \eqref{Lindeberg} both follow.

\paragraph{Irrelevance of boundary term.}
To extend the limiting distribution from the subsequence $\tau(n)$ to all times, we use the decomposition \eqref{sumpiece} to write 
\begin{equation}\label{byeboundary}
  \frac{X_n - \bb{E}[X_n]}{\sqrt{\Var(X_n)}} 
=   \frac{X_{\tau(\ell_n)} - \bb{E}[X_{\tau(\ell_n)}]}{\sqrt{\Var(X_{\tau(\ell_n)})}} \sqrt{\frac{ \Var(X_{\tau(\ell_n)})}{\Var(X_n)} } 
+ \frac{Z^{(\ell_n+1)}_{\bar{T}_n} - \EE[Z^{(\ell_n+1)}_{\bar{T}_n}]}{\sqrt{\Var(X_n)}} .  
\end{equation}
It follows from \eqref{Xtaun-Gauss} that the first term on the right converges in distribution to $\rN$ as $n\to\infty$. Thus, it is enough to show that 
\begin{equation}\label{boundaryclaims}
 \lim_{n\to\infty} \frac{ \Var(X_{\tau(\ell_n)})}{\Var(X_n)} = 1, 
\quad\text{and}\quad 
\frac{Z^{(\ell_n+1)}_{\bar{T}_n} - \EE[Z^{(\ell_n+1)}_{\bar{T}_n}]}{\sqrt{\Var(X_n)}} \underset{n\to\infty}{\limdist} 0. 
\end{equation}

For the first claim in \eqref{boundaryclaims}, since $\Var(X_{\tau(\ell_n)}) \leq \Var(X_n) \leq \Var(X_{\tau(\ell_n)}) + \Var(Z_{T_{\ell_n+1}})$ we will show that $\lim_{\ell\to\infty} \frac{\Var(Z_{T_{\ell+1}})}{\Var(X_{\tau(\ell)})} = 0$. 
To this end,
by \eqref{varzn}
there is a $C<\infty$ such that 
\begin{equation}\label{bc1}
 \frac{\Var(Z_{T_{\ell+1}})}{\Var(X_{\tau(\ell)})} \leq C \frac{T_{\ell+1}^{3-s}}{\sum_{k=1}^{\ell} T_k^{3-s} }. 
\end{equation}
We will then show that the assumption \eqref{Gauss_iff} implies the right side vanishes as $\ell\to\infty$. 
Indeed, given any $\epsilon > 0$, \eqref{Gauss_iff}
implies that for $\ell$ sufficiently large we have $T_{\ell+1}^2 \leq \epsilon \sum_{k=1}^{\ell+1} T_k^{3-s}$.
Using the fact that $T_{\ell+1}^{3-s} \leq T_{\ell+1}^2$
since $s\in (1,2)$ this then implies that
$(1-\epsilon)T_{\ell+1}^{3-s} \leq \epsilon \sum_{k=1}^\ell T_k^{3-s}$.
Combined with \eqref{bc1} this implies that
$\frac{\Var(Z_{T_{\ell+1}})}{\Var(X_{\tau(\ell)})} \leq \frac{C\epsilon}{1-\epsilon}$
for $\ell$ sufficiently large.
Since $\epsilon>0$ was arbitrary, this completes the proof of the first claim in \eqref{boundaryclaims}. 

For the second claim in \eqref{boundaryclaims}, since \eqref{StableLeftSkew} implies that $\{ \frac{Z_{\bar{T}_n} - \EE[Z_{\bar{T}_n}]}{\bar{T}_n^{1/s}} \}_{n\geq 1}$ is tight, it is enough to show that 
\begin{equation}\label{bc2}
\lim_{n\to\infty}  \frac{\bar{T}_n^{1/s} }{\sqrt{\Var(X_n)}} = 0. 
\end{equation}
To obtain an upper bound $\frac{\bar{T}_n^{1/s} }{\sqrt{\Var(X_n)}}$, it follows from \eqref{varzn} that $\Var(X_n) \geq \Var(Z_{\bar{T}_n}) \geq c \bar{T}_n^{3-s}$ for some $c>0$, so that 
\[
 \frac{\bar{T}_n^{1/s} }{\sqrt{\Var(X_n)}} \leq \frac{\bar{T}_n^{1/s} }{\sqrt{c} \bar{T}_n^{\frac{3-s}{2}}} = 
\frac{1}{\sqrt{c}} \bar{T}_n^{\frac{1}{s} - \frac{3-s}{2}}. 
\]
Since $\frac{1}{s} - \frac{3-s}{2} = \frac{-(2-s)(s-1)}{2s} < 0$ for $s \in (1,2)$, this upper bound becomes vanishingly small as $\bar{T}_n$ becomes large. 
Therefore, if we use this upper bound when $\bar{T}_n \geq \Var(X_n)^{s/4}$ and when $\bar{T}_n \leq \Var(X_n)^{s/4}$ we use that $\frac{\bar{T}_n^{1/s} }{\sqrt{\Var(X_n)}} \leq \Var(X_n)^{-1/4}$, then since $\Var(X_n) \to \infty$ we have that \eqref{bc2} follows. 
This completes the proof of Theorem \ref{GeneralConditions-N}.

\subsection{ CLT for polynomial increments}
In this section we prove~\eqref{gauss_slow}.
For two real valued functions we say that $f(n) \sim g(n)$ if
$\lim_{n \to \infty} \frac{f(n)}{g(n)} = 1$
and we write $f(n) = \Theta(g(n))$ if there are constants $c,C>0$
for which $c < \frac{f(n)}{g(n)} < C$ for all $n$.
We consider a polynomially growing cooling increments $T_k \sim A k^{a}$ as in~\eqref{polygrow}.
In this case, it follows that
$\sup_{k\leq n} \frac{T_k}{\sqrt{\sum_{k=0}^n T_k^{(3-s)}}} =\Theta\left( n^{a-\frac{a(3-s)+1}{2}}\right)$,
from which we see that condition~\eqref{Gauss_iff}
is satisfied only when $a-\frac{a(3-s)+1}{2}<0$, i.e. for $a<1/(s-1)$.
Thus, if $a<1/(s-1)$ applying Theorem \ref{GeneralConditions-N}
implies that $\frac{X_n - \EE[X_n]}{\sqrt{\Var(X_n)}} \limdist \rN$. 
To finish the proof of \eqref{gauss_slow} it remains only to show that
$\Var(X_n) \sim B^2 n^{2\beta}$ with the constants $B$ and $\beta$
as given in the statement of Theorem \ref{PolyTransition}. 

It follows from \eqref{varzn} that
$\sum_{k=1}^{\ell_n} \Var(Z_{T_k}) \sim \sum_{k=1}^{\ell_n} \sZ^2 A^{3-s} k^{a(3-s)} \sim \frac{\sZ^2 A^{3-s}}{a(3-s)+1} \ell_n^{a(3-s)+1}$, 
and since $\tau(n) \sim \frac{A}{a+1} n^{a+1}$
implies that $\ell_n \sim (\frac{a+1}{A})^{1/(a+1)} n^{1/(a+1)}$ it follows that 
\begin{equation}\label{VXtl}
 \sum_{k=1}^{\ell_n} \Var(Z_{T_k}) 
\sim \frac{\sZ^2A^{\frac{2-s}{a+1}}(a+1)^{\frac{a(3-s)+1}{a+1}}}{a(3-s)+1} n^{\frac{a(3-s)+1}{a+1}}.
\end{equation}
Another application of \eqref{varzn} implies that $\Var(Z_{\bar{T}_n})=\bigo(\bar{T}_n^{3-s}) =\bigo( n^{a(3-s)})$, and since $a<\frac{1}{s-1}$ and $1<s<2$ imply that $a(3-s) < \frac{a(3-s)+1}{a+1}$, it then follows that
$\Var(X_n) = \sum_{k=1}^{\ell_n} \Var(Z_{T_k}) + \Var(Z_{\bar{T}_n}) \sim \sum_{k=1}^{\ell_n} \Var(Z_{T_k}) $. 
Comparing with \eqref{VXtl} and recalling the formula for $\sZ^2$ in Theorem \ref{lem:VarZn}, this completes the proof of $\Var(X_n) \sim B^2 n^{2\beta}$ with  
with $\beta = \frac{a(3-s)+1}{2(a+1)}$ and 
$B^2  := \frac{2\tc \vs^{3-s}A^{\frac{2-s}{a+1}}(a+1)^{\frac{a(3-s)+1}{a+1}}}{(2-s)(3-s)(a(3-s)+1)}$.
\qed

\section{Proofs: generalized tempered and stable limits }\label{S7}

In this section we will prove the general Theorem~\ref{GeneralConditions-S} and then deduce \eqref{Xnstableliml} and \eqref{Xnstablelim} from it.

The chapter is organized as follows.
We split the proof of the general theorem in two main parts corresponding to the two different statements, \eqref{Xnstablelim1} and \eqref{Xnstablelim2}, respectively, which in particular will require two different proof strategies.

The first part is presented in Section \ref{FastCool}
where we show \eqref{Xnstablelim1} under assumption {\bf (S1)}.
Typical examples that satisfy the first requirement in {\bf (S1)} are cooling sequences that grow
very rapidly (e.g. exponentially fast)\footnote{However, the case $r > 2^s$
in Example~\ref{mixedEx} in Section \ref{rEx} shows
that there are cooling maps that satisfy {\bf(S1)}
but for which $\tau(n)$ grows only polynomially fast.}.
Together with the second requirement, {\bf (S1)} allows one 
to make a replacement argument and approximate
each term in the decomposition of $X_n$ in~\eqref{sumpiece}, after centering and rescaling,
by an independent copy of the stable law $\rS_s$.

We then move in Section~\ref{Poisson}
to the second part in which we show~\eqref{Xnstablelim2}.
Under {\bf (S2)}
increments grow slowly,
so that a growing number of the terms in the decomposition~\eqref{sumpiece}
contribute to the distribution of $X_n$, and the
replacement argument used under assumption {\bf (S1)} no longer works.
In this case we show that the joint distribution of the terms in~\eqref{sumpiece} converge,
after proper centering and scaling,
to that of the atoms of a certain non-homogeneous Poisson process. 
This proof is similar to standard proofs of stable limit laws
for sums of i.i.d.\ random variables. 
Indeed, in the case where $g(\infty) = 0$
(where the limiting distribution is $\rS_s$)
this non-homogeneous Poisson process is exactly
the same as what one would get if one were considering
i.i.d.\ sums of random variables in the domain of attraction of $\rS_s$. 
When $g(\infty) >0$ the Poisson process is slightly
different and leads to the presence of a tempered
stable component in the limiting distribution. 
 
This completes the proof of  the general Theorem \ref{GeneralConditions-S},
and we conclude in Section \ref{polycases} by showing how to use this to
derive the stable and tempered stable limits in~\eqref{Xnstableliml}-\eqref{Xnstablelim}
for polynomial cooling maps. 

\subsection{Fast enough cooling maps: proof of stable limits under {\bf (S1)}}
\label{FastCool}
In this Section we prove~\eqref{Xnstablelim1} under {\bf (S1)}. 
For ease of notation, given two random variables $X,Y$,
we will write $X \overset{\text{Law}}{=}Y$ when the two random variables
have the same distribution, i.e., when
$\bb{E}\crt{\exp(iuX)} = \bb{E}\crt{\exp(iuY)}$ for all $u \in \bb{R}$.

We start with two preliminary observations: first, we note that {\bf (S1)}
implies in particular that \eqref{cag} is satisfied with $g(x)\equiv 0$,
and second, as expressed in Lemma~\ref{stablesubseq} below,
we show that while proving the claim the boundary term can be neglected.

For the first observation, 
if the first condition in {\bf(S1)} holds, then there is $c>0$ such that $\tau(n) \geq c n^s$.
Furthermore, it can also be shown that for all $x>0$, $\theta>0$:
\begin{equation}\label{thetascale}
  \lim_{n\to\infty} \frac{\sum_{k = 1}^n T_k \ind{\tau(n)^\theta<T_k<x\tau(n)^{1/s}}}{\tau(n)} = 0.
\end{equation}
This two conditions imply that we may choose $\theta$ small enough so that $s(1 - \theta )>1$ and so
\begin{align*}
\limsup_{n\to\infty}  \frac{ \sum_{k=1}^{n} T_k \ind{T_k < x \tau(n)^{1/s}} }{\tau(n)}
  &\leq \limsup_{n\to\infty} \frac{ \sum_{k=1}^{n} T_k \ind{T_k < \tau(n)^{\theta}} }{\tau(n)} \\
  &\leq \frac{n \tau(n)^\theta}{\tau(n)}  =\lim_n \frac{n}{\tau(n)^{1 - \theta}}
  \leq \lim_{n\to\infty}\frac{n}{n^{s(1 - \theta)}} = 0 .
\end{align*}
Yet, as this is not required in the proof\footnote{Our proof of \eqref{Xnstablelim1}
uses only {\bf(S1)} and doesn't use \eqref{cag}.
We include the observation that {\bf(S1)} implies \eqref{cag}
with $g\equiv 0$ to help show the consistency of the two parts
of Theorem \ref{GeneralConditions-S}},
we leave to the interested reader to check~\eqref{thetascale}.

The second preliminary observation is captured in the next lemma.
The idea behind it is that as soon as the last term
in the decomposition~\eqref{sumpiece} is large enough to make
a non-negligible contribution to the distribution of $X_n$,
then the distribution of this last term
can be combined with the other terms to give the limit stable law in \eqref{StableLeftSkew}.

\begin{lemma}[{\bf Negligible boundary for pure stable limit}]\label{stablesubseq}
Let $\rS_s$ be the stable random variable which arises
as the limiting distribution of RWRE in~\eqref{StableLeftSkew}. 
If the cooling sequence $\{T_k\}_{k\geq 1}$ is such that 
\begin{equation}\label{taustable}
 \frac{X_{\tau(n)} - \EE[X_{\tau(n)}]}{ \tau(n)^{1/s} } \limdist \rS_s,  
\end{equation}
then it follows that~\eqref{Xnstablelim} holds also. 
\end{lemma}
\begin{proof}
Denote by $\chi(u) := E\left[ e^{i u \rS_s} \right]$
the characteristic function of the stable random variable $\rS_s$ in~\eqref{StableLeftSkew}, and let
\begin{equation}\label{chfc}
 \phi_n(u) = \EE\left[\exp\left\{iu \frac{Z_n - \EE[Z_n]}{n^{1/s}} \right\} \right]
\quad\text{and}\quad
 \psi_n(u) = \EE\left[\exp\left\{ iu \frac{X_{\tau(n)} - \EE[X_{\tau(n)}]}{(\tau(n))^{1/s}} \right\} \right]
\end{equation}
be the characteristic functions of $Z_n$ and $X_{\tau(n)}$
after appropriate centering and scaling.
Also let $q_n := \tfrac{\bar{T}_n}{n}$.
With this notation and using the decomposition in~\eqref{sumpiece},
the characteristic function of $\frac{X_n-\EE[X_n]}{n^{1/s}}$ can be expressed as
\begin{align}
& \EE\left[ \exp\left\{ iu \frac{X_n-\EE[X_n]}{n^{1/s}} \right\} \right]
= \psi_{\ell_n}\left( \left( 1- q_n\right)^{1/s}u \right)
\phi_{\bar{T}_n}\left( \left( q_n \right)^{1/s} u \right) \nonumber\\
&\quad= \chi\left( \left( 1-q_n \right)^{1/s} u \right)
\chi\left( \left( q_n\right)^{1/s} u \right) \nonumber\\
&\quad\qquad + \left\{ \psi_{\ell_n}\left( \left( 1-q_n \right)^{1/s} u \right)
- \chi\left( \left( 1-q_n \right)^{1/s} u \right) \right\}
\phi_{\bar{T}_n}\left( \left( q_n \right)^{1/s} u  \right) \label{psihatL} \\
&\quad\qquad + \chi\left( \left( 1-q_n \right)^{1/s} u \right)
\left\{ \phi_{\bar{T}_n}\left( \left( q_n \right)^{1/s} u \right)
- \chi\left( \left( q_n \right)^{1/s} u \right) \right\} \label{phihatL}. 
\end{align}

It follows from the explicit formula for $\chi(u)$ in~\eqref{stablecf} that 
\begin{equation}\label{chprod}
 \chi\left(\left( 1- t \right)^{1/s} u \right) \chi\left( t^{1/s} u \right) = \chi(u), \quad \forall t \in [0,1].
\end{equation}
To finish the proof, we show that for any fixed $u \in \R$
\eqref{psihatL} and~\eqref{phihatL} vanish as $n\to\infty$. 
This follows from the fact that convergence in distribution
implies uniform convergence of characteristic functions
on compact sets (see for instance Theorem 15 in Chapter 14 of~\cite{FG97}).
Indeed, since $\ell_n \to \infty$ as $n\to\infty$,
\eqref{taustable} implies that~\eqref{psihatL}
vanishes as $n\to\infty$. To control~\eqref{phihatL}, note that
for any fixed $m<\infty$
\begin{align*}
  & \left| \phi_{\bar{T}_n}\left( \left( q_n \right)^{1/s} u \right)
  - \chi\left( \left( q_n \right)^{1/s} u \right) \right| \\&
  \leq \max_{k\leq m} \left| \phi_k\left(\left( \tfrac{k}{n}\right)^{1/s} u \right)
  - \chi\left(\left( \tfrac{k}{n}\right)^{1/s} u \right) \right|
+  \max_{k>m} \sup_{|v| \leq |u|} \left| \phi_k(v) - \chi(v) \right|. 
\end{align*}
The first term on the right vanishes as $n\to\infty$ for any fixed $m$
since all characteristic functions are continuous and equal to one at the origin,
while the second term can be made arbitrarily small by choosing $m$ sufficiently large.
This completes the proof of the lemma. 
\end{proof}

\paragraph{Proof of stable limits under condition {\bf(S1)}.}
In view of the previous lemma, it remains to show~\eqref{taustable}.
We may and do consider the space $(\Omega, \mc{F}, \bb{P})$
to be rich enough to contain an extra infinite sequence of
uniform random variables $\bar{U}: = \{U^{(k)}, k \in \bb{N}\}$
with respect to which we will define auxiliary random variables. 

Given a random variable $X$, let $F^{-1}_X(a):= \inf\{x \colon \bb{P}(X<x)\geq a\}$
represent its generalized inverse function.
Let
\begin{equation}\label{invW}	
  \rS_s^{(k)} := F_{\rS_s}^{-1}(U^{(k)}) \quad \text{and let} \quad
  \Psi^{(k)}_m := F_{\frac{Z_{m} - \bb{E}[Z_{m}]}{m^{1/s}}}^{-1}(U^{(k)}).
\end{equation} 
The limiting distribution in~\eqref{StableLeftSkew}
together with \eqref{meanzn}
implies $\lim_{m\to\infty}\Psi^{(k)}_m = \rS_s^{(k)} $, almost surely, for any $k$. 
Then, the uniform moment bounds in~\eqref{LpMom}
imply that this convergence holds in $L^p$ for any $p \in (0,s)$.
That is, if for each $m \in \bb{N}$ we define the error term 
$ \err^{(k)}_{m}: =  \Psi^{(k)}_{m}- \rS_s^{(k)}$, then 
\begin{equation}\label{pmom0}
\lim_{m \to \infty} \sup_k  \bb{E}\big[\,\vert\,\err^{(k)}_{m}\,\vert^p\, \big]= 0,
\qquad \forall p \in (0,s). 
\end{equation}
In particular, for $p=1$ we have that
for any $\gep>0$ there is an $m_0= m_0(\gep)$ such that 
\begin{equation}\label{mom1gep}
	m>m_0 \Rightarrow \sup_k  \bb{E}\big[\,\vert\,\err^{(k)}_{m}\,\vert\, \big] <\gep.
\end{equation}
Now, by the definition of~\eqref{invW},
if we set $\alpha_{k,n}:= \prt{\frac{T_k}{\tau(n)}}^{1/s}$, then it follows that
\begin{equation}\label{remain}
\frac{X_{\tau(n)} - \EE[X_{\tau(n)}]}{\tau(n)^{1/s}}
= \sum_{k=1}^{n} \frac{Z_{T_k}^{(k)} - \bb{E}[Z_{T_k}^{(k)}]}{\tau(n)^{1/s}}
\overset{\text{Law}}{=}
\sum_{k=1}^{n} \alpha_{k,n}\rS_s^{(k)} + \sum_{k=1}^{n} \alpha_{k,n} \err^{(k)}_{T_k}.
\end{equation} 
Since the law of $\rS_s$ satisifes~\eqref{stablecf},
we obtain that
$\sum_{k=1}^{n} \alpha_{k,n}\rS_s^{(k)}
\overset{\text{Law}}{=}
\left( \sum_{k=1}^n \alpha_{k,n}^s \right)^{1/s} \rS_s = \rS_s.$
To complete the proof, we show that $\sum_{k=1}^{n} \alpha_{k,n} \err^{(k)}_{T_k}$ converges to zero in $L^1$ (and therefore also in distribution). 
For any $\epsilon>0$ fixed and $m_0$ as in~\eqref{mom1gep} there is $C>0$ for which
\begin{equation}\label{ee:02}
\begin{aligned}
\sum_{k=1}^n \alpha_{k,n} \EE\left[ \left|  \err^{(k)}_{T_k} \right| \right] 
&\leq \sum_{k=1}^n \alpha_{k,n} \EE\left[ \left|  \err^{(k)}_{T_k} \right| \right] \ind{T_k \leq m_0}
+ \sum_{k=1}^n \alpha_{k,n} \epsilon \\
&\leq C \left( \sum_{k=1}^n \alpha_{k,n} \ind{T_k \leq m_0} \right) + C \epsilon,  
\end{aligned}
\end{equation}
where the last inequality follows from the first condition in {\bf(S1)} 
and the fact that \eqref{pmom0} implies $\sup_{m,k} \EE\left[ \left|  \err^{(k)}_{m} \right| \right] < \infty$.
Furthermore, the second condition of {\bf(S1)}  implies that 
$\lim_{n\to\infty} \bb{E}\left[ \sum_{k=1}^{n} \alpha_{k,n} \ind{T_k \leq m_0} \right] = 0$.
Since $\epsilon>0$ was arbitrary,
this completes the proof \eqref{taustable} and thus of the stable limit in \eqref{Xnstablelim1} under condition {\bf(S1)}. 
\qed

\subsection{Poisson processes: tempered stable and stable limits under {\bf (S2)}}\label{Poisson} 
Analogously to Lemma~\ref{stablesubseq},
we start with a lemma which allows us to ignore boundary terms.
\begin{lemma}[{\bf Negligible boundary for $n^{1/s}$ scaling}]\label{lem:boundary}
 Suppose that the cooling sequence is such that
\begin{equation}\label{nodominant-2}
 \lim_{n\to\infty} \max_{k\leq n} \frac{T_k}{\tau(n)} = 0. 
\end{equation}
If there exists a random variable $Z$ such that $\frac{X_{\tau(n)} - \EE[X_{\tau(n)}]}{\tau(n)^{1/s}} \limdist Z$, 
then it is also true that  $\frac{X_n - \EE[X_n]}{n^{1/s}} \limdist Z$.
\end{lemma}
\begin{proof}
Using the decomposition in \eqref{sumpiece} we can write
\begin{align*}
 \frac{X_n-\EE[X_n]}{n^{1/s}}
= \left( 1- \frac{\bar{T}_n}{n}\right)^{1/s} \frac{X_{\tau(\ell_n)}-\EE[X_{\tau(\ell_n)}]}{\tau(\ell_n)^{1/s}} + \left( \frac{\bar{T}_n}{n} \right)^{1/s} \frac{Z^{(\ell_n+1)}_{\bar{T}_n}-\EE[Z_{\bar{T}_n}] }{\bar{T}_n^{1/s}}. 
\end{align*}
To conclude the proof, it suffices to show that $\lim_{n\to 0} \frac{\bar{T}_n}{n} = 0$.
Indeed, since \eqref{StableLeftSkew} and~\eqref{meanzn}
together imply that $\frac{Z_n - \EE[Z_n]}{n^{1/s}} \limdist \rS_s$, the sequence
$\frac{Z^{(\ell_n+1)}_{\bar{T}_n}-\EE[Z_{\bar{T}_n}] }{\bar{T}_n^{1/s}}$ is tight. 
This implies that the second term on the right converges to 0 in probability,
while the assumptions of the lemma imply that
the first term on the right converges in distribution to $Z$. 

Since $\tau(\ell_n) \leq n < \tau(\ell_n+1)$ and $\bar{T}_n = n-\tau(\ell_n) < T_{\ell_n+1}$, we have that
$\frac{\bar{T}_n}{n} \leq \frac{T_{\ell_n+1}}{\tau(\ell_n)}$. Thus it is enough to show that  $\lim_{\ell\to\infty} \frac{T_{\ell+1}}{\tau(\ell)} = 0$.  
For any $\epsilon > 0$, \eqref{nodominant-2} implies that for $\ell$ sufficiently large $T_{\ell+1} < \epsilon \tau(\ell+1) = \epsilon \tau(\ell) + \epsilon T_{\ell+1}$, 
and thus $\limsup_{\ell \to \infty} \frac{T_{\ell+1}}{\tau(\ell)} \leq \frac{\epsilon}{1-\epsilon}$. Since $\epsilon>0$ was arbitrary, this completes the proof of the lemma. 
\end{proof}

We can now restrict the analysis to the subsequence $\tau(n)$. For convenience of notation, let 
\begin{equation}\label{xikn}
 \xi_{k,n} = \frac{Z^{(k)}_{T_k} - \EE[Z_{T_k}] }{\tau(n)^{1/s}},
\end{equation}
so that $\frac{X_{\tau(n)} - \EE[X_{\tau(n)}]}{\tau(n)^{1/s}} = \sum_{k=1}^n \xi_{k,n}$.
By Lemma~\ref{lem:boundary}, to prove~\eqref{Xnstablelim2} we need to show that
$\sum_{k=1}^n \xi_{k,n} \limdist \rW_{\lambda_g} + (1-g(\infty))^{1/s} \rS_s$,
where in a slight abuse of notation here and below we interpret $\rW_{\lambda_g} \equiv 0$ when $g\equiv 0$. 
The proof is divided in four steps.
First we show that for any $t>0$,
the truncated point process constructed from $\delta_{\xi_{k,n}}\ind{\xi_{k,n}\leq -t}$
converges in distribution to a certain Poisson point process $N^{(g)}_{t}$.
By the continuous mapping theorem, this implies that
$\sum_{k=1}^n \xi_{k,n} \ind{\xi_{k,n} \leq -t}$
converges in distribution to a functional $\Psi\left(N^{(g)}_{t}\right)$
of the point process $N^{(g)}_{t}$, and in step 2 we prove that the 
corresponding means also converge as $n\to\infty$. 
Step 3 controls the error introduced by omitting the terms
$\xi_{k,n}$ in the sum with $\xi_{k,n} > -t$. 
Finally, in step 4 we combine the previous results
to show first of all that
the limiting distribution of the RWCRE is
$\lim_{t\to 0}  \Psi\left(N^{(g)}_{t}\right) - E\left[\Psi\left(N^{(g)}_{t}\right)\right] $
and that this limit has the same distribution as $\rW_{\lambda_g} + (1-g(\infty))^{1/s} \rS_s$. 

\paragraph{Step 1. Convergence of  $t$-truncated processes.}
For any $t > 0$ and $n\geq 1$, let $N_{n,t}$ be the point process defined by
\begin{equation}\label{Nnedef}
  N_{n,t} := \sum_{k=1}^n \delta_{\xi_{k,n}} \ind{\xi_{k,n} \leq - t  }.
\end{equation}
We will show in this step that 
\begin{equation}\label{Radonvague}
  N_{n,t} \Longrightarrow N^{(g)}_{t} , \qquad \forall t>0, 
\end{equation}
i.e., that  $N_{n,t}$ converges in distribution, as $n\to\infty$,
(on the space of Radon point processes on $[-\infty,0)$ equipped with the vague topology)
to a  non-homogeneous Poisson point process $N^{(g)}_{t}$ with intensity 
$\hlg(x) \ind{x\leq -t}$, where $\hlg(x) = \lambda_g(x) + (1-g(\infty))\lambda_0(x) $ with $\lambda_g(x)$ is defined as in \eqref{lgdef} and $\lambda_0(x) = \tc \vs s |x|^{-s-1}$.

By~\cite[Theorem 11.2.V]{DV08},
since for each $n$ the random variables $\{\xi_{k,n} \}_{k\leq n}$ are independent
 to prove~\eqref{Radonvague} it is enough to check that 
\begin{enumerate}
 \item $\lim_{n\to\infty} \max_{k\leq n} \PP( \xi_{k,n} \leq -t) = 0$, $\forall t>0$, 
 \item and if $\mu_n$ is the measure on $(-\infty,0)$ defined by 
 $\mu_n(A) = \EE\left[ \sum_{k=1}^n \delta_{\xi_{k,n}} \ind{\xi_{k,n} \in A } \right]$, 
then $\mu_n(dx)$ converges weakly to the measure 
$\hlg(x) \, dx$.
\end{enumerate}
Since $\xi_{k,n} = \left( \frac{T_k}{\tau(n)} \right)^{1/s} \frac{Z_{T_k}^{(k)} - \EE[Z_{T_k}]}{T_k^{1/s}} $,
condition 1 above follows from the assumption in {\bf(S2)}
and the fact that the family $\{ (Z_{T_k}^{(k)} - \EE[Z_{T_k}])/T_k^{1/s} \}_{k\geq 1}$
is tight. 

To prove the weak convergence of
$\mu_n(dx)$ to $\hlg(x) \, dx$,
we prove for all $t>0$ that 
\begin{equation}\label{munconv}
\lim_{n\to\infty} \mu_n((-\infty,-t])
=\int_{-\infty}^{-t}\hlg(x) \, dx
= \int_{-\infty}^{-t} \lambda_g(x) \, dx + (1-g(\infty)) \tc \vs t^{-s} ,
\end{equation}
where the last equality follows from the definition of $\hlg(x)$. 
We first notice that for any $t>0$ 
\begin{equation}\label{e:14}
\mu_n((-\infty,-t]) 
= \sum_{k=1}^n \PP( \xi_{k,n} \leq -t  )
= \sum_{k=1}^n \PP\left( Z_{T_k} - \EE[Z_{T_k}] \leq -t \tau(n)^{1/s}  \right).
\end{equation}
In order to control the sum in the right-hand side of \eqref{e:14}, we
rely on the estimates in Corollary~\ref{cor:prectail-c}.
However, these estimates will only apply if
\begin{equation}\label{e:7}
  T_k^{1/s}(\log T_k)^4 \leq t \tau(n)^{1/s} \leq T_k \vs - T_k^{1/s} (\log T_k).
\end{equation}
Since condition {\bf (S2)} implies that
the first of these inequalities holds for all $k\leq n$ when $n$ is large enough,
it will be important to consider when the second inequality holds also.
Therefore, it is natural to define the set
\begin{equation}\label{Antdef}
 A_{n,t} = \left\{k \leq n:  t \tau(n)^{1/s} \leq T_k \vs - T_k^{1/s}(\log T_k) \right\}.
\end{equation}
We collect in the next technical lemma some properties
of this set which will be used in the sequel.
In particular, as expressed in \eqref{main} below,
it turns out that the non-vanishing contribution in
the limit of the sum in~\eqref{e:14}
comes precisely from the terms in this set $A_{n,t}$.

\begin{lemma}[{\bf $A_{n,t}$ and non-vanishing contribution of negative points }]\label{AntContri}\text{}\\
Let $A_{n,t}$ be as in~\eqref{Antdef}.
If conditions~\eqref{cag} and {\bf (S2)} hold,
then for every $t>0$
\begin{equation}\label{main}
\lim_{n\to\infty} \sum_{k\leq n, \, k\notin A_{n,t}} \PP(\xi_{k,n} \leq -t) = 0,
\end{equation}
and
\begin{equation}\label{mainE}
\lim_{n\to\infty} \sum_{k\leq n, \, k\notin A_{n,t}}\EE\left[  \xi_{k,n} \ind{\xi_{k,n}\leq -t} \right]=0.
\end{equation}

Further, for any continuous function $f(x)$ on $[0,\infty)$ with $\lim_{x \to \infty} \frac{f(x)}{x}  = L<\infty$,
\begin{equation}\label{e:24}
  \lim_{n\to\infty} \frac{1}{\tau(n)^{1-\frac{1}{s}}} \sum_{k \in A_{n,t}} f\left( \frac{T_k}{\tau(n)^{1/s}} \right)
  = \int_{t/v}^\infty \frac{f(x)}{x} \, g(dx) + L (1 - g(\infty)), \quad \forall t>0.
\end{equation}
\end{lemma}

The proof of this lemma is postponed to Appendix~\ref{Acp}.
We now conclude step 1. By~\eqref{e:14} and~\eqref{main},
to prove \eqref{munconv} it suffices to consider
the sum over $k$'s in $A_{n,t}$.
For the latter, we can use Corollary~\ref{cor:prectail-c} which implies that 
\begin{equation}\label{e:30}
 \lim_{n\to\infty} \sum_{k \in A_{n,t}} \PP(\xi_{k,n} \leq -t)
= \lim_{n\to\infty} \tc \sum_{k \in A_{n,t}}  \frac{T_k \vs - t \tau(n)^{1/s}}{\tau(n) t^s},
\end{equation}
as long as we can prove the limit on the right exists. 
To this end, we re-write the sum and then apply \eqref{e:24} to conclude that
\begin{align}
 \sum_{k \in A_{n,t}} \tc \frac{T_k \vs - t \tau(n)^{1/s}}{\tau(n) t^s}
&= \tc \frac{t^{-s}}{\tau(n)^{1-\frac{1}{s}}} \sum_{k \in A_{n,t}} \left( \frac{T_k \vs}{\tau(n)^{1/s}} - t  \right) \nonumber \\
&=  \tc t^{-s} \int_{t/\vs}^\infty \left(\vs-\frac{t}{x} \right) \, g(dx) +  (1-g(\infty))\tc\vs t^{-s} . \label{Int3}
\end{align}

In view of~\eqref{main} and~\eqref{Int3},
it remains only to check that the integrals in the right-hand side of~\eqref{munconv}
and \eqref{Int3} are equal.
This follows by the definition of $\lambda_g(x)$  in \eqref{lgdef} from which
we have 
\begin{align*}
 \int_{-\infty}^{-t} \lambda_g(x) \, dx 
= \int_t^\infty \lambda_g(-z) \, dz 
&= \int_t^\infty \tc z^{-s} \int_{z/\vs}^\infty \left( \frac{\vs s}{z} - \frac{s-1}{x} \right) \, g(dx) \, dz \\
&= \tc \int_{t/\vs}^\infty \int_t^{x\vs} \left( \vs s z^{-s-1} - \frac{s-1}{x}z^{-s} \right) \, dz \, g(dx) \\
&= \tc \int_{t/\vs}^\infty \left( \vs t^{-s} - \frac{t^{-s+1}}{x} \right) \, g(dx), 
\end{align*}
and this last expression is equal to the integral in the right side of~\eqref{Int3}.
This completes the proof of~\eqref{e:14} and therefore of Step 1. 

\paragraph{Step 2. Convergence of  $t$-truncated means.}
We first notice that the weak convergence in \eqref{Radonvague} shown in Step 1,
implies in particular that for any $t>0$,
\begin{equation}\label{sumxikn-main}
 \sum_{k=1}^n \xi_{k,n} \ind{\xi_{k,n}\leq -t} = \Psi(N_{n,t}) \underset{n\to\infty}{\Longrightarrow}
 \Psi(N^{(g)}_t), \quad \text{where } \Psi(\nu) = \int x \, \nu(dx), 
\end{equation}
since the functional $\Psi$ is continuous
with respect to the vague topology
on the set of point processes with no atoms at $t$. 
In this step we show that the means in \eqref{sumxikn-main} also converge.
That is, we show that
\begin{equation}\label{cmeanconvergence}
   \lim_{n\to \infty} \EE\left[\sum_{k=1}^n  \xi_{k,n} \ind{\xi_{k,n}\leq -t} \right] = E[\Psi(N_t^{(g)})]. 
\end{equation}

Once again, thanks to \eqref{mainE} in Lemma~\ref{AntContri},
we will restrict the sum in the right-hand side of~\eqref{cmeanconvergence}
only to the indexes in $A_{n,t}$. For the latter, we first re-write 
\begin{align}
 \EE\left[\sum_{k \in A_{n,t}}  \xi_{k,n} \ind{\xi_{k,n}\leq -t} \right]
&= -t \sum_{k \in A_{n,t}} \PP(Z_{T_k} - \EE[Z_{T_k}] \leq -t \tau(n)^{1/s} ) \label{tfm1} \\
 &\qquad - \sum_{k \in A_{n,t}} \int_t^\infty \PP(Z_{T_k} - \EE[Z_{T_k}] \leq -u \tau(n)^{1/s} ) \, du.
 \label{tfm2} 
\end{align}
The asymptotics of~\eqref{tfm1} follow from the same analysis leading to~\eqref{Int3} above. That is, 
\begin{equation}
  \lim_{n\to\infty} \eqref{tfm1} 
  =  - \tc t^{1-s} \int_{t/ \vs}^\infty \left( \vs - \frac{t}{x} \right) \, g(dx)
- (1-g(\infty))\tc\vs t^{1-s}.\label{tfm1-lim}
\end{equation}
For the sum in~\eqref{tfm2},
let $\gamma_{k,n} := \frac{T_k \vs - T_k^{1/s}(\log T_k)}{\tau(n)^{1/s}}$.
For $n$ large enough the probabilities inside the integrals
can be approximated by Corollary~\ref{cor:prectail-c} for $u \in [t,\gamma_{k,n}]$. 
That is, assuming we can show the limit on the right exists, we have 
\begin{align*}
& \lim_{n\to\infty} \sum_{k \in A_{n,t}} \int_t^{\gamma_{k,n}} \PP(Z_{T_k} - \EE[Z_{T_k}] \leq -u\tau(n)^{1/s} ) \, du \\
&\qquad = \lim_{n\to\infty} \sum_{k \in A_{n,t}} \int_t^{\gamma_{k,n}} \tc(T_k \vs - u\tau(n)^{1/s}) u^{-s} \tau(n)^{-1} \, du\\
&\qquad = \lim_{n\to\infty} \sum_{k \in A_{n,t}} \left\{  \frac{\tc \vs T_k}{\tau(n)} \int_t^{\gamma_{k,n}}  u^{-s} \, du -  \frac{\tc}{\tau(n)^{1-\frac{1}{s}}} \int_t^{\gamma_{k,n}}  u^{1-s} \, du \right\} \\
&\qquad = \lim_{n\to\infty} \sum_{k \in A_{n,t}} \left\{  \frac{\tc \vs T_k}{\tau(n)} \left( \frac{t^{1-s}-\gamma_{k,n}^{1-s}}{s-1} \right) -  \frac{\tc}{\tau(n)^{1-\frac{1}{s}}}\left( \frac{\gamma_{k,n}^{2-s}-t^{2-s}}{2-s} \right) \right\} \\
&\qquad =  \lim_{n\to\infty} \frac{\tc}{\tau(n)^{1-\frac{1}{s}}} \sum_{k \in A_{n,t}} \left\{  \frac{T_k }{\tau(n)^{1/s}} \frac{\vs t^{1-s}}{s-1} + \frac{t^{2-s}}{2-s} 
- \frac{T_k}{\tau(n)^{1/s}} \frac{\vs}{s-1} \gamma_{k,n}^{1-s} - \frac{1}{2-s} \gamma_{k,n}^{2-s} \right\}.
\end{align*}
By {\bf (S2)},
we can replace $\gamma_{k,n}$ with $\frac{T_k \vs}{\tau(n)^{1/s}}$,
and apply~\eqref{e:24} to conclude that 
\begin{align}
& \lim_{n\to\infty} \sum_{k \in A_{n,t}} \int_t^{\gamma_{k,n}} \PP(Z_{T_k} - \EE[Z_{T_k}] \leq -u\tau(n)^{1/s} ) \, du \nonumber \\
&\qquad = \tc \int_{t/\vs}^\infty \left( \frac{\vs t^{1-s}}{s-1} + \frac{t^{2-s}}{(2-s)x} - \frac{\vs^{2-s}x^{1-s} }{(s-1)(2-s)} \right) \, g(dx)
 + (1-g(\infty)) \tc \frac{\vs t^{1-s}}{s-1} . \label{tfm2-lim}
\end{align}

This computes the main asymptotic value of the terms in~\eqref{tfm2},
but we still need to control the sum over the integrals in~\eqref{tfm2}
for $t \geq \gamma_{k,n}$.
To this end, first note that the probabilities inside the integrals in~\eqref{tfm2}
are decreasing in $u$ and are zero for $u\geq 2T_k/\tau(n)^{1/s}$.
Thus, we obtain the simple upper bound
\[
 \int_{\gamma_{k,n}}^\infty \PP(Z_{T_k} - \EE[Z_{T_k}] \leq -u \tau(n)^{1/s} ) \, du  \leq  \frac{2T_k}{\tau(n)^{1/s}} \PP\left(Z_{T_k} - \EE[Z_{T_k}] \leq -\gamma_{k,n}\tau(n)^{1/s} \right). 
\]
By Corollary~\ref{cor:prectail-c},
for any $\nu \in (\frac{1}{s},1)$
there exists a constant $C<\infty$
so that the probability on the right above is
bounded above by $C T_k^{-s+\nu}$. 
Therefore, 
\begin{align}
 \limsup_{n\to\infty}  \sum_{k \in A_{n,t}} \int_{\gamma_{k,n}}^\infty &\PP\left(Z_{T_k} - \EE[Z_{T_k}] \leq -u \tau(n)^{1/s} \right) \, du \nonumber \\
& \leq \lim_{n\to\infty} \frac{2C}{\tau(n)^{1-\frac{\nu}{s}}} \sum_{k \in A_{n,t}} \left(\frac{T_k}{\tau(n)^{1/s}} \right)^{1-s+\nu}
= 0, \label{tfm2b-lim}
\end{align}
where the last limit is zero by~\eqref{e:24} and the fact that we chose $\nu<1$.

Applying~\eqref{tfm1-lim},~\eqref{tfm2-lim}, and~\eqref{tfm2b-lim} to~\eqref{tfm1}-\eqref{tfm2}, we obtain that 
\begin{align}
& \lim_{n\to\infty} \EE\left[\sum_{k \in A_{n,t}}  \xi_{k,n} \ind{\xi_{k,n}\leq -t} \right] \nonumber \\
&= -\tc \int_{t/\vs}^\infty \left( \frac{\vs s t^{1-s}}{s-1} + \frac{(s-1)t^{2-s}}{(2-s)x} - \frac{\vs^{2-s}x^{1-s} }{(s-1)(2-s)}   \right) \, g(dx)
 - (1-g(\infty))\frac{\tc \vs s t^{1-s}}{s-1}. \label{tfm-lim}
\end{align}
It remains to show that this right-hand side in \eqref{tfm-lim} equals $E[\Psi(N_t^{(g)})]$. 
To this end, Campbell's Theorem~\cite[Section 3.2]{Kin93}
implies that 
$E[\Psi(N_t^{(g)})] = \int_{-\infty}^{-t} x \hlg(x) \, dx = \int_{-\infty}^{-t} x \lambda_g(x) \, dx + (1-g(\infty)) \int_{-\infty}^{-t} x \lambda_0(x) \, dx$.  
Recalling that $\lambda_0(x) = \tc \vs s |x|^{-s-1}$, a simple calculation yields that the second term on the right equals the second term in \eqref{tfm-lim}. 
On the other hand, the formula for $\lambda_g$ in \eqref{lgdef} yields that $\int_{-\infty}^{-t} x \lambda_g(x) \, dx$ equals
\begin{align*}
- \int_t^\infty u \, \lambda_g(-u) \, du
&= - \tc \int_t^\infty u^{1-s} \int_{u/\vs}^\infty \left( \frac{\vs s}{u} - \frac{s-1}{x} \right) \, g(dx) \, du \\
&= -\tc \int_{t/\vs}^\infty \int_t^{x\vs} \left( \vs s u^{-s} - \frac{s-1}{x} u^{1-s}  \right) \, du \, g(dx) \\
&= -\tc  \int_{t/\vs}^\infty \left( \frac{\vs s t^{1-s}}{s-1} + \frac{(s-1)t^{2-s}}{(2-s)x} - \frac{\vs^{2-s}x^{1-s}}{(s-1)(2-s)} \right) \, g(dx), 
\end{align*}
which matches the first term on the right-hand side of \eqref{tfm-lim}.
This completes the proof of~\eqref{cmeanconvergence} and thus finishes Step 2.

\paragraph{Step 3. Negligible contribution from small points.}
Next, we will show that the contribution of the sum of
the $\xi_{k,n}$ with $\xi_{k,n} > -t$ is essentially negligible
if $n$ is large and $t$ is small.
That is, we will show that 
\begin{equation}\label{e:16}
  \lim_{t\to 0} \limsup_{n\to\infty}
  \PP\left( \left| \sum_{k = 1}^n \left( \xi_{k,n} \ind{\xi_{k,n} > -t}
  - \EE[ \xi_{k,n} \ind{\xi_{k,n} > -t}] \right)
  \right| > \delta \right) = 0, \quad \forall \delta > 0. 
\end{equation}
Since the random variables inside the sum are independent, to prove the above it is enough to show that 
$\lim_{t \to 0} \limsup_{n\to\infty} \sum_{k=1}^n  \Var\left(\xi_{k,n} \ind{\xi_{k,n} > -t} \right) = 0$. 
For this, note first of all that  
\begin{equation}\label{tvarsplit}
 \Var\left(\xi_{k,n} \ind{\xi_{k,n} > -t} \right) \leq \EE[ \xi_{k,n}^2  \ind{\xi_{k,n} > -t}] 
= \EE[ \xi_{k,n}^2  \ind{\xi_{k,n} \in (-t,0)}] +\EE[ \xi_{k,n}^2  \ind{\xi_{k,n} > 0}].
\end{equation}
We can bound the first expectation  above with~\eqref{varzn}
if $k \notin A_{n,t}$ or~\eqref{truncVar}
if $k \in A_{n,t}$, and we can bound the second expectation with Corollary~\ref{cor:L2}.
Therefore, there is a  $C>0$, for which
\begin{equation}\label{e:17}
\sum_{k=1}^n \Var\left(\xi_{k,n} \ind{\xi_{k,n} > -t} \right)
\leq \sum_{k \notin A_{n,t}} \frac{C T_k^{3-s}}{\tau(n)^{2/s}} + \sum_{k \in A_{n,t}} \frac{C T_k t^{2-s} }{\tau(n)}
+ \sum_{k=1}^n \frac{C T_k^{2/s}}{\tau(n)^{2/s}}. 
\end{equation}
Since $s<2$ and $\sum_{k \in A_{n,t}} T_k \leq \sum_{k=1}^n T_k = \tau(n)$ the second term on the right-hand side above can be made arbitrarily small if we take $t \to 0$ and 
the third term can be bounded by $C \max_{k\leq n} \left( T_k/\tau(n) \right)^{\frac{2}{s}-1}$ which vanishes as $n\to \infty$ by {\bf (S2)}. 
Finally for the first term,
if $T_k > \tfrac{2t}{ \vs}\tau(n)^{1/s}$ then  for $n$ sufficiently large $T_k - T^{1/s}_k (\log T_k)> t\tau(n)$  and so  $k \in A_{n,t}$. Thus
\begin{equation}\label{e:37}
  \limsup_{n \to \infty} \sum_{k \notin A_{n,t}} \frac{T_k^{3-s}}{\tau(n)^{2/s}}
  = \limsup_{n \to \infty}\sum_{k \notin A_{n,t}} \frac{T_k}{\tau(n)} \left( \frac{T_k}{\tau(n)^{1/s}} \right)^{2-s}
  \leq (2t/\vs)^{2-s}. 
\end{equation}
Now we take $t \to 0$ and  complete the proof of Step 3.

\paragraph{Step 4. Convergence of the process.} 
Finally, we will show how the above three steps imply that $\sum_{k=1}^n \xi_{k,n} \limdist \rW_{\lambda_g}$. 
First of all, for any $t>0$ we can write $\sum_{k=1}^n \xi_{k,n} = W_{n,t} + \mathfrak{E}_{n,t}$, where 
\begin{equation}\label{WneRne}
 \begin{split}
W_{n,t} &= \sum_{k=1}^n \left( \xi_{k,n} \ind{\xi_{k,n} \leq -t}
- \EE[\xi_{k,n} \ind{\xi_{k,n} \leq -t}] \right)
\\ 
\text{and}\qquad 
\mathfrak{E}_{n,t} &= \sum_{k=1}^n \left( \xi_{k,n} \ind{\xi_{k,n} > -t}
- \EE[\xi_{k,n} \ind{\xi_{k,n} > -t}] \right),
 \end{split}
\end{equation}
note here that we are using that we have centered the walk by the mean
rather than the limiting velocity so that
$\EE[\xi_{k,n} \ind{\xi_{k,n} \leq -t}] + \EE[\xi_{k,n} \ind{\xi_{k,n} > -t}] = \EE[\xi_{k,n}] = 0$.

We have already shown that 
\begin{itemize}
 \item $W_{n,t} \underset{n\to\infty}{\Longrightarrow} \Psi(N_t^{(g)}) - E[\Psi(N_t^{(g)})]$ (by Steps 1 and 2),
 \item and $\lim_{t\to 0} \limsup_{n\to\infty} \PP(|\mathfrak{E}_{n,t}| > \delta ) = 0$, for all $\delta > 0$ (by Step 3). 
\end{itemize}
Therefore, by~\cite[Theorem 3.2]{Bi99} the proof is complete if we show that $\Psi(N_t^{(g)}) - E[\Psi(N_t^{(g)})]$ converges in distribution to $\rW_{\lambda_g} + (1-g(\infty))^{1/s} \rS_s$ as $t\to 0$. This can be seen with the help of Campbell's theorem~\cite[Section 3.2]{Kin93}  if we note that for any fixed $u\in\mathbb{R}$
\begin{align}
 \lim_{t\to 0}E\left[ e^{i u \left( \Psi(N_t^{(g)}) - E[\Psi(N_t^{(g)})] \right) } \right] &= \lim_{t\to 0} \exp\left\{ \int_{-\infty}^{-t} \left( e^{i u x} - 1 - iux \right) \hlg(x) \, dx \right\} \nonumber \\
&= \exp\left\{ \int_{-\infty}^{0} \left( e^{i u x} - 1 - iux \right) \lambda_g(x) \, dx \right\} \label{limcf1} \\
&\quad \times \exp\left\{ (1-g(\infty)) \int_{-\infty}^{0} \left( e^{i u x} - 1 - iux \right) \lambda_0(x) \, dx \right\}. \label{limcf2} 
\end{align}
  For the term in \eqref{limcf2}, standard computations in complex analysis yield that
  \[\int_{-\infty}^0 (e^{i u x}-1-iux)|x|^{-s-1} \, dx = \frac{-\Gamma(1-s)\cos(\frac{\pi s}{2}) |u|^2}{s} \left( 1 + i \frac{u}{|u|} \tan(\frac{\pi s}{2} ) \right),
  \]
  from which (recalling the definition of $\lambda_0(x)$ above and the relation between the constants $b$ and $\tc$ in \eqref{relate}, and the formula for the characteristic function of $\rS_s$ in \eqref{stablecf}) we have that
  \[
  \eqref{limcf2} = \exp\left\{ -(1-g(\infty)) b \left( 1 + i \frac{u}{|u|} \tan(\frac{\pi s}{2} ) \right) \right\}
  = E\left[ e^{i u (1-g(\infty))^{1/s} \rS_s} \right].
  \] 
The term in \eqref{limcf1} clearly equals the characteristic function of $\rW_{\lambda_g}$ as defined in \eqref{Wlambdacritical}, but we still need to justify that the function $\lambda_g(x)$ satisfies the properties required of generalized tempered stable laws in Definition \ref{gtslaw}. 
For this, it is enough to check that
$t\mapsto t^{s+1} \lambda_g(-t)$ is a non-increasing, continuous function on $(0,\infty)$ which vanishes at $\infty$ and has a finite limit as $t\to 0^+$. 
These can be checked from \eqref{lgdef} by re-writing 
 \begin{align}
 t^{s+1} \lambda_g(-t) &= \tc \int_{t/\vs}^\infty \left( \vs s - \frac{(s-1)t}{x} \right) \, g(dx)  \nonumber \\
&= \tc \vs s \left( g(\infty)- g(t/\vs) \right) - \tc(s-1)  \int_{t/\vs}^\infty \frac{t}{x} \, g(dx) \label{f2} \\
&= \tc \vs s g(\infty) - \tc \vs g(t/\vs) -\tc (s-1)t \int_{t/\vs}^\infty \frac{g(x)}{x^2} \, dx, \label{f4}
\end{align}
where the last equality follows from integration by parts.
Re-writing the integral in \eqref{f2} as $\int_0^\infty \frac{t}{x}\ind{x\geq t/\vs} \, g(dx)$, it follows from the dominated convergence theorem that 
\[\lim_{t\to \infty} t^{s+1}\lambda_g(-t) = 0
\quad\text{and}\quad \lim_{t\to 0^+} t^{s+1}\lambda_g(-t) = \tc \vs s g(\infty)<\infty.
\]
    Finally, using the representation \eqref{f4} one can see that $t^{s+1} \lambda_g(-t)$ is continuous and non-increasing as a function of $t$. The only difficulty here is in showing that the last term in \eqref{f4} is non-increasing, but this follows from the fundamental theorem of calculus and then integration by parts:
\[
 \frac{d}{dt} \left\{  -t \int_{t/\vs}^\infty \frac{g(x)}{x^2} \, dx \right\} = - \int_{t/\vs}^\infty \frac{g(x)}{x^2} \, dx + \frac{g(t/\vs)}{t/\vs} = -\int_{t/\vs}^\infty \frac{1}{x} \, g(dx) \leq 0. 
\]
 
\qed

\subsection{Polynomial cooling: proofs of critical and stable limits in Theorem \ref{PolyTransition}  }\label{polycases}
We show here how~\eqref{Xnstableliml} and~\eqref{Xnstablelim}
follow as a corollary of the general Theorem~\ref{GeneralConditions-S}. 
\begin{proof}[Proof of   \eqref{Xnstableliml} -- critical case]
  If $T_k \sim A k^{1/(s-1)}$, then $\tau(n) \sim \frac{A(s-1)}{s} n^{s/(s-1)}$.
  This implies {\bf(S2)}.  
Moreover, if $x > K_{A,s} := A^{\frac{s-1}{s}} \left(\frac{s}{s-1}\right)^{1/s}$ then for $n$ large enough one has $T_k \leq x \tau(n)^{1/s}$ for all $k\leq n$, whereas if 
$0<x\leq K_{A,s}$ then for any $\epsilon>0$ and $n$ large enough
\[
  \left\{k:  k\leq (1-\epsilon) \left(\frac{x}{K_{A,s}}\right)^{s-1} n \right\}\subset \{k\leq n: T_k \leq x \tau(n)^{1/s} \} \subset \left\{k:  k\leq (1+\epsilon) \left(\frac{x}{K_{A,s}}\right)^{s-1} n \right\}. 
\]
This implies that
\begin{align*}
 \lim_{n\to\infty} \frac{ \sum_{k=1}^n T_k \ind{T_k < x \tau(n)^{1/s}} }{ \tau(n) }  
= \begin{cases}
   \left(\frac{x}{K_{A,s}}\right)^s & \text{if } x \leq K_{A,s} \\
   1 & \text{if } x > K_{A,s}.
  \end{cases}
\end{align*}
That is, condition \eqref{cag} holds with 
$g(x) = 1 \wedge   \left(\frac{x}{K_{A,s}}\right)^s$. 
Hence, \eqref{Xnstablelim2} in Theorem \ref{GeneralConditions-S} implies convergence to the generalized stable 
variable  $\rW_{\lambda_g}$. 
Moreover, with this choice of $g$ the function $\lambda_g$ defined in \eqref{lgdef} can be calculated to be  
$\lambda_g(x) = \tc \vs s |x|^{-s-1} \left( 1 + \frac{x}{K_{A,s} \vs} \right)_+$. 
By~\eqref{lcr}, it follows that $\lambda_g = \lambda_{c,r}$ with $c=\tc \vs s$ and $r = \vs K_{A,s} = \vs A^{\frac{s-1}{s}} \left(\frac{s}{s-1}\right)^{1/s}$.
\end{proof}

\begin{proof}[Proof of  \eqref{Xnstablelim} -- supercritical case]
$T_k \sim A k^a$ implies that $\tau(n) \sim \frac{1}{a+1} n^{a+1}$ and so {\bf(S2)} holds. 
This also implies that for any fixed $x\in (0,\infty)$ and $n$ sufficiently large we have that $T_k \leq x \tau(n)^{1/s}$ implies that $k\leq \left(\frac{2x}{A(a+1)^{1/s}} \right)^{1/a} n^{\frac{a+1}{as}}$. 
Since the exponent $\frac{a+1}{as} < 1$ when $a > \frac{1}{s-1}$, it follows that $\sum_{k=1}^n T_k \ind{T_k \leq x \tau(n)^{1/s}} = o(n^{a+1})$, and since $\tau(n) \sim \frac{1}{a+1} n^{a+1}$ it follows that \eqref{cag} holds with $g(x) = 0$ for all $x\in(0,\infty)$. 
\end{proof}

\section{Proofs: mixed limiting distributions}\label{S8}


We show here Theorem~\ref{arbmix}. The basic idea is that
one can combine the polynomial cooling maps in
Theorem~\ref{PolyTransition} into a new polynomial map
so as to obtain a mixture of their limiting laws.
We first show how to obtain mixture of two polynomial cooling maps,
Proposition~\ref{intertwine}. We then use the latter to obtain
a large class of generalized stable laws, see Example \ref{ex:finmix}.
At this point we prove in Lemma~\ref{closure} and Corollary~\ref{CoroClosed}
via a limiting closure argument on the class identified in 
Example \ref{ex:finmix}
that indeed it is possible to obtain any generalized
tempered stable $\rW_{\lambda}$ with $\lambda\in\Lconv$.
Theorem~\ref{PolyTransition} is then readily obtained by
Corollary~\ref{CoroClosed} and Proposition~\ref{intertwine}. 

Let  $I$ be an index set. For each $i \in I$ let the cooling maps  $\tau^{(i)}: \bb{N} \to \bb{N}$ be associated with the increment sequences $(T^{(i)}_k , k \in \bb{N})$ by $\tau^{(i)}(n) = \sum_{k=1}^nT^{(i)}_k$, 
and let $X^{(i)} = (X_n^{(i)}){n\geq 0}$ be a RWCRE corresponding to the cooling map $\tau^{(i)}$.
Given a function  $\sigma: \bb{N} \to I$ we define  the $\sigma$-interweaving of the cooling maps to the map $\tau^\sigma$ that corresponds to the increment sequence $(T^\sigma_k , k \in \bb{N})$ defined by
\begin{equation}\label{Tkmix}
T^\sigma_ k = T^{(\sigma(k))}_{M_{k,\sigma(k)}},
\quad
 \text{where } M_{k,i} := \#\{j \leq k : \sigma(j)  = i\}.
\end{equation}
That is,  $M_{k,i}$
counts the number of times the increment sequence $i$  has been selected 
in the first $k$ cooling increments.
We refer to $\sigma$ as a selection function and we assume
that it selects each increment sequence infinitely many times,
i.e., we assume that
$\lim_{k \to \infty} M_{k,i}  = \infty$ for all $i \in I.$
For our first result in this section we will need to assume a few
conditions on the cooling maps $\tau^{(i)}$.
For each $i \in I$ we will assume that there are constants $b_i, C_i > 0$, $\alpha_i \geq 1$, $\beta_i \in [1/2, 1/s]$, and a random variable $\mathcal{X}_i$ such that
\begin{equation}\label{taumixcond}
  \max_{k\leq n} T_k^{(i)} = o(n^{\alpha_i \beta_i s}), \quad \tau^{(i)}(n) \sim C_i n^{\alpha_i}, \quad \text{and}\quad \frac{X_n^{(i)} - \EE[X_n^{(i)}]}{b_i n^{\beta_i}} \limdist \mathcal{X}_i.
\end{equation}
We are now ready to state the first result.
\begin{proposition}[{\bf Interweaving two polynomial maps}]\label{intertwine}
Let $\tau^{(1)}$ and $\tau^{(2)}$ be cooling maps satisfying the conditions in \eqref{taumixcond}. 
Given any constants  $a_1, a_2 > 0$, there exists a cooling map $\tau$ such that 
for some constants $b,C>0$,
$\alpha = \frac{(\alpha_1 \beta_1) \wedge (\alpha_2 \beta_2)}{\beta_1 \wedge \beta_2}$,
and $\beta = \beta_1 \wedge \beta_2$, we have that 
\begin{align}
 &\max_{k\leq n} T_k = o(n^{\alpha\beta s}), \quad \tau(n) \sim C n^\alpha, \label{Ttaucond}\\
&\text{and}\quad \frac{X_n - \EE[X_n]}{b n^{\beta}} \limdist a_1 \mathcal{X}_1 + a_2 \mathcal{X}_2. \label{mixedlimit}
\end{align}
\end{proposition}

\begin{proof}
First of all, we claim that it is enough to construct a cooling map $\tau$ that satisfies \eqref{Ttaucond} and
also 
\begin{equation}\label{mixedlimit-ss}
 \frac{X_{\tau(n)} - \EE[X_{\tau(n)}]}{b\tau(n)^\beta} \limdist a_1 \mathcal{X}_1 + a_2 \mathcal{X}_2. 
\end{equation}
Indeed, since 
\[
 \frac{X_n - \EE[X_n]}{b n^{\beta}} = \frac{X_{\tau(\ell_n)} - \EE[X_{\tau(\ell_n)}]}{b \tau(\ell_n)^\beta} \left(\frac{\tau(\ell_n)}{n} \right)^\beta
+ \frac{Z^{(\ell_n+1)}_{\bar{T}_n} - \EE[Z_{\bar{T}_n}]}{\bar{T}_n^{1/s}} \frac{\bar{T}_n^{1/s}}{b n^\beta}, 
\]
then arguing as in the proof of Lemma \ref{lem:boundary} condition \eqref{mixedlimit}
will follow from \eqref{mixedlimit-ss} if we can show that conditions in \eqref{Ttaucond} imply
that $\frac{\tau(\ell_n)}{n} \to 1$ and $\frac{\bar{T}_n^{1/s}}{b n^\beta} \to 0$ as $n\to \infty$.
For the first of these note that \eqref{Ttaucond} implies that
$1 \geq \frac{\tau(\ell_n)}{n} \geq 1 - \frac{\bar{T}_n}{n} \geq 1 - \frac{T_{\ell_n+1}}{\tau(\ell_n)} = 1 - o\left( \ell_n^{-\alpha(1-\beta s)} \right) = 1-o(1)$ since $\beta \leq 1/s$. 
The second follows similarly since \eqref{Ttaucond} implies
$\frac{\bar{T}_n^{1/s}}{b n^\beta} \leq \frac{(\bar{T}_{\ell_n+1})^{1/s}}{\tau(\ell_n)^\beta} = o(1)$. 

We still need to construct a cooling map $\tau$ satisfying \eqref{Ttaucond} and \eqref{mixedlimit-ss}.
To this end,
given a selection function $\sigma: \bb{N} \to \{1,2\}$, let $\tau = \tau^\sigma$  and
consider the following decomposition
\begin{equation}\label{mixsplit}
\begin{split}
 \frac{X_{\tau(n)} - \EE[X_{\tau(n)}]}{b\tau(n)^\beta}
&\Law \frac{b_1}{b} \left( \frac{\tau^{(1)}(M_{n,1})^{\beta_1} }{\tau(n)^\beta} \right) \frac{X^{(1)}_{\tau^1(M_{n,1})} - \EE[X^{(1)}_{\tau^1(M_{n,1})}]}{b_1 \tau^{(1)}(M_{n,1})^\beta} \\
&\qquad + \frac{b_2}{b} \left( \frac{\tau^{(2)}(M_{n,2})^{\beta_2} }{\tau(n)^\beta} \right) \frac{X^{(2)}_{\tau^2(M_{n,2})} - \EE[X^{(2)}_{\tau^2(M_{n,2})}]}{b_2 \tau^{(2)}(M_{n,2})^\beta}.
\end{split}
\end{equation}
By the third condition in \eqref{taumixcond} and the assumption that $M_{n,i}\to\infty$ for $i=1,2$,
the last fractions on each of the two terms on the right converge in distribution to $\mathcal{X}_1$ and $\mathcal{X}_2$, respectively. We need to choose the mixing function $\sigma$ and the exponent $\beta$ so that the middle fractions for the terms on the right side converge to constants.
How we do this depends on the relative values of $\beta_1$, $\beta_2$,  $\alpha_1 \beta_1$ and $\alpha_2 \beta_2$.  
Without loss of generality we can assume that $\beta_1 \leq \beta_2$. 
We will describe the mixing function $\sigma$ only in terms of the asymptotics of $M_{n,1}$ (or $M_{n,2}$). It is not hard to then give explicit mixing functions which have these asymptotics. 

\noindent\textbf{Case I: $\alpha_1= \alpha_2$ and $\beta_1 = \beta_2$.}
In this case let $M_{n,1} \sim \theta n$ for a value of $\theta \in (0,1)$ to be chosen later (so that necessarily $M_{n,2} \sim (1-\theta) n$).
Let $\alpha = \alpha_1 = \alpha_2$ and $\beta = \beta_1 = \beta_2$. 
Then, 
$ \max_{k\leq n} T_k = \max_{k\leq M_{n,1}} T_k^{(1)} \vee \max_{k\leq M_{n,2}} T_k^{(2)} 
= o(M_{n,1}^{\alpha_1 \beta_1 s}) \vee o( M_{n,2}^{\alpha_2 \beta_2 s}) = o(n^{\alpha \beta s})$, 
and also 
\begin{align*}
 \tau^{(1)}(M_{n,1}) &\sim C_1 M_{n,1}^{\alpha_1} \sim C_1 \theta^{\alpha_1} n^{\alpha_1} = C_1 \theta^\alpha n^\alpha, \\
 \tau^{(2)}(M_{n,2}) &\sim C_2 M_{n,2}^{\alpha_2} \sim C_2 (1-\theta)^{\alpha_2} n^{\alpha_2} = C_2 (1-\theta)^\alpha n^\alpha \\
\text{and } \quad  \tau(n) &= \tau^{(1)}(M_{n,1}) + \tau^{(2)}(M_{n,2}) \sim \left( C_1 \theta^\alpha + C_2 (1-\theta)^\alpha \right) n^{\alpha}. 
\end{align*}
Thus,
it follows from \eqref{mixsplit} that 
\[
 \frac{X_{\tau(n)} - \EE[X_{\tau(n)}]}{b\tau(n)^\beta}
 \limdist \frac{b_1}{b} \left( \frac{C_1 \theta^{\alpha} }{C_1 \theta^\alpha + C_2 (1-\theta)^\alpha} \right)^{\beta}\bigg[ \mathcal{X}_1 + \frac{b_2}{b_1} \left( \frac{C_2 (1-\theta)^\alpha}{C_1 \theta^\alpha } \right)^\beta \mathcal{X}_2
   \bigg]. 
\]
Finally, we choose  $\theta \in (0,1)$ and $b>0$ so that the right side is equal to $a_1 \mathcal{X}_1 + a_2 \mathcal{X}_2$. 
More explicitly, let 
\[
\theta = \left( 1 + \left( \frac{a_2 b_1}{a_1 b_2} \right)^{1/(\alpha \beta)} \left( \frac{C_1}{C_2} \right)^{1/\alpha} \right)^{-1},
\quad \text{and} \quad 
 b = \frac{b_1}{a_1}\left( 1 + \left( \frac{a_2 b_1}{a_1 b_2} \right)^{1/ \beta} \right)^{-\beta}.
\]

For the remaining four cases we will give fewer details
and leave it to the reader to check that in each case the
parameters can be chosen so that the limiting distribution
is equal to $a_1 \mathcal{X}_1 + a_2 \mathcal{X}_2$.

\noindent\textbf{Case II: $\beta_1 = \beta_2$ and $\alpha_1\beta_1 \neq \alpha_2 \beta_2$.}
Without loss of generality we can assume that $\alpha_1 \beta_1 > \alpha_2 \beta_2$ (or equivalently $\alpha_1 > \alpha_2$). 
In this case we will let $\alpha = \alpha_2$, $\beta = \beta_1 = \beta_2$, 
and $M_{n,1} \sim \theta n^{\frac{\alpha_2}{\alpha_1}}$
for some $\theta>0$ to be chosen later 
(since the exponent $\frac{\alpha_2}{\alpha_1}< 1$ this implies that $M_{n,2} = n-M_{n,1} \sim n$). 
Then one can check that 
$\max_{k\leq n} T_k = o(n^{\alpha_2 \beta_1 s}) = o(n^{\alpha\beta s})$, 
$\tau(n) \sim \left(C_1 \theta^{\alpha_1} + C_2\right) n^{\alpha}$, and 
\[
 \frac{X_{\tau(n)} - \EE[X_{\tau(n)}]}{b\tau(n)^\beta}
\limdist \frac{b_1}{b} \left( \frac{C_1 \theta^{\alpha_1}}{C_1 \theta^{\alpha_1} + C_2} \right)^{\beta} \bigg[\mathcal{X}_1 + \frac{b_2}{b_1}\left( \frac{C_2}{C_1 \theta^{\alpha_1}} \right)^{\beta} \mathcal{X}_2\bigg]. 
\]

\noindent\textbf{Case III: $\beta_1 < \beta_2$ and $\alpha_1\beta_1 = \alpha_2 \beta_2$.}
Note that in this case we necessarily have $\alpha_1 > \alpha_2$. 
In this case we will let $\alpha = \alpha_1$, $\beta= \beta_1$, and
$M_{n,1} \sim \theta n$ for a value of $\theta \in (0,1)$ to be
chosen later (so that necessarily $M_{n,2} \sim (1-\theta) n$). 
Then one can check that 
$\max_{k\leq n} T_k = o(n^{((\alpha_1 \beta_1) \vee (\alpha_2 \beta_2)) s}) = o(n^{\alpha\beta s})$, 
$\tau(n) \sim C_1 \theta^{\alpha}n^{\alpha}$, and 
\[
\frac{X_{\tau(n)} - \EE[X_{\tau(n)}]}{b\tau(n)^\beta}
\limdist \frac{b_1}{b} \bigg[\mathcal{X}_1 + \frac{b_2}{b_1} \frac{C_2^{\beta_2} (1-\theta)^{\alpha_2\beta_2} }{C_1^{\beta} \theta^{\alpha \beta}} \mathcal{X}_2\bigg]. 
\]

\noindent\textbf{Case IV: $\beta_1 < \beta_2$ and $\alpha_1 \beta_1 > \alpha_2 \beta_2$}. 
In this case we will let $\alpha = \frac{\alpha_2 \beta_2}{\beta_1}$, $\beta = \beta_1$,
and  $M_{n,1} \sim \theta n^{\frac{\alpha_2 \beta_2}{\alpha_1 \beta_1}}$ for some $\theta > 0$ to be chosen later (since the exponent $\frac{\alpha_2 \beta_2}{\alpha_1\beta_1}< 1$ this implies that $M_{n,2} = n-M_{n,1} \sim n$). 
Then one can check that 
$\max_{k\leq n} T_k = o(n^{\alpha_2 \beta_2 s}) = o(n^{\alpha\beta s})$, 
$\tau(n) \sim C_1 \theta^{\alpha_1}n^{\alpha}$, and 
\[
 \frac{X_{\tau(n)} - \EE[X_{\tau(n)}]}{b\tau(n)^\beta}
\limdist \frac{b_1}{b} \bigg[\mathcal{X}_1 + \frac{b_2}{b_1}\frac{C_2^{\beta_2}}{C_1^{\beta}\theta^{\alpha_1 \beta}} \mathcal{X}_2\bigg]. 
\]

\noindent\textbf{Case V: $\beta_1 < \beta_2$ and $\alpha_1 \beta_1 < \alpha_2 \beta_2$}. 
In this case we will let $\alpha = \alpha_1$, $\beta = \beta_1$, 
and $M_{n,2} \sim \theta n^{\frac{\alpha_1 \beta_1}{\alpha_2\beta_2}}$ for some $\theta>0$ to be chosen later 
(since the exponent $\frac{\alpha_1\beta_1}{\alpha_2 \beta_2}< 1$ this implies that $M_{n,1} = n-M_{n,2} \sim n$).
Then one can check that 
$\max_{k\leq n} T_k = o(n^{\alpha_1 \beta_1 s}) = o(n^{\alpha\beta s})$, 
$\tau(n) \sim C_1 n^{\alpha}$, and 
\[
 \frac{X_{\tau(n)} - \EE[X_{\tau(n)}]}{b\tau(n)^\beta}
\limdist \frac{b_1}{b} \bigg[\mathcal{X}_1 + \frac{b_2}{b_1}\frac{C_2^{\beta_2}\theta^{\alpha_2 \beta_2}}{C_1^{\beta}}\mathcal{X}_2\bigg]. 
\]
\end{proof}
Since the polynomial cooling maps from Theorem \ref{PolyTransition}
satisfy the conditions in \eqref{taumixcond},
it follows from Lemma \ref{intertwine} that 
by interweaving a finite number of these polynomial
cooling maps we can obtain a cooling map whose
corresponding RWCRE converges (after proper centering and scaling)
to any finite linear combination of the limit laws captured in Theorem~\ref{PolyTransition}. 
In particular, by intertweaving a finite number of critical
polynomial cooling maps we can obtain any limiting distributions
of the form $\sum_{i=1}^\ell a_i \rW_{\lambda_{c,r_i}}$ where $c = \tc \vs s$,
and $a_i, r_i > 0$ for $i=1,2,\ldots,\ell$.
To give a simpler characterization of this type of limiting distribution,
we use the following properties of the generalized tempered stable laws
which are easy to check from the definition: 
(1) $a \rW_{\lambda_{c,r}} \Law \rW_{\lambda_{a^s c, ar}}$, and
(2) if $\rW_{\lambda}$ and $\rW_{\lambda'}$ are independent,
then $\rW_{\lambda} + \rW_{\lambda'} \Law \rW_{\lambda+\lambda'}$.
From this, it follows that 
$\sum_{i=1}^\ell a_i \rW_{\lambda_{c,r_i}} \Law \rW_\lambda$, where 
\[
 \lambda(x) = \sum_{i=1}^\ell \lambda_{a_i^s c, a_i r_i}(x)
 =  |x|^{-s-1} \sum_{i=1}^\ell c a_i^s \left( 1 + \frac{x}{a_i r_i} \right)_+. 
 \]
From this we see that we can characterize the limiting distributions
of this type as generalized tempered stable random variables $\rW_\lambda$
where $\lambda(x) = c|x|^{-s-1}a(x)$ and $a(x)$ is a convex and
piecewise linear function. In fact, as the following example shows,
by choosing the interweaving carefully we can attain a limiting
distribution of this form for any such convex piecewise linear function $a(x)$. 

\begin{example}\label{ex:finmix}
Let 
\begin{equation}\label{aghdef}
 a(x) = \sum_{i=1}^\ell (g_i x + h_i)_+
\end{equation}
with $g_i, h_i > 0$ for $i\leq \ell$ and $\sum_{i=1}^\ell h_i = 1$ be a generic piecewise linear convex funciton on $(-\infty,0]$ that vanishes at $-\infty$ and has $a(0) = 1$. 
Now, for each $i\leq \ell$ let $\tau^{(i)}$ be a critical polynomial cooling map with $T_k^{(i)} \sim A_i k^{1/(s-1)}$, and let $\tau = \tau^\sigma$ be an interweaving of these cooling maps where the mixing function $\sigma$ is chosen so that $\lim_{k\to\infty} \frac{M_{k,i}}{k} = \theta_i \in (0,1)$ for each $i\leq \ell$ with $\sum_{i=1}^\ell \theta_i = 1$. 
If we choose the parameters $\theta_i$ and $A_i$ for constructing this cooling map so that 
\[
 \theta_i = \frac{g_i}{\sum_{j=1}^\ell g_j}, \quad \text{and}\quad 
 A_i = \frac{\left(\frac{s-1}{s}\right)^{1/(s-1)} h_i}{\vs^{s/(s-1)} g_i^{s/(s-1)}}, \quad i\leq \ell, 
\]
then by repeating the sort of computation in the proof of Lemma \ref{intertwine} the reader can check that the corresponding RWCRE $X$ has limiting distribution $\frac{X_n - \EE[X_n]}{n^{1/s}} \limdist \rW_{\lambda}$ with $\lambda(x) = \tc \vs s |x|^{-s-1} a(x)$ with $a(x)$ as in \eqref{aghdef}. 
Moreover, one can also check that the growth rate of this cooling map is $\tau(n) \sim \left( \frac{ (s-1)}{ s \vs \sum_{i=1}^\ell g_i } \right)^{s/(s-1)} n^{s/(s-1)}$, and since $\sum_{i=1}^\ell g_i = a'(0)$ this gives a relation between the growth rate of $\tau(n)$ and $a'(0)$ for this particular type of cooling map that we will use in the proof of Corollary \ref{CoroClosed} below. 
\end{example}

\begin{lemma}[{\bf Closure for critical maps}]\label{closure}
Suppose that for each $j\geq 1$, $\tau^{(j)}$ is a cooling map such that the corresponding RWCRE $X^{(j)}$ has limiting distribution 
\begin{equation}\label{Xjlim}
 \frac{X^{(j)}_n - \EE[X^{(j)}_n]}{n^{1/s}} \limdist \rW_{\lambda_j},
\end{equation}
and the cooling map has asymptotic growth rate $\tau^{(j)}(n) \sim K_j n^{s/(s-1)}$
as $n\to\infty$ for some $K_j \in (0,\infty)$. 
If we also assume that 
\begin{equation}\label{Wljlim}
  \rW_{\lambda_j} \underset{j\to\infty}{\limdist} \rW_{\lambda}
  \quad\text{and}\quad
  \lim_{j\to\infty} \frac{K_{j+1}}{K_j} = 1,
\end{equation}
then there exists a cooling map $\tau$ such that the corresponding RWCRE $X$ 
has a limiting distribution $\frac{X_n - \EE[X_n]}{n^{1/s}} \limdist \rW_{\lambda}$.
\end{lemma}

\begin{proof}
The new cooling map $\tau$ will be constructed from from the cooling maps $\tau^{(j)}$ as follows. We will choose an increasing sequence of integers $0=m_0 < m_1 < m_2 < m_3 < \cdots$  with properties given below and then construct the cooling map $\tau$ by choosing the $k$-th cooling interval from the cooling map $\tau^{(j)}$ if $m_{j-1} < k \leq m_j$. That is, 
\[
 T_k = T^{(j)}_k, \quad \text{if } m_{j-1} < k \leq m_j. 
\]
We will choose the sequence of integers $m_j$ in the following manner. Assuming that $m_{j-1}$ has already been determined, we choose $m_j$ large enough so that 
\begin{align}
 \left| \frac{\tau^{(j+1)}(m_j)}{\tau(m_j)} - 1 \right| \leq 2 \left| \frac{K_{j+1}}{K_j} - 1 \right| =: \eta_j \hspace{1.5in} \label{mj1} \\
 \sup_{n\geq \tau^{(j+1)}(m_j)} \sup_{x \in \R} \left| \PP\left( \frac{X^{(j+1)}_{n} - \EE[X^{(j+1)}_{n}] }{ n^{1/s} } \leq x \right) - P\left( \rW_{\lambda_{j+1}} \leq x \right) \right| < \frac{1}{j} \label{mj2} \\
\text{and} \qquad\qquad \sup_{x \in \R} \left| \PP\left( \frac{X_{\tau(m_j)} - \EE[X_{\tau(m_j)}] }{ \tau(m_j)^{1/s} } \leq x \right) - P\left( \rW_{\lambda_j} \leq x \right) \right| < \frac{1}{j}. \label{mj3}
\end{align}
Condition \eqref{mj2} follows easily from the assumption in \eqref{Xjlim} by taking $m_j$ sufficiently large.
For condition \eqref{mj1}, first note that
$\frac{\tau^{(j+1)}(m_j)}{\tau(m_j)} = \frac{\tau^{(j+1)}(m_j)}{\tau^{(j)}(m_j) - \tau^{(j)}(m_{j-1}) + \tau(m_{j-1})}$. 
The assumptions on the growth rate of $\tau^{(j)}$ and $\tau^{(j+1)}$
imply that by taking $m_j$ sufficiently large,
we can make this fraction arbitrarily close to $\frac{K_{j+1}}{K_j}$
so that \eqref{mj1} is satisfied. 
Finally, for \eqref{mj3} we note that we can expand our
probability space to include all of the walks $X^{(j)}$
and so that we can construct the walk $X_\cdot$ by letting 
\[
\{X_{\tau(m_{j-1})+k} - X_{\tau(m_{j-1})} \}_{k\leq \tau(m_{j}) - \tau(m_{j-1})} = \{ X^{(j)}_{\tau^{(j)}(m_{j-1})+k} - X^{(j)}_{\tau^{(j)}(m_{j-1})} \}_{k\leq \tau^{(j)}(m_{j}) - \tau^{(j)}(m_{j-1})}.
\] 
(Note that we are using here that $\tau(m_j) - \tau(m_{j-1}) = \tau^{(j)}(m_j) - \tau^{(j)}(m_{j-1})$, so that both sequences above have the same number of terms.)
Using this construction we then have that 
\begin{align*}
 \frac{X_{\tau(m_j)} - \EE[X_{\tau(m_j)}] }{ \tau(m_j)^{1/s} }
&= \frac{X^{(j)}_{\tau^{(j)}(m_j)} - \EE[X^{(j)}_{\tau^{(j)}(m_j)}]}{\tau^{(j)}(m_j)^{1/s}} \left( \frac{\tau^{(j)}(m_j)}{\tau(m_j)} \right)^{1/s} \\
&\quad + \frac{X_{\tau(m_{j-1})} - \EE[X_{\tau(m_{j-1})}] }{ \tau(m_j)^{1/s} } - \frac{X^{(j)}_{\tau^{(j)}(m_{j-1})} - \EE[X^{(j)}_{\tau^{(j)}(m_{j-1})}] }{ \tau(m_j)^{1/s} }.
\end{align*}
If $m_{j-1}$ has already been fixed, then \eqref{Xjlim} implies that as $m_j \to \infty$ the right side above converges in distribution to $\rW_{\lambda_j}$. Thus, we can take $m_j$ large enough so that \eqref{mj3} holds. 

Suppose now that $\tau(m_j) < n \leq \tau(m_{j+1})$. Since $X_{\tau(m_j)}$ and $X_n-X_{\tau(m_j)}$ are independent, we can write 
\begin{align}
& \PP\left( \frac{X_n - \EE[X_n]}{n^{1/s}} \leq x \right) \nonumber \\
&= \int_\R \PP\left( \frac{X_{\tau(m_j)} - \EE[X_{\tau(m_j)}]}{n^{1/s}} \leq x-y \right) \PP\left( \frac{ X_n - X_{\tau(m_j)} - \EE[ X_n - X_{\tau(m_j)}] }{n^{1/s}} \in dy \right) \label{convolution}
\end{align}
For the first probability inside the integral, it follows from \eqref{mj2} and \eqref{mj3} that 
\begin{align}
& \PP\left( \frac{X_{\tau(m_j)} - \EE[X_{\tau(m_j)}]}{n^{1/s}} \leq x-y \right) 
=  \PP\left( \left( \frac{\tau(m_j)}{n} \right)^{1/s} \frac{X_{\tau(m_j)} - \EE[X_{\tau(m_j)}]}{\tau(m_j)^{1/s}} \leq x-y \right) \nonumber \\
&\leq \PP\left( \left( \frac{\tau(m_j)}{n} \right)^{1/s}\frac{X^{(j+1)}_{\tau^{(j+1)}(m_j)} - \EE[X^{(j+1)}_{\tau^{(j+1)}(m_j)}] }{ \tau^{(j+1)}(m_j)^{1/s} } \leq x-y \right) + \frac{2}{j} + \delta_j, \label{conv1}
\end{align}
where $\delta_j:= \sup_{x \in \R} \left| P(\rW_{\lambda_j} \leq x ) - P(\rW_{\lambda_{j+1}} \leq x ) \right|$ (note that \eqref{Wljlim} implies that  $\delta_j \to 0$ as $j \to \infty$).
For the second probability in \eqref{convolution}, note that we can replace $X_n - X_{\tau(m_j)}$ with $X^{(j+1)}_{\tau^{(j+1)}(m_j) + n- \tau(m_j)} - X^{(j+1)}_{\tau^{(j+1)}(m_j)}$,
so that applying \eqref{conv1} to \eqref{convolution} we can conclude that 

\begin{align*}
 & \PP\left( \frac{X_n - \EE[X_n]}{n^{1/s}} \leq x \right) \\
&\leq \PP\left( \frac{ X^{(j+1)}_{\tau^{(j+1)}(m_j) + n - \tau(m_j)} - X^{(j+1)}_{\tau^{(j+1)}(m_j)} - \EE[X^{(j+1)}_{\tau^{(j+1)}(m_j) + n - \tau(m_j)} - X^{(j+1)}_{\tau^{(j+1)}(m_j)}] }{n^{1/s}} \right. \\
&\qquad\qquad \left. + \left( \frac{\tau(m_j)}{n} \right)^{1/s}\frac{X^{(j+1)}_{\tau^{(j+1)}(m_j)} - \EE[X^{(j+1)}_{\tau^{(j+1)}(m_j)}] }{ \tau^{(j+1)}(m_j)^{1/s} } \leq x \right) + \frac{2}{j} + \delta_j \\
&=  \PP\left( \frac{ X^{(j+1)}_{\tau^{(j+1)}(m_j) + n - \tau(m_j)} - \EE[X^{(j+1)}_{\tau^{(j+1)}(m_j) + n - \tau(m_j)}] }{n^{1/s}} \right. \\
&\qquad\qquad \left. + \left( \left( \frac{\tau(m_j)}{n} \right)^{1/s} - \left( \frac{\tau^{(j+1)}(m_j)}{n} \right)^{1/s}  \right) \frac{X^{(j+1)}_{\tau^{(j+1)}(m_j)} - \EE[X^{(j+1)}_{\tau^{(j+1)}(m_j)}] }{ \tau^{(j+1)}(m_j)^{1/s} } \leq x \right) + \frac{2}{j} + \delta_j. 
\end{align*}
Note that the above inequality holds for all $x \in \R$ and $n \in (\tau(m_j), \tau(m_{j+1})]$. 
Since $n > \tau(m_j)$ it follows from \eqref{mj1} that 
\[
 \left| \left(\frac{\tau(m_j)}{n} \right)^{1/s} - \left( \frac{\tau^{(j+1)}(m_j)}{n} \right)^{1/s} \right|
< \left| 1 - \left( \frac{\tau^{(j+1)}(m_j)}{\tau(m_j)} \right)^{1/s} \right|
\leq \left| 1 - \frac{\tau^{(j+1)}(m_j)}{\tau(m_j)}  \right|^{1/s} \leq \eta_j^{1/s}, 
\]
Therefore, for any $\epsilon > 0$, $x \in \R$, and $n \in (\tau(m_j), \tau(m_{j+1})]$ we have 
\begin{align*}
 \PP\left( \frac{X_n - \EE[X_n]}{n^{1/s}} \leq x \right) 
&\leq \PP\left( \frac{ X^{(j+1)}_{\tau^{(j+1)}(m_j) + n - \tau(m_j)} - \EE[X^{(j+1)}_{\tau^{(j+1)}(m_j) + n - \tau(m_j)}] }{n^{1/s}} \leq x + \epsilon \right) \\
&\qquad + \PP\left( \left| \frac{X^{(j+1)}_{\tau^{(j+1)}(m_j)} - \EE[X^{(j+1)}_{\tau^{(j+1)}(m_j)}] }{ \tau^{(j+1)}(m_j)^{1/s} } \right| \geq \eta_j^{-1/s} \epsilon \right) + \frac{2}{j} + \delta_j \\
&\leq P\left( \left( \frac{\tau^{(j+1)}(m_j) + n - \tau(m_j)}{n} \right)^{1/s} \rW_{\lambda_{j+1}} \leq x+\epsilon \right) \\
&\qquad + P\left( |\rW_{\lambda_{j+1}}| \geq \eta_j^{-1/s} \epsilon \right) +  \frac{4}{j} + \delta_j \\
&\leq P\left( \rW_{\lambda_{j+1}} \leq \frac{x+\epsilon}{\left(1-\eta_j \right)^{1/s}} \right) + P\left( |\rW_{\lambda_{j+1}}| \geq \eta_j^{-1/s} \epsilon \right) +  \frac{4}{j} + \delta_j, 
\end{align*}
where the second inequality follows from two applications of \eqref{mj2} and the last inequality follows from \eqref{mj1} and the fact that $n>\tau(m_j)$. 
Letting $\epsilon = \eta_j^{1/(2s)}$ and then taking $j$ large enough the right side can be made arbitrarily close to $P\left( \rW_{\lambda} \leq x \right)$ (note that we are using here that \eqref{Wljlim} and the definition of $\eta_j$ in \eqref{mj1} imply that $\eta_j \to 0$ as $j\to\infty$). Therefore, we can conclude that 
\[
 \limsup_{n\to\infty} \PP\left( \frac{X_n - \EE[X_n]}{n^{1/s}} \leq x \right)  \leq P\left( \rW_{\lambda} \leq x \right) , \quad \forall x \in \R. 
\]
A matching lower bound is proved similarly. 
\end{proof}


\begin{corollary}[{\bf Achievable generalized stable laws via interweaving}]\label{CoroClosed}
 For any $\lambda \in \Lconv$ there exists a cooling map $\tau$ and a constant $b>0$ such that $\frac{X_n - \EE[X_n]}{b n^{1/s}} \limdist \rW_{\lambda}$. 
\end{corollary}
\begin{proof}
First of all, note that if $\lambda(x) = c |x|^{-s-1} a(x)$ then $\frac{1}{b} \rW_{\lambda} \Law \rW_{\tilde\lambda}$, where $\tilde{\lambda}(x) = b^{-s}c |x|^{-s-1} a(bx)$. 
Thus, it is enough to prove the statement of the corollary only for $\lambda \in \Lconv$ with leading constant $c = \tc \vs s$.  
To this end, we fix a convex, non-decreasing function $a(x)$ on $(-\infty,0]$ that is vanishing at $-\infty$ and equals $1$ at $x=0$.
Then it is easy to see that there exists a sequence of functions $a_j(x)$ converging pointwise to $a(x)$ where
for each $j\geq 1$, the function $a_j(x)$ is a convex, non-decreasing function on $(-\infty,0]$ with compact support and whose graph consists of finitely many linear pieces. 
Moreover, the derivatives $a_j'(0)$ can be chosen so that 
\begin{itemize}
 \item if $a'(0)<\infty$ then $a_j'(0) = a'(0)$ for all $j\geq 1$, 
 \item while if $a'(0) = \infty$ then the functions $a_j'(0) = j$ for all $j\geq 1$. 
\end{itemize}
As shown in Example \ref{ex:finmix} above,
for each $j\geq 1$ there exists a cooling map $\tau^{(j)}$
such that the corresponding RWCRE $X^{(j)}$ has limiting distribution 
$\frac{X^{(j)}_n - \EE[X^{(j)}_n]}{n^{1/s}} \limdist \rW_{\lambda_j}$
where $\lambda_j(x) = \tc \vs s |x|^{-s-1} a_j(x)$, and 
the asymptotics of the cooling maps are given by
$\tau^{(j)}(n) \sim K_j n^{s/(s-1)}$ where $K_j = \left( \frac{s-1}{\vs s a_j'(0)} \right)^{s/(s-1)}$. 
Since the functions $a_j(x)$ converge pointwise to $a(x)$ and $|a_j(x)|\leq 1$,
it follows from the explicit form of the characteristic functions
in \eqref{Wlambdacritical} and the dominated convergence theorem that
$\rW_{\lambda_j} \limdist \rW_{\lambda}$ as $j\to \infty$. 
Also, the condition on the derivatives of $a_j'(0)$ implies that
$\lim_{j\to\infty} \frac{K_{j+1}}{K_j} = 1$.
Therefore, the sequence of cooling maps $\tau^{(j)}$ satisfy all
of the assumptions of Proposition \ref{closure},
and thus there exists a cooling map $\tau$ such that the corresponding
RWCRE has limiting distribution $\frac{X_n - \EE[X_n]}{n^{1/s}} \limdist \rW_{\lambda}$. 
\end{proof}

\begin{proof}[Proof of Theorem \ref{arbmix}]
By the proof of Corollary~\ref{CoroClosed} it is possible to
construct a polynomial cooling map $\tau$ for which the corresponding
RWCRE $X$ converges weakly after centering and scaling by $n^{1/s}$ to a random variable  $\rW_{\lambda}$ 
with $\lambda\in\Lconv$. By applying twice Proposition~\ref{intertwine},
we can then interweave this cooling map $\tau$ with two other
polynomial cooling maps such that~\eqref{triplemix} is satisfied.
\end{proof}

\section{ Examples of irregular cooling maps}\label{rEx}

Theorem \ref{arbmix} characterizes a large class of limiting distributions
that can be obtained for RWCRE $X$ associated to an $s$-canonical law $\mu$ on environments, and 
Theorems \ref{GeneralConditions-N} and \ref{GeneralConditions-S}
give sufficient conditions which imply the walk has a specified limiting distribution. 
These results, however, are not as complete as the results obtained in the case where the
law $\mu$ on environments is such that the RWRE is recurrent. 
In that case, a general limiting distribution result was obtained which identified
all possible limiting distributions and also identified necessary and sufficient
conditions on the cooling map $\tau$ for each of these distributions to arise as a
(subsequential) limiting distribution of the RWCRE \cite[Theorem 2]{ACdCdH20}.

One reason for the weaker results in the present paper is that while in \cite{ACdCdH20} the limiting distributions in all cases could be obtained by centering by the mean and scaling by the standard deviation of the walk, in the present paper the scaling one should use differs depending on the limiting distribution. For instance, in Theorem \ref{GeneralConditions-N} one obtains Gaussian limits after scaling by the  deviation whereas in Theorem \ref{GeneralConditions-S} one obtains stable or generalized tempered stable laws after scaling by $n^{1/s}$ (and in general $\sqrt{\Var(X_n)}$ doesn't grow like $n^{1/s}$). 
In light of the results for the interweaving of cooling maps in Section \ref{S8}, a natural idea to handle general cooling maps is to divide a cooling map into ``small'' cooling intervals which will give rise to a Gaussian component in the limit and ``large'' cooling intervals which will give rise to a stable or generalized tempered stable component in the limit. 
However, it is not a priori clear how to properly characterize the ``small'' and ``large'' cooling intervals to make this approach work. 

The following example gives a cooling map where it is immediately clear how to divide the ``small'' and ``large'' cooling intervals. This example is quite simple to analyze (given the earlier results in the paper) and shows that even in the cases where a limiting distribution is a pure Gaussian or pure stable one might still have to use this dividing approach to obtain the correct limiting distribution rather than simply applying a general result like Theorem \ref{GeneralConditions-N} or \ref{GeneralConditions-S}. 

\begin{example}[{\bf Parametric coling map for Gaussian \& Stable mixtures}]
  \label{mixedEx}
\text{}\\ 
Let $X$ be a RWCRE associated to an $s$-canonical law
$\mu$, as in Def.~\ref{goodenvironment} and cooling
increment sequence defined as follows.
Fix a parameter $r>1$ and consider the following sequence: 
\begin{equation}\label{rExample}
  T_{2^j} = \fl{r^j}, \text{ for } j\in \N,
  \quad \text{and} \quad
  T_k = 1, \text{ otherwise (that is, if $\log_2(k) \notin \N$)}.
\end{equation}
In this case the cooling intervals with $T_k = 1$ are considered ``small'' and all others are considered ``large''. With this in mind we can decompose $X_n - \EE[X_n]$ as 
\[
\sum_{k=1}^{\ell_n} \left(Z^{(k)}_{1} - \EE[Z_1]\right) \ind{T_k = 1} + \sum_{k=1}^{\ell_n}\left(Z^{(k)}_{T_k} - \EE[Z_{T_k}]\right) \ind{T_k > 1} + \left(Z^{(\ell_n+1)}_{\bar{T}_n} - \EE[Z_{\bar{T}_n}] \right). 
\]
Letting $a_\alpha^2 = \Var(Z_1)$, it follows from the classical CLT that the first term scaled by 
$A_{r,n} = a_\alpha ( \sum_{k=1}^{\ell_n} \ind{T_k = 1} )^{1/2}$
converges in distribution to $\rN$, and Theorem \ref{GeneralConditions-S} (using condition {\bf(S1)}) implies that the last two terms scaled by $B_{r,n} = \left( \sum_{k=1}^{\ell_n}T_k \ind{T_k > 1} + \bar{T}_n \right)^{1/s}$ converges in distribution to $\rS_s$. 	
By 
computing the asymptotics of $A_{r,n}$ and $B_{r,n}$ and comparing their relative sizes as $n\to\infty$, 
one can then obtain the following limiting distributions.  
\begin{itemize}
\item {\bf (Normal)} if $1 < r < 2^{s/2}$, then  
  \begin{equation}\label{pureG}
    \frac{X_{n} - \EE[X_{n}] }{\varz \sqrt{n} }
    \Longrightarrow  \rN.  
 \end{equation}
Note that in this case it can be shown that $\Var(X_n) \sim a_\alpha^2 n$
if and only if $r < 2^{1/(3-s)}$. Thus, it is evident that for
$r \in \left[ 2^{1/(3-s)}, 2^{s/2} \right)$ one cannot derive the Gaussian limiting
distribution using the approach of Theorem \ref{GeneralConditions-N}.
Indeed, it can be checked that this is because the Lindeberg
condition \eqref{Lindeberg} fails when $r \geq 2^{1/(3-s)}$.
On the other hand, when $r < \sqrt{2}$ one can derive \eqref{pureG}
by directly applying Theorem \ref{GeneralConditions-N}
since \eqref{Gauss_iff} holds in this case, whereas if
$r \in [\sqrt{2}, 2^{1/(3-s)})$ then even though \eqref{Gauss_iff}
doesn't hold the same proof idea works since the
Lindeberg condition \eqref{Lindeberg} can be verified when $r < 2^{1/(3-s)}$. 
\item {\bf (Mixture)} if $r = 2^{s/2}$, then the sequence $\frac{X_{n} - \EE[X_{n}] }{ \sqrt{n} }$ is tight and admits multiple limit points 
of the form $a \rN + b \rS_s$ where $\rN$ and $\rS_s$ are independent and $a,b>0$.  
\item {\bf (Stable)} 
If $r>2^{s/2}$ then the scaling limits are always stable, but the scaling is different when $r\leq 2$. If $2^{s/2} < r \leq 2$ then there exists a sequence of numbers $d_n$ which are bounded away from 0 and $\infty$ (and which do not converge to a constant as $n\to\infty$) such that 
\[
 \frac{X_n - \EE[X_n]}{d_n n^{\beta}} \limdist \rS_s, \quad\text{where the scaling exponent } \beta = \frac{\log_2(r)}{s} \in \left(\frac{1}{2}, \frac{1}{s} \right]. 
\]
On the other hand, if $r > 2$ then $\frac{X_n - \EE[X_n]}{n^{1/s}} \limdist \rS_s$. 
The above stable limiting distributions can be obtained by directly applying Theorem \ref{GeneralConditions-S} only when $r>2^s$.  Indeed, it can be shown that condition {\bf{(S1)}} holds if and only if $r > 2^s$. Condition {\bf{(S2)}} holds only when $r<2$ but in this case the limit in \eqref{cag} is $g(x) \equiv 1$ which violates the assumption that $g(0) = 0$ in Theorem \ref{GeneralConditions-S}.  
\end{itemize}
\end{example}

Another natural question regarding the characterization
of limiting distributions for general cooling maps is whether
or not Theorem~\ref{arbmix} identifies all possible limiting
distributions for RWCRE associated to $s$-canonical laws $\mu$ on environments.
In particular, can one obtain generalized tempered stable laws $\rW_\lambda$ with $\lambda \notin \Lconv$? 
Theorem \ref{GeneralConditions-S} provides a strategy
for how one might do this: find function $g\in \mathcal{G}$
for which $\lambda_g \notin \Lconv$, and then construct a cooling map
$\tau$ which satisfies \eqref{cag} for this choice of $g$
(and such that $\tau$ also satisfies condition {\bf (S2)}),
but we have been able to find such examples only along subsequences. 
That is, if we allow for subsequential weak limits, then
an easy modification of the proof of Theorem \ref{GeneralConditions-S}
shows that if there is a subsequence $n_j$ and
a function $g\in \mathcal{G}$ such that the limits in \eqref{cag}
and condition {\bf (S2)} hold along the subsequence $n_j$,
then $\frac{X_{\tau(n_j)}-\EE[X_{\tau(n_j)}]}{\tau(n_j)^{1/s}} \underset{j\to\infty}{\limdist} \rW_{\lambda_g}$.
The following gives an explicit example of how this can
be applied to get a subsequential limiting distribution
which is not of the type included in Theorem~\ref{arbmix}. 
\begin{example}[{\bf ``Exotic" cooling map}]\label{ex:exotic}
 Let $n_j = 2^{2^j}$, and let 
\[
 T_k = \left\lceil k n_j^{\frac{2-s}{s-1}} \right\rceil, \quad \text{for } n_{j-1} < k \leq n_j. 
\]
The interested reader can check that for this example one has $\tau(n_j) \sim \frac{1}{2} n_j^{s/(s-1)}$ and 
\[
 \lim_{j\to\infty} \sum_{k=1}^{n_j} \frac{T_k}{\tau(n_j)} \ind{T_k \leq x \tau(n_j)^{1/s}} = \left( \frac{x}{2^{1/s}} \right)^2 \wedge 1 =: g(x). 
\]
Therefore, it follows that $\frac{X_{\tau(n_j)}-\EE[X_{\tau(n_j)}]}{\tau(n_j)^{1/s}} \underset{j\to\infty}{\limdist} \rW_{\lambda_g}$ with 
\[
 \lambda_g(x) = \tc \vs s |x|^{-s-1}\left( 1 + \frac{(s-1)2^{(s-1)/s}}{\vs s} x - \frac{2-s}{s 2^{2/s} \vs^2} x^2 \right)_+. 
\]
Note that $\lambda_g \notin \Lconv$ in this case since $|x|^{s+1} \lambda_g(x)$ is not convex. 
\end{example}

\appendix

\section{Sums of heavy tailed random variables}\label{A}
This section contains some needed results regarding moment bounds and tail decay for sums of i.i.d.\ heavy tailed random variables. 
The results below seem to be part of the folklore known to experts
in heavy-tailed random variables though we could not find a
convenient reference, but we have included the proofs here both for completeness. 

Throughout this section we will assume that $\{\xi_i\}_{i\geq 1}$ are i.i.d.\ random variables and that $S_n = \sum_{i=1}^n \xi_i$. 

\begin{lemma}\label{Sntail}
Assume that $E[\xi_1] = 0$ and that $P(|\xi_1| > x) = \bigo( x^{-s})$ for some $s>1$.  Then, 
there exists a constant $C>0$ such that $P( |S_n| > t n^{1/s}) \leq C t^{-s}$ for all $t>0$ and $n$ large enough. 
\end{lemma}
\begin{proof}
It is enough to prove a bound $P(|S_n| > t n^{1/s}) \leq C_1 t^{-s}$ for some $C_1>0$ and all $t\geq t_1>0$ since we can then choose $C$ large enough so that $C t^{-s} \geq 1$ for $t \in (0,t_1)$. 
Now, first note that 
\begin{align*}
 P(|S_n| > t n^{1/s})
&\leq n P\left( |\xi_1| > \frac{t}{2} n^{1/s} \right) + P\left( \left| \sum_{k=1}^n \xi_k \ind{|\xi_k| \leq \frac{t}{2} n^{1/s}} \right| > t n^{1/s} \right), 
\end{align*}
and for some $C'>0$ and $n$ sufficiently large the first term can be bounded by $C' t^{-s}$ for all $t\geq 1$. For the second term, first of all note that
\begin{equation}\label{e:56}
  \left| E[\xi_k \ind{|\xi_k| \leq x} ] \right|
= \left| E[\xi_k \ind{|\xi_k| > x} ] \right|
\leq  E[|\xi_k| \ind{|\xi_k| > x} ] 
= \bigo( x^{1-s}) , 
\end{equation}
where the first equality follows from the assumption that $E[\xi_1] = 0$ and the last equality follows from the assumed tail asymptotics of $|\xi_1|$. 
Therefore, there exists a constant $b>0$ such that $\left| E\left[  \sum_{k=1}^n \xi_k \ind{|\xi_k| \leq \frac{t}{2} n^{1/s}} \right] \right| \leq b t^{1-s} n^{\frac{1}{s}}$ for $n$ sufficiently large and $t\geq 1$. 
Therefore, 
\begin{align*}
& P\left( \left| \sum_{k=1}^n \xi_k \ind{|\xi_k| \leq \frac{t}{2} n^{1/s}} \right| > t n^{1/s} \right)\\
&\leq P\left( \left| \sum_{k=1}^n \left\{ \xi_k \ind{|\xi_k| \leq \frac{t}{2} n^{1/s}} - E[\xi_k \ind{|\xi_k| \leq \frac{t}{2} n^{1/s}}] \right\} \right| > \left(1 - b t^{-s} \right) t n^{1/s} \right)\\
&\leq P\left( \left| \sum_{k=1}^n \left\{ \xi_k \ind{|\xi_k| \leq \frac{t}{2} n^{1/s}} - E[\xi_k \ind{|\xi_k| \leq \frac{t}{2} n^{1/s}}] \right\} \right| > \frac{t}{2} n^{1/s} \right), 
\qquad \text{for } t \geq t_0 = \left( 2b \right)^{1/s}. 
\end{align*}
Since the tail decay of $\xi_1$ implies that $\Var( \xi_1 \ind{|\xi_1| \leq x } ) \leq E\left[ \xi_1^2 \ind{|\xi_1| \leq x} \right] = \bigo(x^{2-s})$ as $x\to \infty$, 
then applying Chebychev's inequality to the above bound implies that there exists a constant $C''>0$ such that for $t\geq t_0$ and $n$ large enough
\begin{equation}\label{e:57}
 P\left( \left| \sum_{k=1}^n \xi_k \ind{|\xi_k| \leq \frac{t}{2} n^{1/s}} \right| > t n^{1/s} \right)
\leq \frac{4 n \Var( \xi_1 \ind{|\xi_1| \leq \frac{t}{2} n^{1/s}} ) }{ t^2 n^{2/s} } \leq C'' t^{-s}.
\end{equation}
Letting $t_1 = \max\{1,t_0\}$ and $C_1 = \max\{C', C'' \}$ we have that  $P(|S_n| > t n^{1/s}) \leq C_1 t^{-s}$ for all $t\geq t_1>0$ .
\end{proof}

Since Lemma \ref{Sntail} gives the same tail decay bound
for $S_n/n^{1/s}$ for $n$ sufficiently large, we immediately obtain the following corollary.

\begin{corollary}\label{LpSn}
 Assume that $E[\xi_1] = 0$ and that $P(|\xi_1| > x) = \bigo( x^{-s})$ for some $s>1$. 
 If $p \in (0,s)$, then $E[|S_n|^p] = \bigo( n^{p/s})$. 
\end{corollary}

Our final result in this section gives left tail asymptotics for $S_n$ when the random variables $\xi_1$ are bounded to the left and heavy tailed to the right. 
\begin{lemma}\label{Sntail-left}
 Assume that $\xi_1$ has mean zero, is bounded below (i.e. $P(\xi_1 \geq -L) = 1$ for some $L<\infty$), and has right tail decay $P(\xi_1 \geq x) = \bigo( x^{-s})$ for some $s \in (1,2)$.  
Then, there exists a constant $C>0$ such that 
\begin{equation}\label{Snlefttail}
 P\left( S_n < -t n^{1/s} \right) \leq e^{-C t^{\frac{s}{s-1}}}, \quad \text{for all } t>0. 
\end{equation}
\end{lemma}
\begin{proof}
 We begin by claiming that there is a constant $c'>0$ such that
\begin{equation}\label{leftmgfbound}
 E\left[ e^{-\lambda \xi_1} \right] \leq e^{c' \lambda^s}, \qquad \text{for all } \lambda > 0.
\end{equation}
For ease of notation let $\hat\xi = \xi_1 + L$ so that our assumptions on $\xi_1$ imply that $\hat\xi$ is a non-negative random variable and that $P(\hat{\xi} \geq x) \leq K x^{-s}$ for some $K > 0$ and $s \in (1,2)$. 
Then, 
\begin{align*}
 e^{-\lambda L} E\left[ e^{-\lambda \xi_1} \right] = E\left[ e^{-\lambda \hat\xi} \right]
&= 1-\int_0^\infty \lambda e^{-\lambda x} P(\hat\xi \geq x) \, dx \\
&\leq 1 - \int_0^\infty \lambda \left( 1 - \min\{\lambda x, 1\} \right) P(\hat\xi \geq x) \, dx \\
&= 1 - \lambda L + \int_0^{\lambda^{-1}} \lambda^2 x P(\hat\xi\geq x) \, dx + \int_{\lambda^{-1}}^\infty \lambda P(\hat\xi\geq x) \, dx \\
&\leq 1 - \lambda L + K \lambda^2 \int_0^{\lambda^{-1}} x^{1-s} dx + K \lambda \int_{\lambda^{-1}}^\infty x^{-s} \, dx \\
&= 1 - \lambda L + \frac{K}{(2-s)(s-1)} \lambda^s \\
&\leq e^{-\lambda L +  \frac{K}{(2-s)(s-1)} \lambda^s }. 
\end{align*}
(Note that in the third line above we used that $E[\hat\xi] = E[\xi_1] + L = L$ and in the second to last line we used that $s \in (1,2)$.)
This proves~\eqref{leftmgfbound} with $c' = \frac{K}{(2-s)(s-1)}$. 

The proof of~\eqref{Snlefttail} from~\eqref{leftmgfbound} follows standard large deviation techniques. First of all, it follows from Chebychev's inequality and then~\eqref{leftmgfbound} that 
\begin{equation}\label{e:58}
 P( S_n \leq -t n^{1/s} ) \leq e^{-\lambda t n^{1/s} } E\left[ e^{-\lambda S_n} \right]
\leq e^{-\lambda t n^{1/s} + c' \lambda^s n}, \quad \text{for any } \lambda > 0. 
\end{equation}
Choosing $\lambda = \left( \frac{t}{c' s n^{1-\frac{1}{s}}} \right)^{1/(s-1)}$, this gives the bound 
\begin{equation}\label{e:59}
 P( S_n \leq -t n^{1/s} ) \leq \exp \left\{ - (c')^{\frac{-1}{s-1}}\left( \frac{1}{s^{\frac{1}{s-1}}} - \frac{1}{s^{\frac{s}{s-1}}} \right)  t^{\frac{s}{s-1}} \right\}. 
\end{equation}
Since $s>1$ implies that $ \frac{1}{s^{\frac{1}{s-1}}} - \frac{1}{s^{\frac{s}{s-1}}} > 0$, this finishes the proof of~\eqref{Snlefttail}. 
\end{proof}


\section{RWRE: regeneration times for \texorpdfstring{$s>0$}{s>0}}\label{Regeneration}

We recall and collect some useful facts about \emph{regeneration times} associated to RWRE. 
For more details we refer the reader to~\cite{SZ99} and the references specified below. 

The sequence of regeneration times $(R_k)_{k\in \N}$ is defined as follows. 
\begin{equation}\label{e:60}
 R_1 := \inf\{n>0: \, \max_{\ell < n} Z_\ell < Z_n \leq \min_{m\geq n} Z_m \}, 
\end{equation}
and
\begin{equation}\label{Rk}
 R_k := \inf\{n>R_{k-1}: \, \max_{\ell < n} Z_\ell < Z_n \leq \min_{m\geq n} Z_m \}, \quad k\geq 2.  
\end{equation}
The important facts we will use about regeneration times is that they give an independence structure under the annealed measure $\PP$. 
\begin{itemize}
 \item The sequence of joint random variables 
 \begin{equation}\label{regseq}
  (Z_{R_1},R_1), (Z_{R_2}-Z_{R_1},R_2-R_1), (Z_{R_3}-Z_{R_2},R_3-R_2),\ldots
 \end{equation}
 is independent under the measure $\PP$, and all but the first term are identically distributed. 
 \item The joint sequence $\{(Z_{R_k}-Z_{R_{k-1}},R_k-R_{k-1})\}_{k\geq 2}$
   has the same distribution as that of \eqref{regseq}
   under the measure $\bPP(\cdot) = \PP(\, \cdot \mid Z_n \geq 0, \, \forall n\geq 0)$. 
\end{itemize}

As a consequence, the following identities in mean are valid for any $n\in\N$:
\begin{align}
 \EE[R_n] &= \EE[R_1] + (n-1) \bEE[R_1] = n \bEE[R_1] + \bigo(1), \label{tme} \\
\text{and} \quad  \EE[Z_{R_n}] &= 
 \EE[Z_{R_1}] + (n-1)\bEE[Z_{R_1}]  = n \bEE[Z_{R_1}] + \bigo(1). \label{ime}
\end{align}
Furthermore, it is worth noticing that the limiting speed $\vs$ of RWRE defined in \eqref{speedann} can be expressed in terms of regenerations as 
\begin{equation}\label{vR1} \vs = \frac{\EE[Z_{R_2}-Z_{R_1}]}{\EE[R_2-R_1]} = \frac{\bEE[Z_{R_1}]}{\bEE[R_1]}.\end{equation}

Kesten, Kozlov, and Spitzer~\cite{KKS75} studied transient one-dimensional
RWRE via a related Markov chain $\{V_i\}_{i\geq 0}$ which can be interpreted
as a branching process with immigration where each generation has an offspring
distribution which is a random Geometric distribution. However,
while they did not state their results this way, their analysis of the Markov
chain $V$ gives information on the regeneration structure of the RWRE. 
\begin{lemma}[{\bf Characterization in terms of branching processes}~\cite{KZ08}]\label{lem:rt-bp}
Let the Markov chain $V$ start at $V_0 = 0$ and let $\nu = \inf\{i>0: V_i = 0\}$.
Then the joint distribution of $(Z_{R_2}-Z_{R_1}, R_2-R_1)$ is the same as
the joint distribution of $(\nu, \nu + 2 \sum_{i=0}^{\nu-1} V_i )$. 
\end{lemma}
\begin{proof}
This was proved in~\cite[Lemma 12]{KZ08} for transient,
one-dimensional excited random walks. However, the proof carries over without
any changes to RWRE in i.i.d.\ environments. 
\end{proof}

\begin{corollary}[{\bf Tail control on regenerations and increments}]\label{cor:regtails}
There exist constants $C_1,C_2>0$ such that 
\begin{equation}\label{Xtau-tail}
  \PP(Z_{R_2}-Z_{R_1} > n) \leq C_1 e^{-C_2 n}. 
\end{equation}
Moreover, if $s \in (0,2]$ then there exists a constant $C_3 > 0$ such that 
\begin{equation}\label{tau-tail}
  \PP(R_2-R_1 > n) \sim C_3 n^{-s}, \quad \text{as } n\to \infty. 
\end{equation}
\end{corollary}
\begin{proof}
It was shown in~\cite[Lemmas 2 and 6]{KKS75} that 
\begin{equation}\label{bptail}
 P(\nu > n) \leq C_1 e^{-C_2 n}
 \quad \text{and}\quad
 P\left( \sum_{i=0}^{\nu-1} V_i > n \right) \sim K n^{-s}, \quad \text{as } n\to \infty, 
\end{equation}
for some constants $C_1,C_2,K>0$. 
Then,~\eqref{Xtau-tail} follows from Lemma~\ref{lem:rt-bp} and the above tail decay for $\nu$. 
Regarding~\eqref{tau-tail}, it follows from Lemma~\ref{lem:rt-bp} that 
\begin{equation}\label{e:61}
 P\left( \sum_{i=0}^{\nu-1} V_i > \frac{n}{2} \right)
 \leq \PP(R_2-R_1 > n)  
 \leq P\left( \sum_{i=0}^{\nu-1} V_i > \frac{n-\sqrt{n}}{2} + \sqrt{n} \right) + P( \nu > \sqrt{n} ). 
\end{equation}
Letting $C_3 = K 2^s$, \eqref{bptail} implies that both the lower bound
and upper bound above are asymptotic to $C_3 n^{-s}$ as $n\to \infty$.
\end{proof}

The following lemma gives a control on the 1st ``special' regeneration.
\begin{lemma}[{\bf First regeneration: p-moment and tail of displacement}~\cite{GP19,S04}]\label{lem:taumoments}
If $s>0$, then $\EE[R_1^p] < \infty$ if and only if $p \in (0,s)$. 
Furthermore, there exists a constant $c>0$ such that

\begin{equation}\label{ZR1tail}
  \EE\left[ e^{cZ_{R_1}}\right] <\infty.
\end{equation}

\end{lemma}

We refer the reader to~\cite[Prop. 3.5]{GP19} for the boundedness
of the $p$-moments of the first regeneration time stated above.
A proof of~\eqref{ZR1tail} can be found in~\cite{S04}, see Lemma
2.5 and Eq. (97) therein. For the latter, we stress that even
though this reference deals with high-dimensional setups assuming
directional transience along a given direction, see Eq. (80)
in~\cite{S04}, these statements remain still valid in dimension
one under our assumption~\eqref{trans}.

\section{Proof of Lemma \ref{AntContri}}\label{Acp}

We start by proving \eqref{e:24} and then move to \eqref{main} and \eqref{mainE}.

\begin{proof}[Proof of \eqref{e:24}]
First of all, from the definition of $A_{n,t}$, in~\eqref{Antdef},
we see that the condition in~{\bf(S2)} implies for any fixed $t,\delta>0$ that 
\begin{equation}\label{e:25}
 \ind{ \frac{T_k}{\tau(n)^{1/s}} \geq \frac{t}{\vs} + \delta} \leq \ind{k \in A_{n,t}} \leq \ind{\frac{T_k}{\tau(n)^{1/s}} \geq \frac{t}{\vs}}, \qquad \text{for all $n$ large.}
\end{equation}
To obtain~\eqref{e:24} it is enough to prove the following statement:
\[
\lim_{n\to\infty} \frac{1}{\tau(n)^{1-\frac{1}{s}}}
\sum_{k=1}^{n} f\left( \frac{T_k}{\tau(n)^{1/s}} \right) \ind{T_k/\tau(n)^{1/s} \geq a }
= \int_a^\infty \frac{f(x)}{x} \, g(dx) + L (1 - g(\infty)), \quad \forall a>0.
\]
A technical difficulty arises with the integral on the right
being over an unbounded interval, so will will prove the
following two statements which together imply the above limit:
\begin{equation}\label{dintab}
\lim_{n\to\infty} \frac{1}{\tau(n)^{1-\frac{1}{s}}} \sum_{k=1}^{n} f\left( \frac{T_k}{\tau(n)^{1/s}} \right) \ind{T_k/\tau(n)^{1/s} \in [a,b) }
= \int_a^b \frac{f(x)}{x}\, g(dx) ,\quad \forall 0<a<b<\infty.
\end{equation}
\begin{equation}\label{from_infty}
\lim_{b \to \infty}\limsup_{n\to\infty} \left| \frac{1}{\tau(n)^{1-\frac{1}{s}}} \sum_{k=1}^{n} f\left( \frac{T_k}{\tau(n)^{1/s}} \right) \ind{T_k \geq b \tau(n)^{1/s}} - L (1 - g(\infty)) \right| = 0
\end{equation}

To prove~\eqref{dintab}, 
we fix $0<a<b$ and a large integer $N$ and partition the interval $[a,b]$ into $N$ equally spaced sub-intervals $[x_{j-1,N},x_{j,N}]$ where $x_{j,N} = a + \frac{j(b-a)}{N}$ for $j=0,1,2,\ldots N$. For $j=1,2,\ldots, N$ let $x_{j,N}^*$ be a point in the interval $[x_{j-1,N},x_{j,N}]$ where the function $f(x)/x$ achieves its maximum. 
With this notation we can get the following upper bound on the sum 
\begin{align*}
& \limsup_{n\to\infty} \frac{1}{\tau(n)^{1-\frac{1}{s}}} \sum_{k=1}^{n} f\left( \frac{T_k}{\tau(n)^{1/s}} \right) \ind{T_k/\tau(n_j)^{1/s} \in [a,b) } \\
&\quad =  \limsup_{n\to\infty}\sum_{j=1}^N \sum_{k=1}^{n} \frac{ f\left( \frac{T_k}{\tau(n)^{1/s}} \right) }{\frac{T_k}{\tau(n)^{1/s}}} \frac{T_k}{\tau(n)} \ind{T_k/\tau(n)^{1/s} \in [x_{j-1,N},x_{j,N}) } \\
&\quad \leq  \sum_{j=1}^N \frac{ f\left( x_{j,N}^* \right) }{x_{j,N}^*} \left( \lim_{n\to\infty} \sum_{k=1}^{n} \frac{T_k}{\tau(n)} \ind{T_k/\tau(n)^{1/s} \in [x_{j-1,N},x_{j,N}) } \right) \\
&\quad = \sum_{j=1}^N \frac{ f\left( x_{j,N}^* \right) }{x_{j,N}^*} \left( g(x_{j,N}) - g(x_{j-1,N}) \right), 
\end{align*}
where in the last step we used again assumption \eqref{cag}. Finally, by taking $N\to\infty$, the upper bound above becomes arbitrarily close to the Riemann-Stieltjes integral $\int_a^b \frac{f(x)}{x} \, g(dx)$. 
The proof of the matching lower bound is obtained similarly.

To prove~\eqref{from_infty}, 
fix $\epsilon > 0$ and choose $b$ large enough so that
$\sup_{x\geq b} |\frac{f(x)}{x} - L| \leq \epsilon$. Then, 
\begin{align*}
 &\left| \frac{1}{\tau(n)^{1-\frac{1}{s}}} \sum_{k=1}^{n} f\left( \frac{T_k}{\tau(n)^{1/s}} \right) \ind{T_k \geq b \tau(n)^{1/s}} - L (1 - g(\infty)) \right| \\
&\leq \left| \sum_{k=1}^{n} \bigg(\frac{f\left( \frac{T_k}{\tau(n)^{1/s}} \right)}{\frac{T_k}{\tau(n)^{1/s}}} -  L \bigg) \frac{T_k}{\tau(n)} \ind{T_k \geq b \tau(n)^{1/s}}   \right| \\
&\qquad + L \left| \sum_{k=1}^{n} \frac{T_k}{\tau(n)} \ind{T_k > b \tau(n)^{1/s}} - (1-g(b)) \right| + L(g(\infty)-g(b)) \\
&\leq (L+\epsilon) \left| \sum_{k=1}^{n} \frac{T_k}{\tau(n)} \ind{T_k > b \tau(n)^{1/s}} - (1-g(b)) \right| + L(g(\infty)-g(b))
\end{align*}
It follows from \eqref{cag} that the first term in the last line above vanishes as $n\to \infty$. Then taking $b\to\infty$ finishes the proof of \eqref{from_infty}. 
\end{proof}

\begin{proof}[Proof of \eqref{main}]

Let us first fix a sequence $a(n)\to\infty$ with the properties that
\begin{equation}\label{a(n)}
 \max_{k\leq n} T_k^{1/s} \leq a(n) \quad \text{and} \quad a(n) = o(\tau(n)^{1/s} ),
\end{equation}
where we note that such a sequence is guaranteed by~{\bf(S2)}. 
Further, since $g(x)$ in~\eqref{cag} is 
less than 1 for some $x$ it follows from the first condition in \eqref{a(n)} that
for some $\delta>0$ and for $n$ large enough
\begin{equation}\label{anlb}
 a(n) \geq \delta \tau(n)^{1/s^2}.
\end{equation}

Next, partition the integers from $1$ to $n$ as the disjoint union 
$B_{n,t}^1 \cup B_{n,t}^2 \cup A_{n,t}$, with
$A_{n,t}$ as in \eqref{Antdef},
\begin{align*}
 &B_{n,t}^1 := \{k\leq n: \, T_k \vs - t \tau(n)^{1/s} \leq -a(n) \}, \\
\text{and} \quad &B_{n,t}^2 := \{k\leq n: \, -a(n) < T_k \vs - t\tau(n)^{1/s} < (\log T_k)T_k^{1/s} \}. 
\end{align*}

With this decomposition, the claim \eqref{main} it is equivalent to
\begin{equation}\label{Intj}
 \lim_{n\to\infty} \sum_{k \in B_{n,t}^j} \PP(\xi_{k,n} \leq -t) = 0, \quad\text{for  } j\in\{1,2\},
\end{equation}
which we next show first for $j=1$ and then for $j=2$.

\noindent\textbf{Sum over $B_{n,t}^1$:} First of all, 
note that if $k \in B_{n,t}^1$ then
\begin{align*}
 \PP(\xi_{k,n} \leq -t)
 &= \PP( Z_{T_k}  \leq  \EE[Z_{T_k}] - t \tau(n)^{1/s} )\\
 &\leq \PP( Z_{T_k} \leq -a(n) + \EE[Z_{T_k}] - T_k \vs ).
\end{align*}
The asymptotic behavior of $\EE[Z_{T_k}]$ from~\eqref{meanzn}
and our choice of $a(n)$ in\eqref{a(n)} imply that for $n$
sufficiently large we have $\EE[Z_{T_k}] - T_k \vs \leq a(n)/2$ for $k\leq n$.
Thus, for $n$ sufficiently large we may bound 
\begin{equation}\label{backtrack} 
\max_{k\in B_{n,t}^1} \PP(\xi_{k,n} \leq -t)
 \leq \PP\left( \inf_{n\geq 0} Z_n \leq -\frac{a(n)}{2} \right) \leq C e^{-c a(n)},
\end{equation}
where the last inequality follows from~\cite[Lemma 3.3]{GS02}.
Since there are at most $n$ terms in $B_{n,t}^1$,
each bounded by \eqref{backtrack} this proves \eqref{Intj} for $j=1$.

\noindent\textbf{Sum over $B_{n,t}^2$:}
Eq. \eqref{varzn} implies that 
$\PP(\xi_{k,n} \leq -t) \leq \frac{1}{t^2 \tau(n)^{2/s}} \Var(Z_{T_k}) \leq \frac{C T_k^{3-s}}{t^2 \tau(n)^{2/s}}$,
from which we have that 
\begin{equation}\label{e:27}
 \sum_{k \in B_{n,t}^2} \PP(\xi_{k,n} \leq -t)
 \leq \frac{C}{t^2 \tau(n)^{2/s}} \sum_{k \in B_{n,t}^2} T_k^{3-s} 
 =  \frac{C}{t^2 \tau(n)^{1-\frac{1}{s}}} \sum_{k \in B_{n,t}^2} \left( \frac{T_k}{\tau(n)^{1/s}} \right)^{3-s}.
\end{equation}
Since we have chosen $a(n) = o(\tau(n)^{1/s})$ and we are assuming {\bf(S2)},
then it follows for any fixed $\delta > 0$, that for all $n$ large enough,
$k \in B_{n,t}^2$ implies that $|\frac{T_k}{\tau(n)^{1/s}} - \frac{t}{\vs} | < \delta$. 
Therefore, we have that for $n$ large enough 
\begin{equation}\label{e:28}
 \sum_{k \in B_{n,t}^2} \PP(\xi_{k,n} \leq -t) \leq \frac{C}{t^2 \tau(n)^{1-\frac{1}{s}}} \sum_{k \leq n } \left( \frac{T_k}{\tau(n)^{1/s}} \right)^{3-s} \ind{|\frac{T_k}{\tau(n)^{1/s}} - \frac{t}{\vs} | < \delta }.
\end{equation}
It follows from~\eqref{dintab} that we can compute the limit of this upper bound so that 
\begin{equation}\label{e:29}
 \limsup_{n\to\infty}  \sum_{k \in B_{n,t}^2} \PP(\xi_{k,n} \leq -t)
\leq \frac{C}{t^2} \int_{\frac{t}{\vs}-\delta}^{\frac{t}{\vs}+\delta} x^{2-s} \, g(dx). 
\end{equation}
Since the right-hand side vanishes as $\delta\to 0$ (recall that $g$ is continuous), then it follows that 
\eqref{Intj} holds for $j=2$ as well.
\end{proof}

\begin{proof}[Proof of \eqref{mainE}]
As in the proof of \eqref{main} above, we proceed by showing 
\begin{equation}\label{cmcj} 
 \lim_{n\to\infty} \sum_{k \in B_{n,t}^j} \EE\left[  \xi_{k,n} \ind{\xi_{k,n}\leq -t} \right]  = 0, \quad 
 \text{for } j=1,2. 
\end{equation}

\noindent\textbf{Sum over $B_{n,t}^1$:} 
It follows from the Cauchy-Schwartz inequality, \eqref{varzn}, and \eqref{backtrack} that 
\[
 \left| \sum_{k \in B_{n,t}^1} \EE\left[  \xi_{k,n} \ind{\xi_{k,n}\leq -t} \right] \right|
\leq \sum_{k \in B_{n,t}^1} \EE\left[  \xi_{k,n}^2  \right]^{1/2} \PP(\xi_{k,n}\leq -t)^{1/2}
\leq \frac{C e^{-c a(n)}}{\tau(n)^{1/s}} \sum_{k \in B_{n,t}^1} T_k^{\frac{3-s}{2}}. 
\]
Since the definition of $B_{n,t}^1$ implies that $T_k \leq \frac{t}{v} \tau(n)^{1/s}$ for all $k \in B_{n,t}^1$, and since $|B_{n,t}^1| \leq n$, we have that
\[
 \left| \sum_{k \in B_{n,t}^1} \EE\left[  \xi_{k,n} \ind{\xi_{k,n}\leq -t} \right] \right|
\leq \frac{C e^{-c a(n)/2}}{\tau(n)^{1/s}} n \left( \frac{t}{\vs} \tau(n)^{1/s} \right)^{\frac{3-s}{2}}
= \frac{C\left( \frac{t}{\vs} \right)^{\frac{3-s}{2}} n e^{-c a(n)/2}}{\tau(n)^{\frac{s-1}{2s}}}. 
\]
It follows from \eqref{anlb} that this upper bound vanishes as $n\to\infty$, 
and this proves \eqref{cmcj} for $j=1$.

\noindent\textbf{Sum over $B_{n,t}^2$:} It follows from the Cauchy-Schwartz inequality,
Chebychev's inequality, and~\eqref{varzn} that 
\[
  \left| \EE\left[  \xi_{k,n} \ind{\xi_{k,n}\leq -t} \right] \right|
\leq \EE\left[  \xi_{k,n}^2  \right]^{1/2} \PP(|\xi_{k,n}|\geq t)^{1/2}
\leq \frac{\EE[\xi_{k,n}^2] }{t}
\leq \frac{C}{t} \frac{T_k^{3-s}}{ \tau(n)^{2/s}},
\]
by arguing as in~\eqref{e:27} and right after it, we see that 
the sum of this upper bound over $B_{n,t}^2$
vanishes as $n$ increases. This proves \eqref{cmcj} for $j=2$.  
\end{proof}

\bibliographystyle{plain}
\bibliography{BibRWCREstable}
 
\end{document}